\documentclass[10pt]{article}
\usepackage{latexsym,amssymb,amsmath}
\usepackage{appendix}
\begin{document}
\baselineskip=20pt
\newcommand{\for}{\mbox{for}}
\newcommand{\ptl}{\partial}
\newcommand{\mbb}{\mathbb}
\newcommand{\ol}{\overline}
\newcommand{\stl}{\stackrel}
\newcommand{\vt}{\vartheta}
\newcommand{\lmd}{\lambda}
\newcommand{\la}{\langle}
\newcommand{\ra}{\rangle}
\newcommand{\Dlt}{\Delta}
\newcommand{\al}{\alpha}
\newcommand{\be}{\beta}
\newcommand{\G}{\Gamma}
\newcommand{\gm}{\gamma}
\newcommand{\ves}{\varepsilon}
\newcommand{\psp}{\vspace{0.4cm}}
\newcommand{\pse}{\vspace{0.2cm}}
\newcommand{\td}{\tilde}
\newcommand{\rta}{\rightarrow}

\title{\bf Proof of Kac and Rudakov's Conjecture  on    Generalized
Verma Module  over
Lie \\ Superalgebra  $\mbox{E}(5,10)$\footnote{Research
supported by CSC} }

\author{ $\mbox{Yufeng Zhao}$  \\\\\\  LMAM, School of
Mathematical Sciences \\Peking University,
  Beijing, 100871, P. R. China.}

\date{ }
\maketitle

\begin{abstract}The exceptional infinite-dimensional linearly compact simple Lie superalgebra  $\mbox{E}(5,10)$, which Kac believes, is the algebra of symmetries of the $\mbox{SU}_{5}$ Grand Unified Model.
In this paper, we give a proof of Kac and Rudakov's conjecture about the classification of  all the degenerate generalized Verma module over  $\mbox{E}(5,10)$. Also, we
  work out all the nontrivial singular vectors degree by degree.
It is a potential that the representation theory of
$ \mbox{E}(5,10)$ will shed new light on various features of the
the $\mbox{SU}_{5}$ Grand unified model.
\end{abstract}

\section {Introduction}

\quad\quad A  linearly compact infinite-dimensional Lie algebra
is a topological Lie algebra whose underlying space is a topological
space isomorphic to the space of formal power series over complex field
 in finite number of variables with formal
topology. Cartan's list of linearly compact infinite-dimensional  simple Lie
algebras consists of four series: the Lie algebra of all complex vector fields
and its subalgebras of divergence 0 vector fields, symplectic
vector fields and contact vector fields.

 Kac proved the ``super'' version of this result. In other words, he classified linearly compact infinite-dimensional Lie superalgebras [K1]. There turn out to be 10 families and 5 exceptions, which are called $\mbox{E}(1,6), \mbox{E}(3,6), \mbox{E}(3,8), \mbox{E}(4,4)$ and $ \mbox{E}(5,10)$.
 Many of the families are straightforward ``super" generalizations of the 4 families of linearly compact infinite-dimensional  simple Lie
algebras. Some are stranger. Most important for us today are the 5 exceptions discovered by Irina Shchepochkina [Sh].

 The representation theory of $E(3, 6)$ and $E(3,8)$ was developed by Kac and Rudakov [KR1-KR3], and some further observations were made on
its connections to the Standard Model [K2].
It was  found  quite remarkable that the $\mbox{SU}_{5}$ Grand unified model of Georgi-Glashow combines the left multiplets of fundamental fermions in precisely
the negative part of the consistent gradation of $\mbox{E}(5, 10)$.
This is perhaps an indication of the possibility that an extension from
$\mbox{su}_5$ to algebra of internal symmetries may resolve the difficulties
with the proton decay. It is a potential that the representation theory of
$ \mbox{E}(5,10)$ will shed new light on various features of the
the $\mbox{SU}_{5}$ Grand unified model.

As to the representation theory of $E(5, 10)$, Kac and Rudakov formulate an conjecture [KR3], which can be stated as follows.
 The Lie superalgebra $L = E(5, 10)$ carries a unique consistent irreducible
$\mathbb{Z}$-gradation $L = \bigoplus\limits_{j\geq -2}L_j$, where $L_0$ is isomorphic to simple Lie algebra $ sl_5$.
Given  $L_0$- module  $V$, we extend it to a $L$ module by letting $L_{+}$ acts trivially, and define the induced module
$$M(V)=U(L)\otimes_{U(L_0)}V\cong U(L_{-})V. $$
If $V$ is finite-dimensional irreducible $L_0$
-module, the $L$ module $M(V)$ is called a {\it  generalized Verma module } associated to $V$, and it is called {\it degenerate } if it is not irreducible.

We denote by $\mbox{V} (\lambda_1,\lambda_2,\lambda_3, \lambda_4)$ the finite-dimensional irreducible $L_0$
-module with highest
weight $\sum\limits_{i=1}^{4}\lambda_i\omega_i$, where $\omega_1, \omega_2,\omega_3,\omega$ are the fundamental weights for $ sl_5$. Let
$$M=M(\lambda_1,\lambda_2,\lambda_3, \lambda_4) = M(V (\lambda_1,\lambda_2,\lambda_3, \lambda_4))$$
denote the corresponding generalized Verma module over $\mbox{E}(5, 10)$.
Denote by  $
\mathbb N $ the additive semigroup of nonnegative integers.\vspace{0.1cm}

{\it Kac and Rudakov's Conjecture}  \quad {\it The following is a complete list of degenerate Verma modules over
$\mbox{E}(5, 10)$:
$$\mbox{M}(m, n, 0, 0); \mbox{M}(m, 0, 0, n);  \mbox{M}(0, 0, m, n) \  ( m, n \in \mathbb{N}). $$}

In this paper, we give a proof of this conjecture and  work out all the nontrivial singular vectors for any generalized Verma module over $\mbox{E}(5, 10)$.
The first key point of our proof is   investigating that there exists
 a grading on the generalized Verma module, through which we find that  any singular vector  is controlled by its leading term via an exponential-like differential operator, where the leading term lies in certain tensor product module of $sl_5$; the second one is the observation of
  an irreducible tensor operator of rank $\omega_1+\omega_2$ for simple Lie algebra $sl_5$, which plays the center role in our calculation of all the singular vectors.

The paper is organized as follows: In section 2, we recall Kac's geometric  construction of Lie superalgebra $\mbox{E}(5,10)$ and the  KR conjecture. In Section 3, we provide  some   techniques concerning the irreducible tensor operators and tensor module decomposition  theory of simple Lie algebra. In Section 4, we
prove that
all the nontrivial singular vectors are of degree less than or equal to four. Also,
the leading term of any singular vector must lie in one of the tensor decomposition  of four  tensor product module of $sl_5$ (cf. Theorem 4.6). In Section 5,
 we work out all the nontrivial singular vectors degree by degree (cf. Theorem 5.3, Theorem 5.4, Theorem 5.5, Theorem  5.6).\vspace{0.1cm}

\section{Lie superalgebra $\mbox{E}(5,10)$ and KR conjecture }

 \quad \quad In  this section, we recall Kac's geometric  construction of Lie superalgebra $\mbox{E}(5,10)$ and  KR Conjecture which are stated in [KR3].

 For two
integers $m<n$, we denote $\overline{m,n}=\{m,m+1,\cdots,n\}$.
 Let
$$\mbox{W}_n=\{\sum\limits_{i=1}^{n}p_i(x)\partial_i \ | \ p_i(x) \in \mathbb{C}[[x_1,\cdots,x_n]],  \partial_i=\partial_{x_{i}}\}\eqno(2.1)$$
denote the Lie algebra of formal vector fields in $n$ indeterminates;
$$\mbox{S}_n=\{D=\sum\limits_{i=1}^{n}p_i\partial_i \ | \ \mbox{div}D=\sum\limits_{i=1}^{n}\partial_i(p_i)=0 \}\eqno(2.2)$$denote the
Lie subalgebra of divergenceless formal vector fields; $\Omega^k(n)$ denote the associative algebra of formal differential forms of degree k in $n$ indeterminates, $\Omega^k_{\mbox{cl}}(n)$ denote the
subspace of closed forms.

The exceptional infinite-dimensional linearly compact Lie superalgebra $\mbox{E}(5,10)=\mbox{E}(5,10)_{\underline{0}}
+\mbox{E}(5,10)_{\underline{1}}$ is constructed as follows:
$$\mbox{E}(5,10)_{\underline{0}}=S_5, \ \mbox{E}(5,10)_{\underline{1}}=\Omega_{\mbox{cl}}^{2}(5),\eqno(2.3)$$
where $\mbox{E}(5,10)_{\underline{0}}$ acts on $\mbox{E}(5,10)_{\underline{1}}$ via the Lie derivative,
$$[\omega_2,\omega_2']=\omega_2 \wedge \omega_2'\in \Omega_{\mbox{cl}}^{4}(5)=S_5 \eqno(2.4)$$ for $\omega_2,\omega_2' \in
\mbox{E}(5,10)_{\underline{1}}$.

We use for the odd elements of $\mbox{E}(5,10)$   the notation $d_{ij} = dx_i \wedge dx_j (i, j \in \overline{1,5})$; recall that we have the following commutation relation ($ f, g \in  C[[x_1, \cdots , x_5]]$ ):
$$[fd_{jk}; g d_{lm}] = \varepsilon_{ijklm}fg\partial_i
;\eqno(2.5)$$
where
 $$\varepsilon_{ijklm} =\left\{\begin{array}{ll}
\mbox{the \ sign \ of \  the  \ permutation} \  $(ijklm)$,  & \mbox{ if \  all  \ indices  \  $ijklm$  \  are \  distinct}  , \\
0,  & \mbox{otherwise}
. \end{array}\right.\eqno(2.6)$$\vspace{0.1cm}
And the Lie superalgebra $\mbox{L}=\mbox{E}(5,10)$ carries a unique consistent irreducible
$\mathbb{Z}$-gradation $\mbox{L}= \bigoplus\limits_{j\geq -2}\mbox{L}_j$. It is defined by:
$$deg x_i = 2 = -\partial_i,  deg d_{ij} =-1 \eqno(2.7)$$
One has: $\mbox{L}_0 \simeq \mbox{sl}_5$ and the $\mbox{L}_0$-modules occurring in the negative part are:
$$\mbox{L}_{-1}=\mbox{Span}_{\mathbb{C}}\{ d_{ij} \ | \ i, j \in \overline{1,5}\} \simeq \Lambda^2\mathbb{C}^5,$$
$$\mbox{L}_{-2}=\mbox{Span}_{\mathbb{C}}\{ \partial_i \ | \ i \in \overline{1,5} \}\simeq
\mathbb{C}^{5*}\eqno(2.8)$$
Recall also that $\mbox{L}_1$
consist of closed 2-forms with linear coefficients, that $\mbox{L}_1$ is an irreducible $\mbox{L}_0$-module and $\mbox{L}_j = [\mbox{L}_1[\cdots]] =\mbox{L}_1^{j}$ for $ j \geq 1$.
We take for the Borel subalgebra of $\mbox{L}_0 \simeq \mbox{sl}_5$the subalgebra of the vector fields
$$\mbox{Span}\{x_i\partial_j (1\leq i\leq j \leq 5),  \ x_i\partial_i-x_{i+1}\partial_{i+1} (i \in \overline{1,4}) \}.\eqno(2.9)$$
Given $\mbox{L}_0$ module  $V$, we extend it to a $\mbox{L}$ module by letting $\mbox{L}_{+}$ acts trivially, and define the induced module
$$M(V)=U(\mbox{L})\otimes_{U(\mbox{L}_0)}V\cong U(\mbox{L}_{-})V. \eqno(2.10)$$
 We denote by $\mbox{V} (\lambda_1,\lambda_2,\lambda_3, \lambda_4)$ the finite-dimensional irreducible $\mbox{L}_0$
-module with highest
weight $\sum\limits_{i=1}^{4}\lambda_i\omega_i$, where $\omega_i$ ($i \in \overline{1,4}$) are the fundamental weights for $sl_5$. Let
$$M=M(\lambda_1,\lambda_2,\lambda_3, \lambda_4) = M(V (\lambda_1,\lambda_2,\lambda_3, \lambda_4)) \eqno(2.11)$$
denote the corresponding generalized Verma module over $\mbox{E}(5, 10)$.
 \vspace{0.1cm}

{\it Definition 2.1} \quad {\it  If $\xi \in M$
satisfies :
$$(x_i\partial_{x_{i+1}}).\xi=0 ( i \in \overline{1,4}),\eqno(2.12)$$$$ x_5d_{45}.\xi=0,\eqno(2.13)$$
then we call $\xi$  a singular vector for generalized  verma module $M$ of $E(5,10)$.}\vspace{0.1cm}

The aim of the following sections is to determine all the nontrivial singular vectors for $\mbox{E}(5, 10)$-module $M$.\vspace{0.1cm}

\section {Preliminary }

 \quad \quad  In this section, we give  some preparatory  techniques about the irreducible tensor operators and the decomposition of tensor product module of simple Lie algebra.

 Following the notations of Humphreys [H], let $H$  be the Cartan subalgebra of simple Lie algebra $L $, and let $\Delta=\{
 \alpha_1,\cdots,\alpha_l\}$  be a base for the root system  $ \phi$ of $H^{*}$. The corresponding fundamental dominant weights $\{\omega_1,\cdots,\omega_l\}$ are defined from the root system via the form $ <\cdot,\cdot>$ given by:
 $$<\omega_i,\alpha_j>\equiv \frac{2(\omega_i,\alpha_j)}{(\alpha_j,\alpha_j)}=\delta_{ij},\eqno(3.1)$$
 where $(\cdot,\cdot)$ denotes the inner product induced on
 $H^{*}$ by the Killing form on $H$. Consider a basis $\{h_1,\cdots, h_l, x_{\alpha}, \alpha \in \phi\}$ of $L$ where $h_1,\cdots, h_l$ is a basis for $H$ and $x_{\alpha}$ is a nonzero element of the root space $L_{\alpha}$. The dual basis may therefore be written $\{h^1,\cdots, h^l, x^{\alpha},\alpha \in \phi\} $ where $x^{\alpha}$ is the unique element of $L_{-\alpha}$ which is dual to $x_{\alpha}$ under the Killing form of $L$. Write the universal Casimir element in the form:
 $$c_{L}=\sum\limits_{i=1}^{l}h_ih^{i}+\sum\limits_{\alpha\in \phi}x_{\alpha}x^{\alpha}.\eqno(3.2)$$
Let $V(\mu)$ be an irreducible highest weight module over $L $ and let $\pi_{\mu}$ be the representation afforded by $ V(\mu)$.
 Choose an ordeded basis $\{e_1,\cdots, e_d\}$ of $V(\mu)$, let $\pi_{\mu}(x)$ denote the matrix representing $x\in L$ on
 $V(\mu)$ with respect to this basis.
 \vspace{0.1cm}

{\it Definition  3.1}\quad {\it  We call a collection of linear operators $\{T_i: V\rightarrow W  \ | \ i \in \overline{1,d}\}$
an irreducible tensor operator of rank $\mu$ if these components  transform according to the rule: $$[x, T_i]=\pi_{W}(x)T_i-T_i\pi_{V}(x)=\sum\limits_{j=1}^{d}\pi_{\mu}(x)_{ji}T_j, \ x \in L, \eqno(3.3) $$
where
$V$,$W$ are (possibly infinite dimensional) $L$-modules and $\pi_{V} $ (resp.
$\pi_{W}$) is the representation afforded by $V$ (resp. $W$).}\vspace{0.1cm}

Then we can define  the following  intertwining operator between $L$-modules $V(\mu)\otimes V $ and $ W$:
 $$T:V(\mu)\otimes V \rightarrow W, \ T(e_i\otimes v)=T_i(v),  \ i \in \overline{1,d}, \ v \in V. \eqno(3.4)$$
 In other words, $T\in \mbox{Hom}_{L}(V(\mu)\otimes V,W)$ is an element of the set of all operators from $V(\mu)\otimes V  $ to $W$ commuting with the action of $L$.
 \vspace{0.1cm}

 {\it Remark  3.2}\quad {\it In Section 5, we find an irreducible tensor operator of rank $\omega_1+\omega_2$ for simple Lie algebra $sl_5$, which play the center role in our determining all the singular vectors.}\vspace{0.1cm}

In the following two Lemmas, we record some well-known facts concerning the decomposition of tensor modules:\vspace{0.1cm}

{\it Lemma 3.3}\quad {\it  (1) (cf. [H]) \quad  The $\alpha$- string through any weight $\nu$  of $V(\mu)$  is of length $<\nu,\alpha>$, for  $\alpha\in \phi$.

(2) (cf. [EG])\quad  Denote $\mu_1,\cdots, \mu_m$ the weights occurring in  $V(\mu)$ with multiplicities $n_1,\cdots,n_m$ respectively. For each $i \in \overline{1,m}$, let $V_i(\mu)$ denote the space of weight vectors of weight $\mu_i$. The decomposition  of the tensor product module $V(\mu)\otimes V(\lambda)$ is written:
$$V(\mu)\otimes V(\lambda)=\sum\limits_{i=1}^{m}m(\lambda+\mu_i:\mu \otimes\lambda)V(\lambda+\mu_i), \  \lambda+\mu_i \in \Lambda^{+},\eqno(3.5)$$where the multiplicities are given by $$m(\lambda+\mu_i:\mu \otimes\lambda)=dim V_{i,\lambda}(\mu),\eqno(3.6) $$$$ V_{i,\lambda}(\mu)=\{v \in
V_i(\mu) \ | \ e_{j}^{<\lambda+\delta,\alpha_j>}v=0, j \in \overline{1,l}.
 \}\eqno(3.7)$$

(3) (cf. [EG]) \quad  Assume $\{e_{i,j} \ | \ j \in \overline{1,m(\lambda+\mu_i:\mu \otimes\lambda)}
\}$ is a basis for the space $V_{i,\lambda}(\mu)$ and $v_{\lambda}$ is  the maximal weight vector of $V(\lambda)$.
 A full set of independent maximal weight states of weight $\lambda+\mu_i$ is given by the vectors:
$$\{P_i(e_{i,j}\otimes v_{\lambda}),  \ j \in \overline{1,m(\lambda+\mu_i:\mu \otimes\lambda)}\},$$ where $$
P_i=\prod\limits_{\mu_i < \sigma \leq \mu} \frac{\tilde{c_L}-\chi_{\sigma+\lambda}(\tilde{c_L})}
{\chi_{\mu_i+\lambda}(\tilde{c_L})-\chi_{\sigma+\lambda}(\tilde{c_L})}, $$$$ \chi_{\sigma+\lambda}(\tilde{c_L})=\frac{(\sigma+\lambda,
\sigma+\lambda+2\delta)-
(\mu,\mu+2\delta)-(\lambda,\lambda+2\delta)}{2},$$$$\tilde{c_L}=\sum\limits_{i=1}^{l}
\pi_{\mu}(h_i)\otimes \pi_{\lambda}(h^i)+
\sum\limits_{\alpha\in \phi}\pi_{\mu}(x_{\alpha})\otimes \pi_{\lambda}(x^{\alpha}).
\eqno(3.8)$$

(4) (cf. [MS]) \quad The tensor product module $V(\mu)\otimes V(\lambda)$ is
a cyclic module which is cyclically generated by the vector $v^{\mu}\otimes v_{\lambda}$, where $v^{\mu}$ is the lowest weight vector for $V(\mu)$ and $v_{\lambda}$ is the highest weight  vector for $V(\lambda)$.
}
\vspace{0.1cm}

For $\overrightarrow{a}=(a_1,\cdots,a_{n-1}) \in \mathbb{N}^{n-1}$ and $0 < k \in \mathbb{N}$, we denote$$\overrightarrow{a}^{*}=(a_{n-1},a_{n-2},\cdots,a_{1}),\eqno(3.9)$$
$$ I(\overrightarrow{a},k)=\{(a_{1}+c_{1}-c_{2},a_{2}+c_{2}-c_{3}
,\cdots,a_{n-1}+c_{n-1}-c_{n})
 \ | \ c_{i} \in \mathbb{N} \ $$$$\mbox{such \ that} \ \sum\limits_{i=1}^{n}c_{i}=k \ \mbox{and} \ c_{s+1} \leq a_{s} \mbox{for} \ s \in \overline{1,n-1} \}.\eqno(3.10)$$
Set
$$\omega_{\overrightarrow{a}}=\sum\limits_{i=1}^{n-1}a_{i}\omega_{i}, \  \ \mbox{for} \ \overrightarrow{a}\in \mathbb{N}^{n-1}, \eqno(3.11)$$

{\it Lemma 3.4}\quad {\it (Pieri's formula  cf. [FH]) (1)\quad For any $\overrightarrow{a}\in \mathbb{N}^{n-1}$, the tensor product of $sl_n$-module $V(\omega_{\overrightarrow{a}})$ with $V(k\omega_1)$ decomposes into a direct sum:
$$V(\omega_{\overrightarrow{a}})\otimes V(k\omega_1)=\bigoplus_{\overrightarrow{b}\in I(\overrightarrow{a},k) }V(\omega_{\overrightarrow{b}}).\eqno(3.12)$$

 (2)\quad For $sl_n$,   we have $V(\omega_{\overrightarrow{a}})^{*}=V(\omega_{\overrightarrow{a}^{*}})$ and
$$V(\omega_{\overrightarrow{a}})\otimes V(k\omega_{n-1})
=\bigoplus_{\overrightarrow{b}\in I(\overrightarrow{a}^{*},k) }V(\omega_{\overrightarrow{b}^{*}}).\eqno(3.13)$$}
\vspace{0.1cm}

In the rest of this section, we will concentrate on some special wedge and tensor modules for $sl_5$.
Take  $\{h_i=E_{i,i}-E_{i+1,i+1} (i \in \overline{1,4}), E_{ij} (1 \leq i \neq j \leq 5)\}$ as a basis for Lie algebra $sl_5$.
Then $\{\frac{ h_{i}^{*}}{10} (i \in \overline{1,4}), \frac{E_{ji}}{10} (1 \leq i \neq j \leq 5)\}$ is its dual basis via the Killing form,
where $$ h_{1}^*=\frac{4}{5}h_1+\frac{3}{5}h_2+\frac{2}{5}h_3+\frac{1}{5}h_4, \
h_{2}^*=\frac{3}{5}h_1+\frac{6}{5}h_2+\frac{4}{5}h_3+\frac{2}{5}h_4,$$$$
h_{3}^*=\frac{2}{5}h_1+\frac{4}{5}h_2+\frac{6}{5}h_3+\frac{3}{5}h_4, \
h_{4}^*=\frac{1}{5}h_1+\frac{2}{5}h_2+\frac{3}{5}h_3+\frac{4}{5}h_4.\eqno(3.14)$$
\vspace{0.1cm}
 And
 the Casimir operator $c$ of the universal enveloping algebra of $sl_5$ is
$$c=\frac{1}{10}(\sum\limits_{i=1}^{4}h_{i}h_{i}^{*}+\sum\limits_{i\neq j\in \overline{1,5}}E_{i,j}.E_{j,i}).\eqno(3.15)$$
Relative to the ordered basis $ \omega_1, \omega_2, \omega_3,  \omega_4$, the coordinates of the simple roots $\alpha_i (i \in \overline{1,4})$ are:
 $$\alpha_1=(2, -1,0,0), \ \alpha_2=(-1, 2, -1, 0), \ \alpha_3=(0,-1,2,-1), \ \alpha_4=(0,0,-1,2).\eqno(3.16)$$
  And the  killing form  for the simple root $\alpha_i (i \in \overline{1,4})$ are:
$$(\alpha_i,\alpha_j)=\left\{\begin{array}{lll}
0,  & |i-j| > 1,\\
\frac{1}{5}, &
i
=j,\\
-\frac{1}{10}, & |i-j|=1. \end{array}\right.\eqno(3.17)$$
{\it Lemma 3.5 }\quad {\it Assume $L=sl_5$, $C_{L}=c$ and $\sigma=\mu-\sum\limits_{i=1}^{4}k_i\alpha_i$ in Lemma 3.4.
Then $\chi_{\sigma+\lambda}(\tilde{c})$ in (3.6) is explicitly given by:
\begin{eqnarray}&&
\chi_{\sigma+\lambda}(\tilde{c})=\frac{\lambda_1(4\mu_1+3\mu_2+2\mu_3+\mu_4)}{50}
+\frac{\lambda_2(3\mu_1+6\mu_2+4\mu_3+2\mu_4)}{50}
+\frac{\lambda_3(2\mu_1+4\mu_2+6\mu_3+3\mu_4)}{50}\nonumber\\&&+
\frac{\lambda_4(\mu_1+2\mu_2+3\mu_3+4\mu_4)}{50}
+\frac{
\sum\limits_{i=1}^{4}k_{i}^{2}-k_{1}k_{2}-k_{2}k_{3}-k_{3}k_{4}
-\sum\limits_{i=1}^{4}k_i-k_i(\lambda_i+\mu_i)}{10}. \hspace{2.0cm} (3.18)\nonumber
\end{eqnarray}}

From (2.7) and (2.9), we know that $L_0\simeq sl_5$. And $L_0$-module $L_{-1}$ is isomorphic to fundamental module $V(\omega_2)=W$. The set of its weights
and the basis for the corresponding weight space are tabulated in Table 1.
The $L_0$-module $L_{1}$ is isomorphic to highest weight module $V(\omega_1+\omega_2)$ with  lowest weight vector $x_{5}d_{45}$ (cf. Table 9).
\vspace{0.1cm}

{\it Lemma 3.6 } \quad {\it The   wedge module $\Lambda^{k}W$ ($k \in \overline{1,10}$) for $sl_5$ are decomposed multiplicity freely into  irreducible components, which are listed in Table 2.}

{\it Proof} \quad
By Weyl's dimension formula,  we get:
$\mbox{dim}V(\omega_1+\omega_3)=45$,
$\mbox{dim}V(2\omega_3)=50$,
$\mbox{dim}V(2\omega_1+\omega_4)=70$,
$\mbox{dim}V(3\omega_1)=35$,
$\mbox{dim}V(\omega_1+\omega_3+\omega_4)=175$,
$\mbox{dim}V(2\omega_1+\omega_3)=126$,
$\mbox{dim}V(\omega_2+2\omega_4)=126$.
Since $\mbox{dim}W=10$, $\mbox{dim}\wedge^{k}W=C_{10}^{k}$. Thus the decomposition follows through comparing the dimensions of both sides.
\hspace{1cm} $\Box$\vspace{0.1cm}

{\it Lemma 3.7} \quad {\it The  tensor  module $V(k\omega_4)\otimes \Lambda^{n}W$ ($k\in \mathbb{N}, \ n \in \overline{1,10}$) for $sl_5$ are decomposed  into  irreducible components, which are listed in Table 3.}
\vspace{0.1cm}

For any  highest weight module $ V(\mu)$ of simple Lie algebra $sl_{5}$, denote the set of  its weights  by $\Pi(\mu)$, which are listed by $\{\overrightarrow{w}_{j}^{\mu} \ |  \  j \in \overline{1,|\Pi(\mu)|}\}$. Let $\{v_{j,k}^{\mu}  \ |  \  j \in \overline{1,|\Pi(\mu)|} , k \in \overline{1,\mbox{mult}(\overrightarrow{w}_{j}^{\mu} )}\}$ be the Verma basis for the weight space of weight $\overrightarrow{w}_{j}^{\mu}$, where $ \mbox{mult}(\overrightarrow{w}_{j}^{\mu} )$ denotes the multiplicity of the weight $ \overrightarrow{w}_{j}^{\mu} $. \vspace{0.1cm}

{\it Lemma 3.8}\quad {\it For
$\mu \in \{  \omega_1+\omega_3, 2\omega_1+\omega_4, 3\omega_1\}$, the set  $\Pi(\mu)$ of weights for  $V(\mu)$ and their corresponding Verma bases for every weight space  are listed in Table 5-Table 8 in the Appendix.}

{\it Proof}\quad Assume $\mu=\sum\limits_{i=1}^{4}m_i\omega_i$. The set  $\Pi(\mu)$ is obtained by the algorithm from [W]. The Verma bases for the weight space with weight
$\mu-\sum\limits_{i=1}^{4}k_i\alpha_i$ are (cf. [LMNP], [RS]):
$$(f_1^{a_{10}}f_{2}^{a_9}
f_{3}^{a_8}f_{4}^{a_7})(f_{1}^{a_6}f_{2}^{a_5}f_{3}^{a_4})(f_{1}^{a_3}f_{2}^{a_2})
f_{1}^{a_1}v_{\mu},\eqno(3.19)
$$where $$a_{10}+a_6+a_3+a_1=k_1,a_{9}+a_5+a_2=k_2,a_{8}+a_4=k_3, a_7=k_4,$$
$$0\leq a_1 \leq m_1, 0 \leq a_2 \leq m_2+a_1, 0\leq a_3 \leq \mbox{min}(m_2,a_2),$$$$
0\leq a_4 \leq m_3+a_2, 0 \leq a_5 \leq\mbox{min}(m_3+a_3,a_4),
0 \leq  a_6 \leq \mbox{min}(m_3,a_5),$$$$ 0 \leq a_7 \leq m_4+a_4, 0 \leq a_8 \leq \mbox{min}(m_4+a_5, a_7),
0 \leq a_9 \leq \mbox{min}(a_4+a_6, a_8), 0 \leq a_{10} \leq \mbox{min}(m_4,a_9).\eqno(3.20)
$$
\hspace{1cm} $\Box$\vspace{0.1cm}

{\it Remark 3.9} \quad {The coordinates of the weights appearing in Table1, Table5-Table9 are with respect to the ordered basis $\omega_1,\omega_2,\omega_3,\omega_4$. The basis of every weight space appearing in these tables  are Verma basis.}

\section{Singular vectors for GVM of $\mbox{E}(5, 10)$ }

\quad\quad   In Section 4.1, we analyze the
detailed structure of the
generalized Verma module $M$ over $\mbox{E}(5, 10)$.
It turns out that  there is a grading on  $M$ and each graded subspace is a finite dimensional $sl_{5}$- module (cf. Equation (4.6) and (4.7)). Moreover, any singular vector for $M$  is controlled by its leading term through an exponential-like differential operator (cf. Equation (4.27) ). In section 4.2,  we inductively prove that any leading term must satisfy three equations, i.e. (4.28), (4.35) and (4.37). Based on the Lemmas in Section 3, we
simplify these three differential equations  and  prove that
 any singular vector   is of degree less than or equal to four. Also,
the leading term of any singular vector must lie in one of the tensor decomposition  of four  tensor product module for  $sl_5$ (cf. Theorem 4.6 ).

\subsection{Gradation for GVM }

\quad  Set $$T=\{0,1\}, \ \  T'=\{(45), (35) , (25) , (15), (34), (24), (14), (23) , (13), (12)\}. \eqno(4.1)$$
 Define  order $ ``\prec''$ on the set $T'$ by:
$$(45)\prec (35) \prec (25) \prec (15)\prec (34)\prec (24)\prec (14)\prec (23) \prec (13)\prec (12).\eqno(4.2) $$
For $\underline{n}=(n_{12},n_{13},n_{14},n_{24},n_{34},n_{15},n_{25},
n_{35},n_{45})\in T^{10}$ and $ \underline{m} \in \mathbb{N}^{5}$, we take the following notations:
$$ \underline{n}\pm\varepsilon_{ij}=( n_{12},\cdots,n_{ij}\pm 1,\cdots, n_{45}), \ \underline{m}\pm\varepsilon_{i}=( m_1,\cdots,m_i \pm 1,\cdots, m_5). \eqno(4.3)$$
 Let$$
d^{\underline{n}}=d_{12}^{n_{12}}d_{13}^{n_{13}}d_{23}^{n_{23}}
d_{14}^{n_{14}}d_{24}^{n_{24}}d_{34}^{n_{34}}d_{15}^{n_{15}}d_{25}^{n_{25}}
d_{35}^{n_{35}}d_{45}^{n_{45}},\eqno(4.4)$$$$  \partial^{\underline{m}}=\partial_{1}^{m_1}
\partial_{2}^{m_2}\partial_{3}^{m_3}
\partial_{4}^{m_4}\partial_{5}^{m_5}.\eqno(4.5)$$
Then the induced module $M$ is spanned by $\{\partial^{\underline{m}}d^{\underline{n}}v_{\nu} \ | \ \underline{n}\in  T^{10}, \underline{m} \in \mathbb{N}^5, \nu \in \Pi(\lambda) \}$.
Define
 $$\partial^{m}\wedge^{n}V=\mbox{Span}
 \{\partial^{\underline{m}}d^{\underline{n}}v \ | \  |\underline{m}|=m, |\underline{n}|=n\}, \ M_{k}=\mbox{Span}\{\partial^{\underline{m}}d^{\underline{n}}v \ | \ 2m+n=k\}.\eqno(4.6)
 $$
Then $$M=\bigoplus\limits_{k \in \mathbb{N}}M_k.\eqno(4.7)$$

{\it Definition 4.1} {\it We say any nonzero vector of $M_{k}$ is of degree $k$.}
\vspace{0.1cm}

The equations $$[x_i\partial_{x_j},d_{kl}]=\delta_{j,k}d_{il}-\delta_{jl}d_{ik},[x_5d_{45},d_{12}]=x_5\partial_{x_3},  \ [x_5d_{45},d_{13}]=-x_5\partial_{x_2}, \ [x_5d_{45},d_{23}]
=x_5\partial_{x_1},$$$$
[x_5d_{45},d_{i4}]=0(i \in \overline{1,3}), \ [x_5d_{45},d_{i5}]=0(i \in \overline{1,4}).\eqno(4.8)
$$
yield
$$ L_0 .\partial^{m}\wedge^{n}V\subseteq \partial^{m}\wedge^{n}V+\partial^{m+1}\wedge^{n-2}V, \ x_{5}d_{45}. M_{k} \subseteq M_{k-1}. \eqno(4.9)$$
That is to say, every graded vector subspace  $M_{k}$ is an $sl_5$-module and every singular vector for $\mbox{E}(5, 10)$-module $M$ is in a certain graded subspace $M_{k}$.

 In the following of this section, we  consider the
maximal vectors for $sl_5$-module $M_{k}$. On any linear vector space $\partial^{m}\wedge^{n}V$, we define the following linear operators:
$$(-1)^{|ij|}:\partial^{m}\wedge^{n}V \rightarrow \partial^{m}\wedge^{n}V;\partial^{\underline{m}}d^{\underline{n}}v
\mapsto (-1)^{\sum\limits_{(kl)\prec (ij)}n_{kl}}
\partial^{\underline{m}}d^{\underline{n}}v,$$
$$(-1)^{|ij,kl|}:\partial^{m}\wedge^{n}V \rightarrow \partial^{m}\wedge^{n}V;
\partial^{\underline{m}}d^{\underline{n}}v
\mapsto(-1)^{\sum\limits_{(kl)\prec (pq)\prec (ij)}n_{pq}}
\partial^{\underline{m}}d^{\underline{n}}v,$$
$$y_{ij}\partial_{y_{kl}}:
\partial^{m}\wedge^{n}V \rightarrow \partial^{m}\wedge^{n}V;
\partial^{\underline{m}}d^{\underline{n}}v
\mapsto n_{kl}\partial^{\underline{m}}
d^{\underline{n}+\varepsilon_{ij}-\varepsilon_{kl}}v,$$$$
z_i:\partial^{m}\wedge^{n}V \rightarrow \partial^{m+1}\wedge^{n}V;
\partial^{\underline{m}}d^{\underline{n}}v\mapsto
\partial^{\underline{m}+\varepsilon_i}d^{\underline{n}}v,$$$$
\partial_{z_i}:\partial^{m}\wedge^{n}V \rightarrow \partial^{m-1}\wedge^{n}V;
\partial^{\underline{m}}d^{\underline{n}}v\mapsto m_i
\partial^{\underline{m}-\varepsilon_i}d^{\underline{n}}v,
$$
$$E_{i,j}:\partial^{m}\wedge^{n}V \rightarrow \partial^{m}\wedge^{n}V;
\partial^{\underline{m}}d^{\underline{n}}v
\mapsto \partial^{\underline{m}}d^{\underline{n}}(E_{i,j}.v).\eqno(4.10)$$
Set
$$ (x_i\partial_{x_j})_0^{'}=\sum\limits_{k\in \overline{1,5},k \neq i,j}(-1)^{|ki,kj|}y_{ki}\partial_{y_{kj}}, \ (x_i\partial_{x_j})_0=-z_j\partial_{z_i}+(x_i\partial_{x_j})_0^{'}+E_{i,j} (i \neq j);\eqno(4.11)$$
$$(x_3\partial_{x_4})_{-2}=z_5\partial_{y_{14}}
\partial_{y_{24}},\eqno(4.12)$$
$$(x_4\partial_{x_5})_{-2}=z_1\partial_{y_{25}}\partial_{y_{35}}
+(-1)^{1+|15,35|}
z_2\partial_{y_{15}}\partial_{y_{35}}
+z_3\partial_{y_{15}}\partial_{y_{25}}.\eqno(4.13)
$$
Using these settings, we could formulate  the equation (2.12) in  the following explicit form:
$$x_1\partial_{x_{2}}=(x_1\partial_{x_{2}})_0, \quad
x_2\partial_{x_{3}}=(x_2\partial_{x_{3}})_0,
$$
$$x_3\partial_{x_{4}}=(x_3\partial_{x_{4}})_0
+(x_3\partial_{x_{4}})_{-2}, \ x_4\partial_{x_{5}}=(x_4\partial_{x_{5}})_0
+(x_4\partial_{x_{5}})_{-2}.\eqno(4.14)$$
 According to the Cartan subalgebra of $L_0$, $M$ can be decomposed  to the following direct sum of subspaces:
$$M=\bigoplus\limits_{\mu\in \Gamma}M^{\mu}, \ M^{\mu}=\mbox{Span}\{ \partial^{\underline{m}}d^{\underline{n}}v_{\nu} \ | \
(x_{i}\partial_{x_{i}}-x_{i+1}\partial_{x_{i+1}}).\partial^{\underline{m}}d^{\underline{n}}v_{\nu}=
\mu_i\partial^{\underline{m}}d^{\underline{n}}v_{\nu}
 \}, \eqno(4.15)$$where
 $$\mu_i=m_{i+1}-m_i+t_i(\underline{n})+\nu_i,$$$$
 t_1(\underline{n})=n_{13}+n_{14}+n_{15}-n_{23}-n_{24}-n_{25}, \
t_2(\underline{n})=n_{12}+n_{24}+n_{25}-n_{13}-n_{34}-n_{35},$$$$
t_3(\underline{n})=n_{13}+n_{23}+n_{35}-n_{14}-n_{24}-n_{45}, \
t_4(\underline{n})=n_{14}
+n_{24}+n_{34}-n_{15}-n_{25}-n_{35}.\eqno(4.16)$$
\vspace{0.1cm}
For any vectors  $v \in M^{\mu}$, we say that it is of weight $\mu$ and denote $\mbox{wt}(v)=\mu, \ |\mbox{wt}(v)|=|\mu|=\sum\limits_{i=1}^{4}\mu_i.$

{\it Proposition 4.2}\quad {\it The differential operators  $(x_i\partial_{x_j})_0 (1\leq i \neq j \leq 5)$ and $ x_{i}\partial_{x_{i}}-x_{i+1}\partial_{x_{i+1}} (i \in \overline{1,4}) $
give
 every vector space
$\partial^{m}\wedge^{n}V$  an $sl_5$-module structure, which is
 isomorphic to tensor module $ V(m\omega_4)\otimes \wedge^{n}W \otimes V$ for $sl_5$. }

 {\it Proof} \quad The module isomorphism is given by:
 $$\phi: V(m\omega_4)\otimes \wedge^{n}W \otimes V \rightarrow \partial^{m}\wedge^{n}V ;
 \partial^{\underline{m}}\otimes (d_{i_1j_1}\wedge \cdots \wedge d_{i_nj_n} ) \otimes v
 \mapsto
 \partial^{\underline{m}}d_{i_1j_1}\cdots d_{i_nj_n}v.\eqno(4.17)
 $$ \hspace{1cm} $\Box$\vspace{0.1cm}

Denote $$\Gamma_k=\{(m,n)\in \mathbb{N}^2 \ | \ 2m+n=k\}\eqno(4.18)$$
 For any $(m,n) \in \Gamma_k$, let $$\Gamma_k^{(m,n)}=\{(m',n') \in \Gamma_k \ | \ m' \geq m\}\eqno(4.19)$$
 Assume $\xi \in M_k$ is any $E(5,10)$ singular vector. Then there exists $(m,n) \in \Gamma_k$ such that $$\xi \in \bigoplus\limits_{(m'n') \in \Gamma_k^{(m,n)}}\partial^{m'}\wedge^{n'}V\eqno(4.20)$$
  For emphasis, we write $$
  \xi=\xi^{m,n}=\sum\limits_{(m'n') \in \Gamma_k^{(m,n)}}\xi_{m',n'}\eqno(4.21)
  $$
We say that $\xi_{m,n}$ is {\it  the leading term } of   $\xi^{m,n}$.
It follows from (4.14) that
$\xi^{m,n}$ must satisfy the following equations inductively:
$$(x_i\partial_{x_{i+1}})_0.\xi_{m,n}=0 ( i \in \overline{1,4}),$$
$$(x_i\partial_{x_{i+1}})_{-2}.\xi_{m',n'}+
(x_i\partial_{x_{i+1}})_{0}.\xi_{m'+1,n'-2}=0, \   i \in \overline{1,4}, \
  (m',n') \in \Gamma_k^{(m,n)}.\eqno(4.22)$$
\vspace{0.1cm}

{\it Remark 4.3}\quad  {\it From  Proposition 4.2 and (4.22),  we derive that the leading term $\xi_{m,n}$ of any singular vector $\xi=\xi^{m,n}$   is also a singular vector of the tensor product module $ V(m\omega_4)\otimes \wedge^n V(\omega_2)\otimes V(\lambda)$ for simple Lie algebra $sl_5$. In the following, we will point out that any singular vector $\xi=\xi^{m,n}$  is completely controlled by its leading term $\xi_{m,n}$ through certain exponential-like differential operator. }
\vspace{0.1cm}

Set $$P=\sum\limits_{
(kl)\prec (ij) \in T', m \in \overline{1,5}
  }\varepsilon_{mijkl}(-1)^{|ij,kl|}z_m
\partial_{y_{ij}}\partial_{y_{kl}},\eqno(4.23)$$where $\varepsilon_{mijkl}$ is defined in (2.6).
The operator $P$ is checked to satisfy the following equations:
 $$[(x_3\partial_{x_{4}})_{-2},P]=0, [(x_4\partial_{x_{5}})_{-2},P]=0, \ [(x_1\partial_{x_{2}})_0,P]=0, [(x_2\partial_{x_{3}})_0,P]=0, $$
 $$[(x_3\partial_{x_{4}})_0,P]=2(x_3\partial_{x_{4}})_{-2}, \ [(x_4\partial_{x_{5}})_0,P]=2(x_4\partial_{x_{5}})_{-2}
.\eqno(4.24)$$
Inductively,
$$[(x_i\partial_{x_{i+1}})_0,P^k]=[x_i\partial_{x_{i+1}},P^k]
=2kP^{k-1}(x_i\partial_{x_{i+1}})_{-2}, k \in \mathbb{N}.\eqno(4.25)$$
It implies
$$x_i\partial_{x_{i+1}}.e^{-\frac{1}{2}P}\xi_{m,n}=0, \quad
(x_i\partial_{x_{i+1}})_{0}.e^{\frac{1}{2}P}\xi^{m,n}=0.\eqno(4.26)$$
Thus  we prove the following formula:

{\it Proposition 4.3}\quad {\it
Assume $ \xi^{m,n}=\sum\limits_{(m'n') \in \Gamma_k^{(m,n)}}\xi_{m',n'} \in M_k$ is any  singular vector for $E(5,10)$-module $M$, then
 $$\xi^{m,n}=e^{-\frac{1}{2}P}\xi_{m,n}.\eqno(4.27)$$}

\subsection{Singular vectors for GVM }

 \quad
 \quad  In this section, we continue  the discussion concerning the equation (2.13) in Definition 2.1.
Recall the notations in (4.10), set
$$(x_5d_{45})_1
=(-1)^{1+|45|}\partial_{z_5}y_{45},\eqno(4.28)$$
\begin{eqnarray}
(x_5d_{45})_{-1}&
=&-z_3\partial_{z_5}\partial_{y_{12}}+(-1)^{|13|}
z_2\partial_{z_5}\partial_{y_{13}}+(-1)^{1+|23|}
z_1\partial_{z_5}\partial_{y_{23}}\nonumber\\&+&
\partial_{y_{12}}E_{53}+(-1)^{1+|13|}\partial_{y_{13}}E_{52}+(-1)^{|23|}
\partial_{y_{23}}E_{51}\nonumber
\\&+&
(-1)^{|13,15|}y_{15}\partial_{y_{12}}\partial_{y_{13}}+
(-1)^{|23,25|}y_{25}\partial_{y_{12}}\partial_{y_{23}}+
(-1)^{|13|+|23,35|}y_{35}\partial_{y_{13}}\partial_{y_{23}}\nonumber\\&+&
(-1)^{1+|34,45|}y_{45}\partial_{y_{12}}\partial_{y_{34}}+
(-1)^{1+|23|+
|14,45|}y_{45}\partial_{y_{23}}\partial_{y_{14}}+
(-1)^{|13|+|24,45|}y_{45}\partial_{y_{13}}\partial_{y_{24}},\nonumber
\end{eqnarray}$$\eqno(4.29)$$
\begin{eqnarray}
(x_5d_{45})_{-3}&
=&(-1)^{|23,34|}z_1\partial_{y_{12}}\partial_{y_{23}}\partial_{y_{34}}
+(-1)^{1+|23,24|+|13|}
z_1\partial_{y_{13}}\partial_{y_{23}}\partial_{y_{24}}\nonumber\\&+&(-1)^{|13|}
z_2\partial_{y_{13}}\partial_{y_{23}}\partial_{y_{14}}+
(-1)^{1+|13,34|}
z_2\partial_{y_{12}}\partial_{y_{13}}\partial_{y_{34}}\nonumber
\\&+&(-1)^{|13,24|}z_3\partial_{y_{12}}\partial_{y_{13}}\partial_{y_{24}}
-z_3\partial_{y_{12}}\partial_{y_{23}}\partial_{y_{14}}-z_4
\partial_{y_{12}}\partial_{y_{13}}\partial_{y_{23}}.
\nonumber \hspace{3.8cm}(4.30)
\end{eqnarray}
It follows from the equation (4.8) that
$$
 x_5d_{45}.\partial^{m}\wedge^{n}V\subseteq \partial^{m-1}\wedge^{n+1}V+
 \partial^{m}\wedge^{n-1}V+\partial^{m+1}\wedge^{n-3}V, \ x_5d_{45}.M_{k}\subseteq M_{k-1},$$
$$
x_5d_{45}=(x_5d_{45})_1+(x_5d_{45})_{-1}+(x_5d_{45})_{-3}.\eqno(4.31)
$$
Furthermore,  $\xi^{m,n}$ must satisfy the following equations inductively:
  $$(x_5d_{45})_1.\xi_{m,n}=0, \ (x_5d_{45})_{-1}.\xi_{m,n}+(x_5d_{45})_1.\xi_{m+1,n-2}=0,$$$$
  (x_5d_{45})_{-3}.\xi_{m,n}
  +(x_5d_{45})_{-1}.\xi_{m+1,n-2}+(x_5d_{45})_1.\xi_{m+2,n-4}=0,$$
$$(x_5d_{45})_{-3}.\xi_{m',n'}
  +(x_5d_{45})_{-1}.\xi_{m'+1,n'-2}+(x_5d_{45})_1.\xi_{m'+2,n'-4}=0, \
 \mbox{ for \ any} \ ( m',n') \in \Gamma_k^{(m,n)}.
  \eqno(4.32)$$\vspace{0.3cm}
 Applying (4.27), the leading term $\xi_{m,n}$ should be killed by the following three operators:
$$ (x_5d_{45})_1,  \ (x_5d_{45})_{-1}+(x_5d_{45})_1(-\frac{1}{2}P), \
(x_5d_{45})_{-3}+(x_5d_{45})_{-1}(-\frac{1}{2}P) +(x_5d_{45})_1(\frac{1}{8}P^2).\eqno(4.33)$$
We can reduce the last two differential operators to  be of more explicit forms.
Indeed, the following relations are easily checked:
\begin{eqnarray}
[(x_5d_{45})_1,P]&=&-z_3\partial_{z_5}\partial_{y_{12}}+(-1)^{|13|}
z_2\partial_{z_5}\partial_{y_{13}}+(-1)^{1+|23|}
z_1\partial_{z_5}\partial_{y_{23}}\nonumber\\&+&
(-1)^{1+|34,45|}y_{45}\partial_{y_{12}}\partial_{y_{34}}+
(-1)^{1+|23|+
|14,45|}y_{45}\partial_{y_{23}}\partial_{y_{14}}+
(-1)^{|13|+|24,45|}y_{45}\partial_{y_{13}}\partial_{y_{24}},\nonumber
\end{eqnarray}
$$[[(x_5d_{45})_1,P],P]=2 (x_5d_{45})_{-3}+2z_4\partial_{y_{12}}\partial_{y_{13}}\partial_{y_{23}},
\ [(x_5d_{45})_{-1},P]=3(x_5d_{45})_{-3}.\eqno(4.34)$$
Therefore,
\begin{eqnarray}
&&[(x_5d_{45})_{-1}+(x_5d_{45})_1(-\frac{1}{2}P)]\xi_{m,n}\nonumber\\
&&=
\{(x_5d_{45})_{-1}-\frac{1}{2}P(x_5d_{45})_1-\frac{1}{2}
[(x_5d_{45})_1,P]\}\xi_{m,n}
\nonumber
\\&&\stackrel{\mbox{by}(4.34)}{=}[-\frac{1}{2}z_3\partial_{z_5}\partial_{y_{12}}+\frac{1}{2}(-1)^{|13|}
z_2\partial_{z_5}\partial_{y_{13}}+(-1)^{1+|23|}\frac{1}{2}
z_1\partial_{z_5}\partial_{y_{23}}\nonumber\\
&&+\partial_{y_{12}}E_{53}+(-1)^{1+|13|}\partial_{y_{13}}E_{52}+(-1)^{|23|}
\partial_{y_{23}}E_{51}\nonumber\\&&+
(-1)^{|13,15|}y_{15}\partial_{y_{12}}\partial_{y_{13}}
+
(-1)^{|23,25|}y_{25}\partial_{y_{12}}\partial_{y_{23}}
\nonumber\\&&+
(-1)^{|13|+|23,35|}y_{35}\partial_{y_{13}}\partial_{y_{23}}
+\frac{1}{2}(-1)^{1+|34,45|}y_{45}\partial_{y_{12}}\partial_{y_{34}}
\nonumber\\&&+\frac{1}{2}
(-1)^{1+|23|+
|14,45|}y_{45}\partial_{y_{23}}\partial_{y_{14}}+\frac{1}{2}
(-1)^{|13|+|24,45|}y_{45}\partial_{y_{13}}\partial_{y_{24}}]\xi_{m,n}.
\nonumber\\&&\stackrel{\mbox{by}(4.11)}{=}[(E_{53}+\frac{1}{2}
(x_{5}\partial_{x_{3}})_{0}'-\frac{1}{2}z_3\partial_{z_5})
\partial_{y_{12}}+(E_{52}+\frac{1}{2}(x_{5}\partial_{x_{2}})_{0}'
-\frac{1}{2}z_2\partial_{z_5})
(-1)^{1+|13|}\partial_{y_{13}}\nonumber\\&&+(E_{51}+\frac{1}{2}
(x_{5}\partial_{x_{1}})_{0}'
-\frac{1}{2}z_1\partial_{z_5})(-1)^{|23|}\partial_{y_{23}}
]\xi_{m,n}\nonumber\hspace{7.2cm}(4.35)
\end{eqnarray}
Hence,
$$P(x_5d_{45})_{-1}\xi_{m,n}=\frac{1}{2}P(x_5d_{45})_{1}P\xi_{m,n}=
\frac{1}{2}P[(x_5d_{45})_1,P]\xi_{m,n}.\eqno(4.36)$$
Furthermore, (4.34) and (4.35) imply that
\begin{eqnarray}
&&[(x_5d_{45})_{-3}+(x_5d_{45})_{-1}(-\frac{1}{2}P) +(x_5d_{45})_1(\frac{1}{8}P^2)]\xi_{m,n}\nonumber\\&&=
\{(x_5d_{45})_{-3}-\frac{3}{2}(x_5d_{45})_{-3}
-\frac{1}{2}P(x_5d_{45})_{-1}\nonumber\\&&+\frac{1}{8}[[(x_5d_{45})_{1},P],P]
+\frac{1}{8}P[(x_5d_{45})_{1},P]+\frac{1}{8}P(x_5d_{45})_{1}P\}\xi_{m,n}
\nonumber\\&&
\stackrel{\mbox{by}(4.34),(4.35),(4.36)}{=}[-\frac{1}{4}(x_5d_{45})_{-3}+\frac{1}{4}
z_4\partial_{y_{12}}\partial_{y_{13}}\partial_{y_{23}}]\xi_{m,n}
\nonumber\\&&
=\frac{1}{4}[z_1(-1)^{|13|+|23,24|}\partial_{y_{13}}\partial_{y_{23}}\partial_{y_{24}}
-z_1(-1)^{|23,34|}\partial_{y_{12}}\partial_{y_{23}}\partial_{y_{34}}
\nonumber\\&&
+z_2(-1)^{|13,34|}\partial_{y_{12}}\partial_{y_{13}}\partial_{y_{34}}
-z_2(-1)^{|13|}\partial_{y_{13}}\partial_{y_{23}}\partial_{y_{14}}
\nonumber\\&&
+z_3\partial_{y_{12}}\partial_{y_{23}}\partial_{y_{14}}
-z_3(-1)^{|13,24|}\partial_{y_{12}}\partial_{y_{13}}
\partial_{y_{24}}+
2z_4\partial_{y_{12}}\partial_{y_{13}}\partial_{y_{23}}]\xi_{m,n}.
\nonumber\hspace{4.2cm}(4.37)
\end{eqnarray}
 Denote the set of all the highest weight vectors for tensor modules
$V(m\omega_4)\bigotimes \wedge^{n}W$ by $$S_{m,n}=\{e_{m,n}^1,\cdots,e_{m,n}^{\nu(m,n)}\}\eqno(4.38)$$
By Lemma 3.3, any singular vector $\xi_{m,n}$ of the tensor product module $ V(m\omega_4)\bigotimes \bigwedge^n V(\omega_2)\bigotimes V(\lambda)$ for $sl_5$ can be written by the following form:
$$\xi_{m,n}=e_{m,n}^{i}\otimes v_{\vartheta}+\cdots.\eqno(4.39)$$
We consider the set
$$S_{m,n}'=\{e_{m,n}^i\in S_{m,n} \ | \ (-1)^{1+|45|}\partial_{z_5}y_{45}.\phi(e_{m,n}^{i}\otimes v)=0,$$$$
 [(x_5d_{45})_{-3}+(x_5d_{45})_{-1}(-\frac{1}{2}P) +(x_5d_{45})_1(\frac{1}{8}P^2)]
.\phi(e_{m,n}^{i}\otimes v)=0, \ \forall \
 v \in V(\lambda)\}.\eqno(4.40) $$

{\it Proposition 4.4 } \quad {\it All the non empty set of $S_{m,n}'$
are listed in the following:
$$S_{0,0}'=\{1\},
S_{0,1}'=\{d_{12}\},
S_{0,2}'=\{d_{12}\wedge d_{13}\},
 S_{0,3}'=\{d_{12}\wedge d_{13} \wedge d_{14}\},
S_{0,4}'=\{d_{12}\wedge d_{13} \wedge d_{14} \wedge d_{15}\}.$$}
{\it Proof}\quad Let $V(\mu)$ be any highest weight module appearing in the decomposition of the  $sl_5$ wedge module $\wedge^{n}V(\omega_2)$ (cf. Table 2). And the highest weights appearing in the decomposition of $V(m\omega_4)\otimes V(\mu)$ are listed in Table 3.
By Lemma 3.3, the maximal vector in the tensor module $V(m\omega_4)\otimes V(\mu)$  is written as:
$$\partial_{5}^{m}\otimes l_{\mu}+\sum\limits_{\underline{q}\in \mathbb{N}^{5}}\partial^{\underline{q}}\otimes v_{\underline{q}},\eqno(4.41)$$
where $l_{\mu}$ satisfies $E_{12}l_{\mu}=0, E_{23}l_{\mu}=0, E_{34}l_{\mu}=0,
E_{45}^{m+1}l_{\mu}=0$. By detailed calculation, we get all the $l_{\mu}$, which are listed in Table-4.
A straightforward but messy check case by case shows that
the assertion holds.
\hspace{1cm} $\Box$\vspace{0.1cm}

{\it Remark 4.5}\quad {\it For the 10-tuple $ d_{12}\wedge d_{13} \wedge  \cdots \wedge d_{45}$, we use the notation ${\hat d_{i_{1}j_{1}} }\wedge {\hat d_{i_{2}j_{2}} }\wedge \cdots \wedge {\hat d_{i_{k}j_{k}} } $  to denote the (10-k)-tuple where $d_{i_{1}j_{1}}, \cdots,
d_{i_{k}j_{k}}$ have been omitted in Table 4. }\vspace{0.1cm}

To summarize Proposition 4.3 and Proposition 4.4, we have proved the following statement in this section:
\vspace{0.1cm}

{\it Theorem  4.6}\quad {\it
 Any singular vector for $E(5,10)$-module $M$ is of the form:
 $$\xi^{0,n}=e^{-\frac{1}{2}P}\xi_{0,n}, \ \ n \in \overline{1,4}$$
where the leading term $ \xi_{0,n}$ satisfies the equation:
$$[(E_{53}+\frac{1}{2}(x_{5}\partial_{x_{3}})_{0}')
\partial_{y_{12}}+(E_{52}+\frac{1}{2}(x_{5}\partial_{x_{2}})_{0}')
(-1)^{1+|13|}\partial_{y_{13}}+(E_{51}+\frac{1}{2}(x_{5}\partial_{x_{1}})_{0}'
)(-1)^{|23|}\partial_{y_{23}}
].\xi_{0,n}=0.\eqno(4.42)$$
  Moreover, $ \xi_{0,n}$ is the maximal vector lying in one of  the following $sl_5$-tensor modules:
 $$\xi_{0,1}\in V(\omega_2)\otimes V(\lambda), \xi_{0,2}\in V(\omega_1+\omega_3)\otimes V(\lambda), \xi_{0,3}\in V(2\omega_1+\omega_4)\otimes V(\lambda),\xi_{0,4}\in V(3\omega_1)\otimes V(\lambda).\eqno(4.43) $$
 }

\section {Singular vectors degree  by degree }
\quad \quad In this section, we work out   all the singular vectors in
 Theorem  4.6 explicitly degree  by degree. Before turning to the calculation, we introduce some formula which we are going to use in the remainder of this section.

Recall that  we could endow any vector space $\Lambda^{m}V$ an $sl_5$- module structure with the action $(x_i\partial_{x_i})_0^{'}-(x_{i+1}\partial_{x_{i+1}})_0^{'} \  (i \in \overline{1,4}), \ (x_{i}\partial_{x_{j}})_{0}' \ (i \neq j)  $, which is isomorphic to the tensor product module $\wedge^{m}V(\omega_2)\otimes V(\lambda)$ in Section 4.1. Now we
define the following   differential operator on the $sl_5$- module $\Lambda^{m}V$:
$$\tilde{c}=\frac{1}{10}[\sum\limits_{i=1}^{4}
((x_i\partial_{x_i})_0^{'}-(x_{i+1}\partial_{x_{i+1}})_0^{'})h_{i}^{*}
+\sum\limits_{1 \leq i\neq j \leq 5}(x_i\partial_{x_j})_0^{'}E_{ji}
],\eqno(5.1)$$
$$T_{i,jkl}=[E_{ij}+\frac{1}{2}(x_{i}\partial_{x_{j}})_{0}'](-1)^{|kl
|}\partial_{y_{kl}}
+[E_{ik}+\frac{1}{2}(x_{i}\partial_{x_{k}})_{0}']
(-1)^{1+|jl|}\partial_{y_{jl}}+
[E_{il}+\frac{1}{2}(x_{i}\partial_{x_{l}})_{0}']
(-1)^{|jk|}
\partial_{y_{jk}}.\eqno(5.2)$$\vspace{0.1cm}

{\it Lemma  5.1}\quad {
\quad Assume $Q_{ij}^{0}\in \mbox{Span}_{\mathbb{F}}\{y_{ij}\partial_{y_{kl}} \ | \ 1 \leq i < j \leq 5, \ 1 \leq k < l \leq 5  \}$. Then $$\sum\limits_{1 \leq i < j \leq 5}Q_{ij}^{0}.(-1)^{|ij|}\partial_{y_{ij}}.\tilde{c}|_{\Lambda^{m}V}=\sum\limits_{1 \leq i < j \leq 5}Q_{ij}^{1}(-1)^{|ij|}\partial_{y_{ij}}|_{\wedge^{m}V},$$where
$$Q_{ij}^{1}=Q_{ij}^{0}(\tilde{c}+\sum\limits_{k=1}^{4}
\frac{s_{k}^{ij}}{10}h_{k}^{*})+\frac{1}{10}\sum\limits_{m\neq i, j}(Q_{im}^{0} E_{jm}-Q_{jm}^{0} E_{im}),$$
 $$s_{k}^{ij}\stackrel{\mbox{by}(4.16)}{=}t_k(\underline{n})- t_k(\underline{n}-\epsilon_{ij})\eqno(5.3)$$ for any
$k \in \overline{1,4},  \underline{n}\in T^{10}, (ij) \in S'$.}

{\it Proof}\quad Indeed, the formula (5.3) follows from:
$$[(-1)^{|ij|}\partial_{y_{ij}},\tilde{c}]|_{\wedge^{m}V}
=\frac{1}{10}(\sum\limits_{k\neq i,j}(-1)^{|ik|}\partial_{y_{ik}}E_{kj}- \sum\limits_{k\neq i,j}(-1)^{|jk|}\partial_{y_{jk}} E_{ki}+\sum\limits_{k=1}^{4}s_{k}^{ij}(-1)^{|ij|}\partial_{y_{ij}} h_{k}^{*})|_{\wedge^{m}V}.\eqno(5.4)$$\hspace{1cm} $\Box$\vspace{0.1cm}

 {\it Lemma  5.2}\quad { We could define the following intertwining operators between the $sl_5$-module $V(\omega_1+\omega_2)\otimes \Lambda^{m}V$ and  $\Lambda^{m-1}V$ by:
 $$T^m: V(\omega_1+\omega_2)\otimes \Lambda^{m}V \rightarrow \Lambda^{m-1}V;
 v_{30,1}^{\omega_1+\omega_2}\otimes \xi \mapsto T_{5,123}(\xi),\eqno(5.5)$$where $v_{30,1}^{\omega_1+\omega_2}$ is the lowest weight vector for $V(\omega_1+\omega_2)$ (cf. Table 9) and $\xi$ is any maximal vector in $sl_5$-module  $\Lambda^{m}V$.}

{\it Proof} \quad Since the $sl_5$-module $V(\omega_1+\omega_2)\otimes \Lambda^{m}V$ is generated by such vectors of $v_{30,1}^{\omega_1+\omega_2}\otimes \xi$ by part (4) of Lemma 3.3, the assertion follows from the following formula:
$$[(x_{i+1}\partial_{x_{i}})_{0}',T_{5,123}]|_{\wedge^{m}V}=0 , \ i \in \overline{1,4}, $$$$[(x_{s}\partial_{x_{t}})_{0}', T_{i,jkl}]|_{\wedge^{m}V}=\delta_{t,i}T_{s,jkl}-
\delta_{s,j}T_{i,tkl}-\delta_{s,k}T_{i,jtl}-\delta_{s,l}T_{i,jkt}. \eqno(5.6)$$
\hspace{1cm} $\Box$\vspace{0.1cm}

\subsection{Singular vectors of degree one }

\quad \quad
{\it Theorem 5.3}\quad {All the possible  degree one  singular vectors are listed in the following:$$d_{12}v_{\lambda},
 \ \mbox{where} \ \lambda=(m,n,0,0),  \ m, n \in \mathbb{N};$$
$$\prod\limits_{\mbox{wt}(d_{15}) < \sigma \leq \omega_2} \frac{\tilde{c}-\chi_{\sigma+\lambda}(\tilde{c})}
{\chi_{\mbox{wt}(d_{15})+\lambda}(\tilde{c})-
\chi_{\sigma+\lambda}(\tilde{c})}.d_{15}v_{\lambda},
 \ \mbox{where} \ \lambda=(m,0,0,n), \  m \in \mathbb{N},  \ 1 \leq n \in \mathbb{N}; $$
$$\prod\limits_{\mbox{wt}(d_{45}) < \sigma \leq \omega_2} \frac{\tilde{c}-\chi_{\sigma+\lambda}(\tilde{c})}
{\chi_{\mbox{wt}(d_{45})+\lambda}(\tilde{c})-
\chi_{\sigma+\lambda}(\tilde{c})}.d_{45}v_{\lambda}, \ \mbox{where} \ \lambda=(0,0,m,n),  \ 1 \leq m \in \mathbb{N}, \ n \in \mathbb{N}.  $$
}\vspace{0.1cm}

{\it Proof} \quad The leading term of any singular vector  of degree one can be written  as $$\xi_{0,1}=\sum\limits_{1 \leq i < j \leq 5}d_{ij}v_{ij}, \ v_{ij} \in V(\lambda), \eqno(5.7)$$ which should  satisfy  :
$$T_{5,123} .\xi_{0,1}= [\partial_{y_{12}}E_{53}+(-1)^{1+|13|}\partial_{y_{13}}E_{52}+(-1)^{|23|}
\partial_{y_{23}}E_{51}].\xi_{0,1}=0,\eqno(5.8)$$
i.e. $$E_{53}v_{12}-E_{52}v_{13}+E_{51}v_{23}=0.\eqno(5.9)$$
Note that $(x_i\partial_{x_j})_{0}.\xi_{0,1}=0$ $(1\leq i < j \leq 5)$ imply that
$$ v_{13}=-E_{23}v_{12}, \ v_{23}=E_{13}v_{12}=-E_{12}v_{13}, \ v_{14}=-E_{34}v_{13}, \  v_{15}=-E_{25}v_{12}=-E_{35}v_{13}, \ $$$$ v_{25}=E_{15}v_{12}=-E_{35}v_{23}, \ v_{35}=E_{15}v_{13}=-E_{23}v_{25}=E_{25}v_{23}, \ v_{45}=-E_{34}v_{35}. \eqno(5.10)$$ Obviously, $v_{12}\neq 0$.

Case 1. wt($\xi_{0,1})=$ wt$(d_{12}v_{\lambda})$.

In this case, $v_{13}=v_{23}=0$, $v_{12}=v_{\lambda}$. And (5.9) implies that $E_{5,3}.v_{12}=E_{5,3}.v_{\lambda}=0$. That is  to say,
$\lambda=(m,n,0,0), (m,n) \in \mathbb{N}^2$.

Case 2.  wt($\xi_{0,1})\in \{ \mbox{wt}(d_{13}v_{\lambda}), \mbox{wt}(d_{14}v_{\lambda}), \mbox{wt}(d_{15}v_{\lambda})\}$.

In these three  cases, we have $v_{23}=0, v_{13}\neq 0$.

Case 2.1 \quad $wt(\xi_{0,1})\in \{ \mbox{wt}(d_{13}v_{\lambda}),  \mbox{wt}(d_{14}v_{\lambda})\}$

In these two cases, $v_{15}=0, v_{13}\neq 0$,
  $(h_2+h_3+h_4).v_{13}=(\lambda_2+\lambda_3+\lambda_4)v_{13}.$
Hence,
\begin{eqnarray}
\nonumber&&0=E_{25}(E_{53}.v_{12}-E_{52}.v_{13})=(E_{23}+E_{53}E_{25})v_{12}
-(h_2+h_3+h_4+E_{52}E_{25})v_{13}\\
\nonumber
&=&-(1+h_2+h_3+h_4)v_{13}-E_{53}v_{15}=-(1+h_2+h_3+h_4)v_{13}=-
(1+\lambda_2+\lambda_3+\lambda_4)v_{13}
\nonumber  \hspace{1.0cm}(5.11)
\end{eqnarray}
provides a contradiction.

Case 2.2 \quad $wt(\xi_{0,1})= \mbox{wt}(d_{15}v_{\lambda})$

In this case, $v_{15}=v_{\lambda}$ and $\mbox{wt}(v_{13})=\lambda-\alpha_{3}-\alpha_4$.
And $$0=E_{25}(E_{53}.v_{12}-E_{52}.v_{13})=-(\lambda_2+\lambda_3+\lambda_4)
v_{13}-E_{53}v_{15},$$
$$0=E_{35}E_{25}(E_{53}.v_{12}-E_{52}.v_{13})=-E_{35}(\lambda_2+\lambda_3+\lambda_4)v_{13}-E_{35}
E_{53}v_{15}=(\lambda_2+\lambda_3+\lambda_4-h_{3}-h_4)v_{15}
=\lambda_{2}v_{15} \eqno(5.12)$$
forces $\lambda_2=0$.

Case 2.2.1. \quad $\lambda_{3}= 0, \ \lambda_{4}> 0$.

Suppose $$T_{5,123}.\tilde{c}^{k}=
\sum\limits_{1 \leq i < j \leq 5}Q_{ij}^{k-1}(-1)^{|ij|}\partial_{y_{ij}}
.\eqno(5.13)$$By (5.3), $Q_{15}^{1}v_{\lambda}=(Q_{12}^{0}E_{52}+Q_{13}^{0}E_{53})v_{\lambda}=0$.
Then $$T_{5,123}.\xi_{0,1}=[\partial_{y_{12}}E_{53}+(-1)
^{1+|13|}\partial_{y_{13}}E_{52}+(-1)^{|23|}
\partial_{y_{23}}E_{51}].\tilde{c}.d_{15}v_{\lambda}=Q_{15}^{1}.v_{\lambda}=0.
\eqno(5.14)$$ Therefore,
$\lambda=(m,0,0,n), m \in \mathbb{N}, n > 0$.

 Case 2.2.2. \quad $\lambda_{3} > 0 , \ \lambda_{4}> 0$.

Note that $E_{54}v_{\lambda} \neq 0 $. Since $T_{5,123}.\tilde{c}.d_{15}v_{\lambda}=Q_{15}^{1}.v_{\lambda}=0$, the equation (3.8) implies $$T_{5,123}.\xi_{0,1}=T_{5,123}.\tilde{c}^{2}.d_{15}v_{\lambda}=
Q_{15}^{2}.v_{\lambda}
=(E_{53}E_{42}E_{54}-E_{52}E_{43}E_{54})v_{\lambda}=0;\eqno(5.15)$$ which yields
$$0=E_{34}E_{25}(E_{53}E_{42}E_{54}-E_{52}E_{43}E_{54})v_{\lambda}=-\lambda_3
(1+\lambda_3+\lambda_4)E_{54}v_{\lambda}.\eqno(5.16)$$ A contradiction arises.

Case 3.  \quad wt($\xi_{0,1})\in \{ \mbox{wt}(d_{23}v_{\lambda}), \mbox{wt}(d_{24}v_{\lambda}),\mbox{wt}(d_{34}v_{\lambda}), \mbox{wt}(d_{25}v_{\lambda}), \mbox{wt}(d_{35}v_{\lambda}), \mbox{wt}(d_{45}v_{\lambda})\}$.

In these cases, $v_{13}\neq 0, v_{23}\neq 0$.
Set
$$Q=E_{53}.v_{12}-E_{52}.v_{13}+E_{51}.v_{23}, \  E_{15}Q=Q_{1},  \ E_{25}Q_{1}=Q_{2}, E_{35}Q_{1}=Q_{2}'\ .\eqno(5.17)$$
 Then $$0=Q_{1}=E_{13}v_{12}+E_{53}E_{15}v_{12}-E_{12}v_{13}-E_{52}E_{15}v_{13}+
(h_1+h_2+h_3+h_4)v_{23}.\eqno(5.18)$$

Case 3.1. \quad $wt(\xi_{0,1})\in \{\mbox{wt}(d_{23}v_{\lambda}), \mbox{wt}(d_{24}v_{\lambda}), \mbox{wt}(d_{34}v_{\lambda})\}$.

In these three cases, we have $v_{23}\neq 0$, $v_{25}=v_{35}=0$.
So $$0=Q_{1}=(2+|\mbox{wt}(v_{23})|)v_{23}\eqno(5.19)$$ induces a contradiction.

Case 3.2. \quad $wt(\xi_{0,1})= \mbox{wt}(d_{25}v_{\lambda})$.

In this case, $v_{3 5}=0$ and $\mbox{wt}(v_{23})=\lambda-\alpha_3-\alpha_4 $. So $$0=Q_{1}=(2+|\mbox{wt}(v_{23})|)v_{23}+E_{53}v_{25},0=Q_{2}'=-(|\lambda|+2-h_{3}-h_4)v_{25}=-(\lambda_1+\lambda_2+1)v_{25}\eqno(5.20)$$
force $\lambda_1+\lambda_2+1=0$.
A contradiction arises.

Case 3.3.  \quad $wt(\xi_{0,1})\in \{\mbox{wt}(d_{35}v_{\lambda}),\mbox{wt}(d_{45}v_{\lambda}) \}$.

  In these two  cases, $v_{35}\neq 0$. And the equations $$0=Q_{1}=(2+|\mbox{wt}(v_{23})|)v_{23}+E_{53}v_{25}-E_{52}v_{35},
 0=Q_{2}=(1+|\mbox{wt}(v_{23})|-h_{2}-h_{3}-h_{4})v_{35}=\lambda_1
v_{35}\eqno(5.21)$$ imply $\lambda_1=0$.
Recall the intertwining  operator defined in Lemma 5.2.
In these two cases, $T_{5,123}(\xi_{0,1})=0$
is equivalent to $T^{1}|_{V(\omega_1+\omega_2)\otimes V(\mbox{wt}(\xi_{0,1}))}=0$.
 Assume $v_{(\omega_1+\omega_2)\otimes (\mbox{wt}(\xi_{0,1}))}^{\lambda}$ is any maximal vector of weight $\lambda$ appearing in the tensor decomposition
$V(\omega_1+\omega_2)\otimes V(\mbox{wt}(\xi_{0,1}))$. Then $T^{1}|_{V(\omega_1+\omega_2)\otimes V(\mbox{wt}(\xi_{0,1}))}=0$ iff
$T^{1}(v_{(\omega_1+\omega_2)\otimes (\mbox{wt}(\xi_{0,1}))}^{\lambda})=0$.

Case 3.3.1  \quad $wt(\xi_{0,1})=\mbox{wt}(d_{35}v_{\lambda})
=(0,\lambda_{2}-1,\lambda_{3}+1,\lambda_{4}-1)
$.

Indeed, $T^{1}(v_{(1,1,0,0)\otimes (0,\lambda_{2}-1,\lambda_{3}+1,\lambda_{4}-1)}^{\lambda})$ in this case could be written as:
\begin{eqnarray}
\nonumber&&T^{1}(v_{(1,1,0,0)\otimes (0,\lambda_{2}-1,\lambda_{3}+1,\lambda_{4}-1)}^{\lambda})\\
\nonumber
&=&(T_{4,345}+T_{2,235})\xi_{0,1}+\frac{1}{2}(T_{1,135}-T_{2,235})\xi_{0,1}
-\frac{1}{1+\lambda_3}(T_{3,345}-T_{2,245}).(x_{4}\partial_{x_{3}})_{0}.\xi_{0,1}
\\
\nonumber&-&\frac{1}{2+2\lambda_3}(T_{2,245}-T_{1,145}).(
x_{4}\partial_{x_{3}})_{0}.
\xi_{0,1}
-\frac{3}{1+\lambda_2+\lambda_3}T_{2,345}.x_{4}\partial_{x_{3}}.x_{3}
\partial_{x_{2}}\xi_{0,1}\\
\nonumber&+&\frac{6+3\lambda_3}{(1+\lambda_2+\lambda_3)
(1+\lambda_3)}T_{2,345}.(x_{3}
\partial_{x_{2}})_{0}.(x_{4}\partial_{x_{3}})_{0}.\xi_{0,1}
\\
\nonumber&-&
\frac{6+3\lambda_3}{(1+\lambda_2+\lambda_3)
(1+\lambda_3)}T_{1,345}.(x_{2}
\partial_{x_{1}})_{0}.(x_{3}
\partial_{x_{2}})_{0}.(x_{4}\partial_{x_{3}})_{0}.\xi_{0,1}\\
\nonumber&+&
\frac{3}{1+\lambda_2+\lambda_3}
T_{1,345}.(x_{4}\partial_{x_{3}})_{0}.(x_{2}
\partial_{x_{1}})_{0}.(x_{3}
\partial_{x_{2}})_{0}.\xi_{0,1}=\frac{(2+\lambda_3)(\lambda_2+\lambda_3+7)
)}{(1+\lambda_3)(\lambda_2+\lambda_3+1)
)}v_{35}\neq 0.
\nonumber \hspace{1.2cm}(5.22)
\end{eqnarray}

Case 3.3.2  \quad $wt(\xi_{0,1})=\mbox{wt}(d_{45}v_{\lambda})
=(0,\lambda_{2},\lambda_{3}-1,\lambda_{4})$.

Suppose  $\lambda_2\neq 0$.
Then
\begin{eqnarray}
\nonumber&&T^{1}(v_{(1,1,0,0)\otimes (0,\lambda_{2},\lambda_{3}-1,\lambda_{4})}^{\lambda})=
\frac{-2\lambda_2}{3}(T_{3,345}-T_{2,245})\xi_{0,1}-
\frac{\lambda_2}{3}(T_{2,245}-T_{1,145})\xi_{0,1}\nonumber\\&+&T_{2,345}.
(x_{3}
\partial_{x_{2}})_{0}.\xi_{0,1}-
T_{1,345}.(x_{2}
\partial_{x_{1}})_{0}.
(x_{3}
\partial_{x_{2}})_{0}.\xi_{0,1}=\frac{2\lambda_2(\lambda_2+3)}{3}v_{45}\neq 0
\nonumber \hspace{2.8cm}(5.23)
\end{eqnarray}
induces a contradiction.
Assume  $\lambda_2=0$.
 Then  it is easily checked that
 $$T^{1}(v_{(1,1,0,0)\otimes (0,\lambda_{2},\lambda_{3}-1,\lambda_{4})}^{\lambda})
=[2(T_{3,345}-T_{2,245})+(T_{2,245}-T_{1,145})]\xi_{0,1}=0.\eqno(5.24)$$
Thus $\lambda=(0,0,m,n)$.
The proof is complete by Lemma 3.3.
\hspace{1cm} $\Box$\vspace{0.3cm}

\subsection{Singular vectors of degree two }

{\it Theorem 5.4}\quad {All the possible  degree two  singular vectors are listed in the following:
$$\prod\limits_{\mbox{wt}(d_{12}d_{15}) < \sigma \leq \omega_1+\omega_3} \frac{\tilde{c}-\chi_{\sigma+\lambda}(\tilde{c})}
{\chi_{\mbox{wt}(d_{12}d_{15})+\lambda}(\tilde{c})-
\chi_{\sigma+\lambda}(\tilde{c})}.d_{12}d_{15}v_{\lambda},
 \ \mbox{where} \
 \  \lambda=(m,0,0,1), \ m \in \mathbb{N}. $$
}

{\it Proof} \quad The leading term of any singular vector of
degree two  could be written as:
   $$\xi_{0,2}=\sum\limits_{j \in \overline{1,35}, k \in  \overline{1,\mbox{mult}(\overrightarrow{w}_{j}^{\omega_1+\omega_3} )}}v_{j,k}^{\omega_1+\omega_3}
v_{j,k}^{\lambda}, \ v_{j,k}^{\lambda} \in V(\lambda), \eqno(5.27)$$
 which should satisfy  $T_{5,123}.\xi_{0,2}=0$. Assume
  $$T_{5,123}.\xi_{0,2}=
 \sum\limits_{1 \leq i  < j  \leq 5}d_{ij}t_{ij},  t_{ij} \in V(\lambda). \eqno(5.28)$$
Then we could derive the following equations:
$$t_{12}=E_{52}v_{1,1}^{\lambda}-E_{51}v_{2,1}^{\lambda}=0,
t_{15}=E_{53}v_{5,1}^{\lambda}-E_{52}v_{10,1}^{\lambda}+
E_{51}(v_{11,2}^{\lambda}+v_{11,3}^{\lambda})+v_{1,1}^{\lambda}=0,$$
$$
t_{13}=E_{53}v_{1,1}^{\lambda}-E_{51}v_{6,1}^{\lambda}=0,
t_{25}=E_{53}v_{9,1}^{\lambda}+
E_{52}(v_{11,1}^{\lambda}+v_{11,3}^{\lambda})
+E_{51}v_{19,1}^{\lambda}+
v_{2,1}^{\lambda}
=0,$$$$
t_{23}=E_{53}v_{2,1}^{\lambda}-E_{52}v_{6,1}^{\lambda}=0,
t_{35}=E_{53}(v_{11,1}^{\lambda}+v_{11,2}^{\lambda})-E_{52}v_{17,1}^{\lambda}
+E_{51}v_{22,1}^{\lambda}+v_{6,1}^{\lambda}
=0.\eqno(5.29)$$
It follows from $v_{1,1}^{\lambda}\neq 0$ that one of $v_{5,1}^{\lambda}, \ v_{10,1}^{\lambda}, \ v_{11,2}^{\lambda}+v_{11,3}^{\lambda} $ should be nonzero.
Hence, the information of the weights in Table 5 implies that $\mbox{wt}(\xi_{0,2})$ should be restricted to the following cases:
\begin{eqnarray}
&&\mbox{wt}(\xi_{0,2})\in \{ \lambda+ \overrightarrow{w}_{5}^{\omega_1+\omega_3},
\lambda+ \overrightarrow{w}_{9}^{\omega_1+\omega_3},
\lambda+ \overrightarrow{w}_{10}^{\omega_1+\omega_3},
\lambda+ \overrightarrow{w}_{11}^{\omega_1+\omega_3},
\lambda+ \overrightarrow{w}_{15}^{\omega_1+\omega_3},
\lambda+ \overrightarrow{w}_{16}^{\omega_1+\omega_3},\nonumber\\&&
\lambda+ \overrightarrow{w}_{17}^{\omega_1+\omega_3},
\lambda+ \overrightarrow{w}_{19}^{\omega_1+\omega_3},
\lambda+ \overrightarrow{w}_{21}^{\omega_1+\omega_3},
\lambda+ \overrightarrow{w}_{22}^{\omega_1+\omega_3},
\lambda+ \overrightarrow{w}_{i}^{\omega_1+\omega_3}(i \in \overline{24,35})
 \}\nonumber \hspace{2.6cm}(5.30)
\end{eqnarray}

Case 1  \quad $ \mbox{wt}(\xi_{0,2})\in \{ \lambda+ \overrightarrow{w}_{5}^{\omega_1+\omega_3},
\lambda+ \overrightarrow{w}_{10}^{\omega_1+\omega_3},
\lambda+ \overrightarrow{w}_{15}^{\omega_1+\omega_3}\}$

In these three cases,
$v_{2,1}^{\lambda}=v_{6,1}^{\lambda}=0$.

Case 1.1 \quad  $wt(\xi_{0,2})=\lambda+ \overrightarrow{w}_{5}^{\omega_1+\omega_3}$

 We have $wt(v_{1,1}^{\lambda})=\lambda-\alpha_{3}-\alpha_{4}$ and $\lambda_4 > 0$.
  So $0=E_{25}.t_{12}=E_{25}.E_{52}v_{1,1}^{\lambda}=(h_2+h_3+h_4)v_{1,1}^{\lambda}=
  (\lambda_2+\lambda_3+\lambda_{4}-1)v_{1,1}^{\lambda}$
 yields $(\lambda_2,\lambda_3,\lambda_4)=(0,0,1).$
  Then $\chi_{\overrightarrow{w}_{1}^{\omega_1+\omega_3}}
  (\tilde{c})=\frac{3\lambda_{1}+2}{25}$ by Lemma 3.5.
  And$$T_{5,123}.\xi_{0,2}=T_{5,123}.(\tilde{c}- \chi_{\overrightarrow{w}_{1}^{\omega_1+\omega_3}}(\tilde{c}))
.d_{12}d_{15}.v_{\lambda}=T_{5,123}.\tilde{c}.d_{12}d_{15}.v_{\lambda}
-\chi_{\overrightarrow{w}_{1}^{\omega_1+\omega_3}}(\tilde{c}
)d_{15}E_{53}.v_{\lambda}=0
,\eqno(5.31)$$
since$$T_{5,123}.\tilde{c}.d_{12}d_{15}.v_{\lambda}\stackrel{\mbox{by} \ (5.13)}{=}
Q_{12}^{1}d_{15}v_{\lambda}
-Q_{15}^{1}d_{12}v_{\lambda};$$where
$$Q_{12}^{1}\stackrel{\mbox{by} (5.3)}{=}Q_{12}^{0}(\tilde{c}+\frac{h_{2}^{*}}{10})
+\frac{1}{10}Q_{13}^{0}E_{23}-\frac{1}{10}Q_{23}^{0}E_{13}, \ Q_{15}^{1}\stackrel{\mbox{by} (5.3)}{=}\frac{1}{10}(Q_{12}^{0}E_{52}
+Q_{13}^{0}E_{53}),$$
$$Q_{12}^{0}= E_{53}+\frac{(x_5\partial_{x_3})_{0}'}{2},
\ Q_{13}^{0}=- E_{52}-\frac{(x_5\partial_{x_2})_{0}'}{2},Q_{23}^{0}= E_{51}+\frac{(x_5\partial_{x_1})_{0}'}{2},$$
$$\tilde{c}.d_{15}v_{\lambda}\stackrel{\mbox{by} (5.1)}{=}\frac{1}{10}[d_{15}(h_{1}^{*}-h_{4}^{*}).v_{\lambda}+
\sum\limits_{i=2}^{4}d_{1i}E_{5i}v_{\lambda}]=
\frac{1}{10}[\frac{3\lambda_{1}-3}{5}d_{15}v_{\lambda}+
\sum\limits_{i=2}^{4}d_{1i}E_{5i}v_{\lambda}].\eqno(5.32)
$$
That is, $\lambda=(m,0,0,1).$

Case 1.2 \quad  $wt(\xi_{0,2})=\lambda+ \overrightarrow{w}_{10}^{\omega_1+\omega_3}$

We have $\mbox{wt}(v_{1,1}^{\lambda})=\lambda-\alpha_{2}-\alpha_{3}-\alpha_{4}$ and $\lambda_2 > 0, \lambda_4 > 0$,
$E_{25}.v_{1,1}^{\lambda}=-v_{10,1}^{\lambda}$.
 So $$0=E_{25}t_{12}=(h_2+h_3+h_4)v_{1,1}^{\lambda}+E_{52}.E_{25}.v_{1,1}^{\lambda}=
  (\lambda_2+\lambda_3+\lambda_{4}-2)v_{1,1}^{\lambda}-E_{52}.v_{10,1}^{\lambda},
  \eqno(5.33)$$
   $$0=E_{25}^{2}t_{12}=E_{25}.[(\lambda_2+\lambda_3+\lambda_{4}-2)
  v_{1,1}^{\lambda}-E_{52}.v_{10,1}^{\lambda}]
  =-2(\lambda_2+\lambda_3+\lambda_{4}-1)v_{10,1}^{\lambda}\eqno(5.34)$$
yields a contradiction.

Case 1.3 \quad $wt(\xi_{0,2})=\lambda+ \overrightarrow{w}_{15}^{\omega_1+\omega_3}$

We have  $\mbox{wt}(v_{10,1}^{\lambda})=\lambda-\alpha_{3}$ and $\lambda_3 > 0$, $E_{23}.v_{1,5}^{\lambda}=-v_{10,1}^{\lambda}$.
Then
 $$E_{53}v_{5,1}^{\lambda}-E_{52}v_{10,1}^{\lambda}+v_{1,1}^{\lambda}=0, \eqno(5.35)$$
$$0=E_{25}(E_{53}v_{5,1}^{\lambda}-E_{52}v_{10,1}^{\lambda}+v_{1,1}^{\lambda})=
E_{23}v_{5,1}^{\lambda}-(h_2+h_3+h_4)
v_{10,1}^{\lambda}+E_{25}v_{1,1}^{\lambda}\eqno(5.36)$$
imply $(\lambda_2,\lambda_3,\lambda_4)=(0,0,0).$
A contradiction arises.

Case 2 \quad $\mbox{wt}(\xi_{0,2})\in \{
\lambda+ \overrightarrow{w}_{9}^{\omega_1+\omega_3},
\lambda+ \overrightarrow{w}_{11}^{\omega_1+\omega_3},
\lambda+ \overrightarrow{w}_{16}^{\omega_1+\omega_3},
\lambda+ \overrightarrow{w}_{17}^{\omega_1+\omega_3},
\lambda+ \overrightarrow{w}_{19}^{\omega_1+\omega_3},
\lambda+ \overrightarrow{w}_{21}^{\omega_1+\omega_3},
\lambda+ \overrightarrow{w}_{22}^{\omega_1+\omega_3},
\lambda+ \overrightarrow{w}_{i}^{\omega_1+\omega_3}(i \in \overline{24,35})
 \}$

Case 2.1 \quad $wt(\xi_{0,2})=\lambda+ \overrightarrow{w}_{9}^{\omega_1+\omega_3}$

We have $ v_{2,1}^{\lambda} \neq 0$, since $E_{35}v_{2,1}^{\lambda}=-v_{9,1}^{\lambda}$.
Note that   $v_{6,1}^{\lambda}=0$, $E_{15}v_{1,1}^{\lambda}=E_{15}v_{2,1}^{\lambda}=0$, $E_{12}v_{1,1}^{\lambda}=-v_{2,1}^{\lambda}$,  $\mbox{wt}(v_{2,1}^{\lambda})=\lambda-\alpha_{3}-\alpha_{4}$ and $\lambda_1 > 0, \lambda_4 > 0$. Then
   $$0=E_{15}.t_{12}=E_{15}.(E_{52}v_{1,1}^{\lambda}-E_{51}v_{2,1}^{\lambda})=
(E_{12}+E_{52}E_{15})v_{1,1}^{\lambda}
-(\sum\limits_{i=1}^{4}h_{i}+E_{51}E_{15})v_{2,1}^{\lambda} \eqno(5.37)$$
implies $|\lambda|=0$. A contradiction arises.

Case 2.2 \quad  $\mbox{wt}(\xi_{0,2})\in \{
\lambda+ \overrightarrow{w}_{11}^{\omega_1+\omega_3},
\lambda+ \overrightarrow{w}_{16}^{\omega_1+\omega_3},
\lambda+ \overrightarrow{w}_{17}^{\omega_1+\omega_3},
\lambda+ \overrightarrow{w}_{19}^{\omega_1+\omega_3},
\lambda+ \overrightarrow{w}_{21}^{\omega_1+\omega_3},
\lambda+ \overrightarrow{w}_{22}^{\omega_1+\omega_3},
\lambda+ \overrightarrow{w}_{i}^{\omega_1+\omega_3}(i \in \overline{24,35})
 \}$

Case 2.2.1 \quad  $\mbox{wt}(\xi_{0,2})\in \{
\lambda+ \overrightarrow{w}_{11}^{\omega_1+\omega_3},
\lambda+ \overrightarrow{w}_{16}^{\omega_1+\omega_3},
\lambda+ \overrightarrow{w}_{17}^{\omega_1+\omega_3},
\lambda+ \overrightarrow{w}_{19}^{\omega_1+\omega_3},
\lambda+ \overrightarrow{w}_{21}^{\omega_1+\omega_3},
\lambda+ \overrightarrow{w}_{22}^{\omega_1+\omega_3},
\lambda+ \overrightarrow{w}_{24}^{\omega_1+\omega_3},
\lambda+ \overrightarrow{w}_{26}^{\omega_1+\omega_3},
\lambda+ \overrightarrow{w}_{27}^{\omega_1+\omega_3},
\lambda+ \overrightarrow{w}_{29}^{\omega_1+\omega_3},
\lambda+ \overrightarrow{w}_{30}^{\omega_1+\omega_3},
\lambda+ \overrightarrow{w}_{33}^{\omega_1+\omega_3}.
 \}$

In these cases, $v_{25,1}^{\lambda}= v_{28,1}^{\lambda}= 0$.
Assume $v_{11,2}^{\lambda}+v_{11,3}^{\lambda}\neq 0$. Then
$$0=E_{15}t_{15}\stackrel{\mbox{by} (5.29)}{=}(2+\sum\limits_{i=1}^{4}h_{i})(v_{11,2}^{\lambda}+v_{11,3}^{\lambda})
=0\eqno(5.38)$$yields a contradiction. Hence, $v_{11,2}^{\lambda}+v_{11,3}^{\lambda}=0$.
Furthermore, either the assertion $v_{10,1}^{\lambda}\neq 0, \ (h_2+h_3+h_4)v_{10,1}^{\lambda}=0$ or the assertion
$v_{10,1}^{\lambda}=0, \ v_{5,1}^{\lambda}\neq 0, \ (h_3+h_4-1)v_{5,1}^{\lambda}=0$ holds.
By detailed check case by case, only  the cases $\mbox{wt}(\xi_{0,2})\in \{
\lambda+ \overrightarrow{w}_{16}^{\omega_1+\omega_3},
\lambda+ \overrightarrow{w}_{24}^{\omega_1+\omega_3},
\lambda+ \overrightarrow{w}_{26}^{\omega_1+\omega_3}\}$
 satisfy this assertion.

For the case $\mbox{wt}(\xi_{0,2})=
\lambda+ \overrightarrow{w}_{16}^{\omega_1+\omega_3}$, we get $
v_{10,1}^{\lambda}=0, \ v_{5,1}^{\lambda}\neq 0, \ (\lambda_3,\lambda_4)=(1,0)$. And we could write $\xi_{0,2}=y_{12}\xi_{0,2}'$, where $\xi_{0,2}'=\sum\limits_{1\leq i < j \leq 5}d_{ij}v_{ij}$. Hence,
$T_{5,123}.\xi_{0,2}=T_{5,123}.y_{12}\xi_{0,2}'=([T_{5,123},y_{12}]+
y_{12}T_{5,123}).\xi_{0,2}'=(x_{5}\partial_{x_{3}}+\frac{1}{2}(-1)^{|34,45|}
y_{45}\partial_{34}).\xi_{0,2}'=\frac{1}{2}(-1)^{|34,45|}
y_{45}\partial_{34}.\xi_{0,2}'\neq 0$.

For the case $\mbox{wt}(\xi_{0,2})=
\lambda+ \overrightarrow{w}_{24}^{\omega_1+\omega_3}$, we get $
v_{19,1}^{\lambda}\neq 0$. The equation $E_{15}t_{25}=0$ implies that
$$0=E_{15}v_{2,1}^{\lambda}+(E_{13}+E_{53}E_{15})v_{9,1}^{\lambda}
+(E_{12}+E_{52}E_{15})(v_{11,1}^{\lambda}+v_{11,3}^{\lambda})+\sum\limits_{i=1}^{4}
h_{i}v_{19,1}^{\lambda}=(3+|\lambda|)v_{19,1}^{\lambda}.\eqno(5.39)$$
For the case $\mbox{wt}(\xi_{0,2})=
\lambda+ \overrightarrow{w}_{26}^{\omega_1+\omega_3}$, one of $v_{19,1}^{\lambda} $ and
$v_{22,1}^{\lambda}$ should be nonzero, otherwise $v_{26,i}^{\lambda}=0$.
Then,
the equation $E_{15}t_{35}=0$ implies that
$$0=E_{15}v_{6,1}^{\lambda}+(E_{13}+E_{53}E_{15})
(v_{11,1}^{\lambda}+v_{11,2}^{\lambda})
-(E_{12}+E_{52}E_{15})v_{17,1}^{\lambda})+\sum\limits_{i=1}^{4}
h_{i}.v_{22,1}^{\lambda}=(1+|\lambda|)v_{22,1}^{\lambda}.\eqno(5.40)$$

Case 2.2.2 \quad  $\mbox{wt}(\xi_{0,2})=
\lambda+ \overrightarrow{w}_{25}^{\omega_1+\omega_3}$

In this case, $v_{25,1}^{\lambda}\neq 0$, $  v_{28,1}^{\lambda}= 0$ and
 $wt(\xi_{0,2})=\lambda+ \mbox{wt}(d_{15}d_{25})$. Then $$0=E_{15}t_{15}=(2+\sum\limits_{i=1}^{4}h_{i})(v_{11,2}^{\lambda}-v_{11,3}^{\lambda}
)-E_{53}v_{25,1}^{\lambda},$$
$$E_{35}.E_{15}t_{15}=(2+|\mbox{wt}(\overrightarrow{w}_{25}
^{\omega_1+\omega_3})|-h_{3}-h_{4})v_{25,1}^{\lambda}
=0.\eqno(5.41)$$
So $\lambda_1+\lambda_2+1=0$. A contradiction arises.

Case 2.2.3 \quad   $\mbox{wt}(\xi_{0,2})\in \{
\lambda+ \overrightarrow{w}_{28}^{\omega_1+\omega_3},
\lambda+ \overrightarrow{w}_{31}^{\omega_1+\omega_3},
\lambda+ \overrightarrow{w}_{32}^{\omega_1+\omega_3},
\lambda+ \overrightarrow{w}_{34}^{\omega_1+\omega_3},
\lambda+ \overrightarrow{w}_{35}^{\omega_1+\omega_3},
 \}$

In these cases, $v_{25,1}^{\lambda}\neq 0$ and $  v_{28,1}^{\lambda}\neq 0$.

Case 2.2.3.1  \quad $\mbox{wt}(\xi_{0,2})\in \{
\lambda+ \overrightarrow{w}_{28}^{\omega_1+\omega_3},
\lambda+ \overrightarrow{w}_{31}^{\omega_1+\omega_3}\}$

The equations $$0=E_{15}t_{12}=(E_{12}+E_{52}E_{15})v_{1,1}^{\lambda}-(\sum\limits_{i=1}^{4}
h_{i}+
E_{51}E_{15})v_{2,1}^{\lambda}
=-(1+|\mbox{wt}(\overrightarrow{w}_{2}^{\omega_1+\omega_3})|)
v_{1,1}^{\lambda}+E_{52}(v_{11,2}^{\lambda}+v_{11,3}^{\lambda}
),\eqno(5.42)$$
$$E_{25}^{2}E_{15}t_{12}=2(\lambda_{2}
+\lambda_3+\lambda_{4}-1)v_{28,1}^{\lambda}=0\eqno(5.43)$$
induce a contradiction.

Case 2.2.3.2  \quad $\mbox{wt}(\xi_{0,2})\in \{
\lambda+ \overrightarrow{w}_{32}^{\omega_1+\omega_3},
\lambda+ \overrightarrow{w}_{34}^{\omega_1+\omega_3},
\lambda+ \overrightarrow{w}_{35}^{\omega_1+\omega_3}\}$

First, we have $v_{32,1}^{\lambda}\neq 0$ in these cases, since $E_{34}
v_{32,1}^{\lambda}=-v_{34,1}^{\lambda},  \ E_{24}
v_{32,1}^{\lambda}=v_{35,1}^{\lambda} $. Then
\begin{eqnarray}
&&0=E_{15}t_{35}=E_{15}v_{6,1}^{\lambda}+
(E_{13}+E_{53}E_{15})(v_{11,1}^{\lambda}+v_{11,2}^{\lambda})
-(E_{12}+E_{52}E_{15})v_{17,1}^{\lambda}\nonumber \\&&+
(\sum\limits_{i=1}^{4}h_{i})v_{22,1}^{\lambda}
=(1+|\mbox{wt}(\overrightarrow{w}_{22}^{\omega_1+\omega_3})|)
v_{22,1}^{\lambda}+E_{53}v_{32,1}^{\lambda},\nonumber
\end{eqnarray}
$$E_{35}E_{15}t_{35}=(-1-|\mbox{wt}
(\overrightarrow{w}_{22}^{\omega_1+\omega_3})|+h_3+h_4)v_{32,1}
^{\lambda}=0.\eqno(5.44)$$
If $
\mbox{wt}(\xi_{0,2})=
\lambda+ \overrightarrow{w}_{32}^{\omega_1+\omega_3}$,
then $\mbox{wt}(\overrightarrow{w}_{32}^{\omega_1+\omega_3})=\lambda$ and
$\mbox{wt}(\overrightarrow{w}_{22}^{\omega_1+\omega_3})
=\lambda-\alpha_3-\alpha_4$. If $
\mbox{wt}(\xi_{0,2})=
\lambda+ \overrightarrow{w}_{34}^{\omega_1+\omega_3}$,
then $\mbox{wt}(\overrightarrow{w}_{32}^{\omega_1+\omega_3})=\lambda-\alpha_3$ and $\mbox{wt}(\overrightarrow{w}_{22}^{\omega_1+\omega_3})
=\lambda-2\alpha_3-\alpha_4$. If $\mbox{wt}(\xi_{0,2})=
\lambda+ \overrightarrow{w}_{35}^{\omega_1+\omega_3}$,
then $\mbox{wt}(\overrightarrow{w}_{32}^{\omega_1+\omega_3})=\lambda-\alpha_2-\alpha_3$
and $\mbox{wt}(\overrightarrow{w}_{22}^{\omega_1+\omega_3})
=\lambda-\alpha_{2}-2\alpha_3-\alpha_4$. Thus, (5.44) yields  $\lambda_1+\lambda_2=0$ or $\lambda_1+\lambda_2+1=0$. A contradiction arises.
\hspace{1cm} $\Box$\vspace{0.1cm}

\subsection{Singular vectors of degree three }

{\it Theorem 5.5}\quad {All the possible  degree three  singular vectors  are listed in the following:
$$ e^{-\frac{P}{2}}.
\prod\limits_{\mbox{wt}(d_{15} d_{25}d_{45}) < \sigma \leq 2\omega_1+\omega_4} \frac{\tilde{c}-\chi_{\sigma+\lambda}(\tilde{c})}
{\chi_{\mbox{wt}(d_{15} d_{25} d_{45})+\lambda}(\tilde{c})-
\chi_{\sigma+\lambda}(\tilde{c})}.d_{15} d_{25} d_{45}v_{\lambda},
 $$$$\ \mbox{where} \
\ \lambda=(0,0,m,n),  \ 1 \leq m \in \mathbb{N}, \  2 \leq n \in \mathbb{N} .  $$}

{\it Proof} \quad The leading term of any singular vector of degree three
could be written as:
 $$\xi_{0,3}=\sum\limits_{j \in \overline{1,55}, k \in  \overline{1,\mbox{mult}(\overrightarrow{w}_{j}^{2\omega_1+\omega_4} )} }v_{j,k}^{2\omega_1+\omega_4}v_{j,k}^{\lambda}, \ v_{j,k}^{\lambda} \in V(\lambda),\eqno(5.45)$$which should satisfy $T_{5,123}.\xi_{0,3}=0$ and $(x_i\partial_{x_j})_{0}.\xi_{0,3}=0$ ( $1 \leq i < j \leq 5$).
 Assume $$T_{5,123}.\xi_{0,3}=\sum\limits_{i_1 < j_1,i_2 < j_2,(i_1,j_1)\neq  (i_2,j_2)}d_{i_1j_1}d_{i_2j_2}t_{i_1j_1,i_2j_2}, \ t_{i_1j_1,i_2j_2} \in V(\lambda).\eqno(5.46)$$
  Since $v_{1,1}^{\lambda}\neq 0 $, the equation $$0=t_{14,15}=-v_{1,1}^{\lambda}+E_{53}v_{7,1}^{\lambda}
  -E_{52}v_{13,1}^{\lambda}
   +E_{51}(v_{18,1}^{\lambda}-v_{18,4}^{\lambda})\eqno(5.47)$$
 implies  that one of the terms $v_{7,1}^{\lambda}$, $v_{13,1}^{\lambda}$ and
 $v_{18,1}^{\lambda}-v_{18,4}^{\lambda}$ should be nonzero. Hence, the information of the weights of  Table 6-7 induces that
  $wt(\xi_{0,3})$
could be restricted to the following cases:
\begin{eqnarray}&&wt(\xi_{0,3})\in \{\lambda+ \overrightarrow{w}_{7}^{2\omega_1+\omega_4},
\lambda+ \overrightarrow{w}_{12}^{2\omega_1+\omega_4},
\lambda+ \overrightarrow{w}_{13}^{2\omega_1+\omega_4},
\lambda+ \overrightarrow{w}_{17}^{2\omega_1+\omega_4},
\lambda+ \overrightarrow{w}_{18}^{2\omega_1+\omega_4},\nonumber\\&&
\lambda+ \overrightarrow{w}_{21}^{2\omega_1+\omega_4},
\lambda+ \overrightarrow{w}_{22}^{2\omega_1+\omega_4},
\lambda+ \overrightarrow{w}_{23}^{2\omega_1+\omega_4},
\lambda+ \overrightarrow{w}_{24}^{2\omega_1+\omega_4},
\lambda+ \overrightarrow{w}_{i}^{2\omega_1+\omega_4}
(i \in \overline{26,55})
   \}.\nonumber \hspace{1.8cm}(5.48)
   \end{eqnarray}

Case 1 \quad $wt(\xi_{0,3})=\lambda+ \overrightarrow{w}_{7}^{2\omega_1+\omega_4}$\vspace{0.1cm}

In this case, $v_{6,1}^{\lambda}=0$, $v_{3,1}^{\lambda} \neq 0$, $wt(v_{3,1}^{\lambda})=\lambda-\alpha_3$, $E_{34}v_{3,1}^{\lambda}=-v_{7,1}^{\lambda}$. Then
$$E_{25}t_{12,15}=E_{25}(E_{52}v_{3,1}^{\lambda}-E_{51}v_{6,1}^{\lambda})
=(h_2+h_3+h_4)v_{3,1}^{\lambda}=(\lambda_2+\lambda_3+\lambda_4)v_{3,1}^{\lambda}=0
\eqno(5.49)$$ contradicts $\lambda_3 > 0$.\vspace{0.1cm}

Case 2 \quad $wt(\xi_{0,3})\in \{\lambda+ \overrightarrow{w}_{12}^{2\omega_1+\omega_4},\lambda+ \overrightarrow{w}_{17}^{2\omega_1+\omega_4}\}$\vspace{0.1cm}

Note  that
$E_{34}v_{6,1}^{\lambda}=-v_{12,1}^{\lambda},
E_{12}v_{12,1}^{\lambda}=-v_{17,1}^{\lambda},
E_{12}v_{3,1}^{\lambda}=-2v_{6,1}^{\lambda}.$
Hence, $v_{3,1}^{\lambda} \neq 0$ and $v_{6,1}^{\lambda},\neq 0$.
Since $wt(v_{6,1}^{\lambda})=
\lambda-\alpha_3$ for $wt(\xi_{0,3})=\lambda+ \overrightarrow{w}_{12}^{2\omega_1+\omega_4}$;
$wt(v_{6,1}^{\lambda})=
\lambda-\alpha_1-\alpha_3 $ for $
wt(\xi_{0,3})=\lambda+ \overrightarrow{w}_{17}^{2\omega_1+\omega_4}$.
The equation $E_{15}t_{12,15}=-(2+|\mbox{wt}v_{6,1}^{\lambda}|)v_{6,1}^{\lambda}=0$
induces that $ |\lambda| < 0$ in both cases.\vspace{0.1cm}

Case 3 \quad $ wt(\xi_{0,3})=\lambda+ \overrightarrow{w}_{13}^{2\omega_1+\omega_4}\vspace{0.1cm}
$

We have $E_{25}v_{1,1}^{\lambda}=-v_{13,1}^{\lambda}$, $
E_{23}v_{7,1}^{\lambda}=-v_{13,1}^{\lambda}$. Then
$$0=E_{25}t_{14,15}\stackrel{\mbox{by} (5.47)}{=}-E_{25}v_{1,1}^{\lambda}+E_{23}v_{7,1}^{\lambda}-\sum\limits_{i=2}^{4}h_i
v_{13,1}^{\lambda}=-(\lambda_2+\lambda_3+\lambda_4)v_{13,1}^{\lambda}=0,
\eqno(5.50)$$
 contradicts $\lambda_2 > 0$.\vspace{0.1cm}

 Case 4 \quad $wt(\xi_{0,3})=\lambda+
  \overrightarrow{w}_{22}^{2\omega_1+\omega_4}$\vspace{0.1cm}

We have $v_{11,1}^{\lambda}\neq 0$, since $E_{24}v_{11,1}^{\lambda}=v_{22,1}^{\lambda}$. Then  $$
0=E_{15}t_{13,15}=E_{15}(E_{53}v_{3,1}^{\lambda}-E_{51}v_{11,1}^{\lambda})
=-(2+|\lambda|)v_{11,1}^{\lambda}\eqno(5.51)$$yields a contradiction.
\vspace{0.1cm}

 Case 5  \quad $wt(\xi_{0,3})\in \{\lambda+ \overrightarrow{w}_{24}^{2\omega_1+\omega_4},
\lambda+ \overrightarrow{w}_{28}^{2\omega_1+\omega_4},\lambda+
\overrightarrow{w}_{30}^{2\omega_1+\omega_4},\lambda+
\overrightarrow{w}_{35}^{2\omega_1+\omega_4},\lambda+
\overrightarrow{w}_{36}^{2\omega_1+\omega_4}
 \}$
\vspace{0.1cm}

For $wt(\xi_{0,3})\in \{\lambda+ \overrightarrow{w}_{24}^{2\omega_1+\omega_4},
\lambda+ \overrightarrow{w}_{30}^{2\omega_1+\omega_4},\lambda+
\overrightarrow{w}_{36}^{2\omega_1+\omega_4}\}$, we have
$v_{28,1}^{\lambda}=0$.
Then $$0=E_{25}t_{15,25}=
E_{25}(E_{52}v_{24,1}^{\lambda}
      -E_{51}v_{28,1}^{\lambda})=
(h_2+h_3+h_4)v_{24,1}^{\lambda}\eqno(5.52)$$ yields
$\lambda_2+\lambda_3+\lambda_4=0$ or $1$.
 For $wt(\xi_{0,3})\in \{\lambda+ \overrightarrow{w}_{28}^{2\omega_1+\omega_4},
\lambda+ \overrightarrow{w}_{35}^{2\omega_1+\omega_4}\},$ the equation
$$0=E_{15}t_{15,25}=-(h_1+h_2+h_3+h_4)v_{28,1}^{\lambda}+(E_{12}+E_{52}E_{15})v_{24
,1}^{\lambda}\eqno(5.53)$$ yields
$1+|\lambda|=0$.\vspace{0.1cm}

 Case 6 \quad $wt(\xi_{0,3})\in \{\lambda+
  \overrightarrow{w}_{23}^{2\omega_1+\omega_4}, \lambda+
  \overrightarrow{w}_{29}^{2\omega_1+\omega_4}\}$\vspace{0.1cm}

We have $E_{35}v_{10,1}^{\lambda}=v_{23,1}^{\lambda}, E_{25}v_{10,1}^{\lambda}=-v_{29,1}^{\lambda},
E_{23}v_{23,1}^{\lambda}=-v_{29,1}^{\lambda}$. Consider the equation
 $$t_{14,45}=-v_{10,1}^{\lambda}+E_{53}v_{23,1}^{\lambda}
   -E_{52}v_{29,1}^{\lambda}-E_{51}(v_{32,1}^{\lambda}-v_{32,4}^{\lambda}).
   \eqno(5.54)
   $$
For $wt(\xi_{0,3})=\lambda+
  \overrightarrow{w}_{23}^{2\omega_1+\omega_4}$, $E_{35}t_{14,45}=(\lambda_3+\lambda_4-1)v_{23,1}^{\lambda}=0$ induces $\lambda_3+\lambda_4=1$, which contradicts $\lambda_3 > 1$.
For $wt(\xi_{0,3})=\lambda+
  \overrightarrow{w}_{29}^{2\omega_1+\omega_4}$, we have $E_{25}t_{14,45}=-(\lambda_2+\lambda_3+\lambda_4)v_{29,1}^{\lambda}=0$ induces $\lambda_2+\lambda_3+\lambda_4=0$, which contradicts $\lambda_2 > 0, \lambda_3 > 0 $.\vspace{0.1cm}

 Case 7 \quad $wt(\xi_{0,3})\in \{\lambda+
  \overrightarrow{w}_{27}^{2\omega_1+\omega_4}, \lambda+
  \overrightarrow{w}_{39}^{2\omega_1+\omega_4},   \lambda+
  \overrightarrow{w}_{42}^{2\omega_1+\omega_4}  \}$\vspace{0.1cm}

Note that in these cases, the equations are derived: $E_{35}v_{15,1}^{\lambda}=v_{27,1}^{\lambda}, E_{15}v_{15,1}^{\lambda}=v_{39,1}^{\lambda},
E_{13}v_{27,1}^{\lambda}=v_{39,1}^{\lambda}, E_{12}v_{32,4}^{\lambda}=-v_{39,1}^{\lambda},
 E_{12}v_{32,1}^{\lambda}=0=
E_{12}v_{32,3}^{\lambda},E_{25}v_{15,1}^{\lambda}=
-\sum\limits_{i=1}^{4}v_{32,i}^{\lambda},
  E_{23}v_{27,1}^{\lambda}=-v_{32,1}^{\lambda}-2v_{32,3}^{\lambda}
  -v_{32,4}^{\lambda}, E_{34}v_{32,1}^{\lambda}=-v_{42,1}^{\lambda},
  E_{34}v_{32,2}^{\lambda}=v_{42,1}^{\lambda},E_{34}v_{32,3}^{\lambda}=0=
  E_{34}v_{32,4}^{\lambda}
 .$
Consider the equation
$$t_{24,45}=-v_{15,1}^{\lambda}+E_{53}v_{27,1}^{\lambda}
   -E_{52}(v_{32,1}^{\lambda}+v_{32,3}^{\lambda}+v_{32,4}^{\lambda})
   +E_{51}v_{39,1}^{\lambda}.\eqno(5.55)
   $$
For $wt(\xi_{0,3})=\lambda+
  \overrightarrow{w}_{27}^{2\omega_1+\omega_4}$, $E_{35}t_{24,45}=(\lambda_3+\lambda_4-1)v_{27,1}^{\lambda}=0$ induces $\lambda_3+\lambda_4=1$, which contradicts $\lambda_3 > 1$.
For $wt(\xi_{0,3})=\lambda+
  \overrightarrow{w}_{39}^{2\omega_1+\omega_4}$, $E_{15}t_{24,45}=(1+|\lambda|)v_{39,1}^{\lambda}=0$ induces $1+|\lambda|=0$.
  For $wt(\xi_{0,3})=\lambda+
  \overrightarrow{w}_{42}^{2\omega_1+\omega_4}$,
   Then
  $E_{15}t_{24,45}=v_{32,2}^{\lambda}-v_{32,3}^{\lambda}-
  (\lambda_2+\lambda_3+\lambda_4)(v_{32,1}^{\lambda}+
  v_{32,3}^{\lambda}+v_{32,4}^{\lambda})=0$, and $E_{34}E_{15}t_{24,45}=
    (1+\lambda_2+\lambda_3+\lambda_4)v_{42,1}^{\lambda}=0$ induces $1+\lambda_2+\lambda_3+\lambda_4=0$.
\vspace{0.1cm}

 Case 8 \quad $wt(\xi_{0,3})\in \{\lambda+
  \overrightarrow{w}_{31}^{2\omega_1+\omega_4}, \lambda+
  \overrightarrow{w}_{33}^{2\omega_1+\omega_4} \}$\vspace{0.1cm}

Observe these equations are derived:
 $E_{15}v_{4,1}^{\lambda}=v_{31,1}^{\lambda}, E_{13}v_{17,1}^{\lambda}=v_{31,1}^{\lambda},
E_{12}v_{21,1}^{\lambda}=-v_{31,1}^{\lambda},
E_{12}v_{21,2}^{\lambda}=v_{31,1}^{\lambda},
 E_{12}v_{21,3}^{\lambda}=0,
E_{12}v_{21,4}^{\lambda}=0.$
Consider the equation
 $$t_{24,25}=-v_{4,1}^{\lambda}+E_{53}v_{17,1}^{\lambda}
   -E_{52}(2v_{21,1}^{\lambda}+v_{21,2}^{\lambda}+v_{21,3}^{\lambda})
       +E_{51}v_{31,1}^{\lambda}.\eqno(5.56)
   $$
It follows from $E_{15}t_{24,25}=0$ that $(1+|\lambda|) v_{31,1}^{\lambda}=0$.
\vspace{0.1cm}

 Case 9 \quad $wt(\xi_{0,3})\in \{\lambda+
  \overrightarrow{w}_{37}^{2\omega_1+\omega_4}, \lambda+
  \overrightarrow{w}_{38}^{2\omega_1+\omega_4} \}$\vspace{0.1cm}

Note
 $E_{15}v_{14,1}^{\lambda}=v_{38,1}^{\lambda}, E_{13}v_{26,2}^{\lambda}=v_{38,1}^{\lambda},
E_{12}v_{37,1}^{\lambda}=-v_{38,1}^{\lambda},
E_{25}v_{14,1}^{\lambda}=-v_{37,1}^{\lambda},
 E_{23}v_{26,2}^{\lambda}=-v_{37,1}^{\lambda}
.$
Consider the equation
 $$t_{34,35}=-v_{14,1}^{\lambda}+E_{53}v_{26,2}^{\lambda}
      -E_{52}v_{37,1}^{\lambda}
       +E_{51}v_{38,1}^{\lambda}.\eqno(5.57)
   $$
For $wt(\xi_{0,3})=\lambda+
  \overrightarrow{w}_{37}^{2\omega_1+\omega_4}$, $E_{25}t_{34,35}=-
 ( \lambda_2+\lambda_3+\lambda_4)v_{37,1}^{\lambda}=0$ implies
 $ \lambda_2+\lambda_3+\lambda_4=0$, which
 which contradicts $\lambda_2 > 2$.
  For $wt(\xi_{0,3})=\lambda+
  \overrightarrow{w}_{38}^{2\omega_1+\omega_4}$,
  $E_{15}t_{34,35}=
 (1+|\lambda|)v_{38,1}^{\lambda}=0$.
\vspace{0.1cm}

Case 10 \quad $wt(\xi_{0,3})\in \{\lambda+
  \overrightarrow{w}_{40}^{2\omega_1+\omega_4}, \lambda+
  \overrightarrow{w}_{44}^{2\omega_1+\omega_4},
\overrightarrow{w}_{48}^{2\omega_1+\omega_4},\overrightarrow{w}_{51}^{2\omega_1+\omega_4}
   \}$\vspace{0.1cm}

Note
 $E_{25}v_{19,1}^{\lambda}=-v_{40,1}^{\lambda}, E_{23}v_{32,3}^{\lambda}=-v_{40,1}^{\lambda}, E_{15}v_{19,1}^{\lambda}=v_{44,1}^{\lambda}, E_{13}v_{32,3}^{\lambda}=v_{44,1}^{\lambda}, E_{12}v_{40,1}^{\lambda}=-v_{44,1}^{\lambda}.$
Consider the equation
$$t_{34,45}=-v_{19,1}^{\lambda}+E_{53}v_{32,3}^{\lambda}
      -E_{52}v_{40,1}^{\lambda}
       +E_{51}v_{44,1}^{\lambda},\eqno(5.58)
   $$
For $wt(\xi_{0,3})=\lambda+
  \overrightarrow{w}_{40}^{2\omega_1+\omega_4}, \lambda+
  \overrightarrow{w}_{48}^{2\omega_1+\omega_4}, $
 the equation $E_{25}t_{34,45}=0$ implies $-(\lambda_2+\lambda_3+\lambda_4)v_{40,1}^{\lambda}=0$.
  For $wt(\xi_{0,3})=\lambda+
  \overrightarrow{w}_{44}^{2\omega_1+\omega_4}, \lambda+
  \overrightarrow{w}_{51}^{2\omega_1+\omega_4}, $
 the equation $E_{15}t_{34,45}=0$ implies $(1+|\lambda|)v_{44,1}^{\lambda}=0$.
\vspace{0.1cm}

Case 11 \quad $wt(\xi_{0,3})\in \{\lambda+
  \overrightarrow{w}_{45}^{2\omega_1+\omega_4}, \lambda+
  \overrightarrow{w}_{52}^{2\omega_1+\omega_4}
   \}$\vspace{0.1cm}

Note
$E_{15}v_{21,3}^{\lambda}=-v_{45,1}^{\lambda},
 E_{15}v_{21,4}^{\lambda}=v_{45,1}^{\lambda},
 E_{13}v_{35,1}^{\lambda}=v_{45,1}^{\lambda}, E_{12}v_{41,1}^{\lambda}=-v_{45,1}^{\lambda},
  E_{12}v_{41,i}^{\lambda}=0 \ (i \in \overline{2,4}), E_{15}v_{18,3}^{\lambda}=v_{41,1}^{\lambda}-v_{41,4}^{\lambda},
E_{15}v_{18,4}^{\lambda}=-v_{41,1}^{\lambda}-v_{41,3}^{\lambda},
E_{13}v_{30,1}^{\lambda}=v_{41,1}^{\lambda}-v_{41,2}^{\lambda}-
v_{41,3}^{\lambda},
E_{12}v_{36,1}^{\lambda}=-2v_{41,1}^{\lambda}-v_{41,2}^{\lambda},
E_{15}v_{30,1}^{\lambda}=E_{35}v_{41,3}^{\lambda}=
-E_{35}v_{41,4}^{\lambda}=-v_{52,1}^{\lambda}
.$
Consider the equation
$$t_{25,45}=-v_{21,3}^{\lambda}-v_{21,4}^{\lambda}+E_{53}v_{35,1}^{\lambda}-
   E_{52}(v_{41,1}^{\lambda}+v_{41,2}^{\lambda}+v_{41,3}^{\lambda})
  +
   E_{51}v_{45,1}^{\lambda}=0.\eqno(5.59)
   $$
For $wt(\xi_{0,3})=\lambda+
  \overrightarrow{w}_{45}^{2\omega_1+\omega_4}, $
 the equation $E_{15}t_{25,45}=0$ implies $(2+|\lambda|)v_{45,1}^{\lambda}=0$.
Now suppose $wt(\xi_{0,3})=\lambda+
  \overrightarrow{w}_{52}^{2\omega_1+\omega_4}. $
The equation $t_{15,45}=0$ implies that
  $$E_{15}t_{15,45}=0=(|\lambda|+2)v_{41,1}^{\lambda}+(1-|\lambda|
  )v_{41,3}^{\lambda}+v_{41,4}^{\lambda}-E_{53}v_{52,1}^{\lambda},$$
 $$0=E_{35}E_{15}t_{15,45}=(\lambda_1+\lambda_2)v_{52,1}^{\lambda},\eqno(5.60)$$
i.e. $\lambda_1=\lambda_2=0$. Then $wt(\xi_{0,3})=\lambda+
  \overrightarrow{w}_{52}^{2\omega_1+\omega_4}=(0,1,\lambda_{3}-1,\lambda_{4}-2). $
Observe that
$$\Pi((\omega_1+\omega_2)\otimes \mbox{wt}(\xi_{0,3}))\bigcap \Pi(
(\omega_{1}+\omega_{3})\otimes \lambda
)=\{\mbox{wt}(\xi_{0,3})+(0,0,1,0), \mbox{wt}(\xi_{0,3})+(1,-1,0,1)\}.\eqno(5.61)$$
$$\Pi((\omega_1+\omega_2)\otimes (\mbox{wt}(\xi_{0,3})+(0,0,1,0)))\bigcap \{\lambda\}=\emptyset, \ \Pi((\omega_1+\omega_2)\otimes (\mbox{wt}(\xi_{0,3})+(1,-1,0,1)))\bigcap \{\lambda\}=\emptyset .\eqno(5.62)$$
Then (5.62) implies that
$$T^{2}|_{V(\omega_1+\omega_2)\otimes V(\mbox{wt}(\xi_{0,3})+(0,0,1,0))}=0,
T^{2}|_{V(\omega_1+\omega_2)\otimes V(\mbox{wt}(\xi_{0,3})+(1,-1,0,1))}=0.\eqno(5.63)$$
Note that $T_{5,123}\xi_{0,3}=0$ iff
$T^{3}|_{V(\omega_1+\omega_2)\otimes V(\mbox{wt}(\xi_{0,3}))}=0$ by Lemma 5.2.
Assume $v_{(\omega_1+\omega_2)\otimes \mbox{wt}(\xi_{0,3})}^{(0,1,\lambda_{3},\lambda_{4}-2)}$
(resp. $v_{(\omega_1+\omega_2)\otimes \mbox{wt}(\xi_{0,3})}^{(1,0,\lambda_{3}-1,\lambda_{4}-1)}$ ) is any maximal vector of weight $(0,1,\lambda_{3},\lambda_{4}-2) $
(resp. $(1,0,\lambda_{3}-1,\lambda_{4}-1) $ )
appearing in the tensor decomposition
$V(\omega_1+\omega_2)\otimes V(\mbox{wt}(\xi_{0,3}))$.
Since $T^{3}{V(\omega_1+\omega_2)\otimes V(\mbox{wt}(\xi_{0,3}))}\subseteq
V(\omega_1+\omega_3)\otimes V(\lambda)$,
$T^{3}|_{V(\omega_1+\omega_2)\otimes V(\mbox{wt}(\xi_{0,3}))}=0$ iff
$T^{3}.(v_{(\omega_1+\omega_2)\otimes \mbox{wt}(\xi_{0,3})}^{(0,1,\lambda_{3},\lambda_{4}-2)})=0,
T^{3}(v_{(\omega_1+\omega_2)\otimes \mbox{wt}(\xi_{0,3})}^{(1,0,\lambda_{3}-1,\lambda_{4}-1)})=0$.
Otherwise, we could get a singular vector of degree two with weights
$(0,1,\lambda_{3},\lambda_{4}-2) $
and  $(1,0,\lambda_{3}-1,\lambda_{4}-1) $  respectively by (5.63),
which contradicts proposition 5.3.
Hence, we get $\lambda=(0,0,m,n)$.\vspace{0.1cm}

Case 12  \quad $wt(\xi_{0,3})\in \{\lambda+
  \overrightarrow{w}_{43}^{2\omega_1+\omega_4}, \lambda+
  \overrightarrow{w}_{47}^{2\omega_1+\omega_4}
   \}$\vspace{0.1cm}

Consider the equation
$$t_{15,45}=-v_{18,3}^{\lambda}-v_{18,4}^{\lambda}+E_{53}v_{30,1}^{\lambda}
   -
   E_{52}v_{36,1}^{\lambda}
   +
   E_{51}v_{41,1}^{\lambda}-E_{51}v_{41,3}^{\lambda},\eqno(5.64)
   $$Note
$E_{15}v_{18,3}^{\lambda}=v_{41,1}^{\lambda}-v_{41,4}^{\lambda},
E_{15}v_{18,4}^{\lambda}=-v_{41,1}^{\lambda}-v_{41,3}^{\lambda},
 E_{13}v_{30,1}^{\lambda}=
v_{41,1}^{\lambda}-v_{41,2}^{\lambda}-v_{41,3}^{\lambda},
 E_{12}v_{36,1}^{\lambda}=
 -2v_{41,1}^{\lambda}-v_{41,2}^{\lambda}.$
Then $E_{15}t_{15,45}=0$ induces that
$3v_{41,1}^{\lambda}+v_{41,4}^{\lambda}+\sum\limits_{i=1}^{4}h_i(v_{41,1}^{\lambda}-
v_{41,3}^{\lambda})=0$. Since $E_{45}v_{41,4}^{\lambda}=-v_{43,1}^{\lambda},
E_{45}v_{41,i}^{\lambda}=0 (i \in \overline{1,3}) $ and $E_{34}v_{41,3}^{\lambda}=-v_{47,1}^{\lambda}, E_{34}v_{41,i}^{\lambda}=0 (i \in \{1,2,4\})$. Hence, $E_{45}E_{15}t_{15,45}=0$ and $E_{34}E_{15}t_{15,45}=0$ imply $v_{41,3}^{\lambda}=0$ and $ |\lambda|v_{47,1}^{\lambda}=0$ in these two cases, respectively.
\vspace{0.1cm}

Case 13 \quad $wt(\xi_{0,3})\in \{\lambda+
  \overrightarrow{w}_{46}^{2\omega_1+\omega_4}, \lambda+
  \overrightarrow{w}_{49}^{2\omega_1+\omega_4},
\lambda+
  \overrightarrow{w}_{50}^{2\omega_1+\omega_4}, \lambda+
  \overrightarrow{w}_{53}^{2\omega_1+\omega_4},\lambda+
  \overrightarrow{w}_{54}^{2\omega_1+\omega_4}, \lambda+
  \overrightarrow{w}_{55}^{2\omega_1+\omega_4}
   \}$\vspace{0.1cm}

Observe
$E_{25}v_{26,3}^{\lambda}=2v_{46,1}^{\lambda},
 E_{25}v_{26,4}^{\lambda}=-v_{46,1}^{\lambda},
 E_{23}v_{41,2}^{\lambda}=-v_{46,1}^{\lambda},
 E_{15}v_{26,3}^{\lambda}=-2v_{49,1}^{\lambda},
 E_{15}v_{26,4}^{\lambda}=v_{49,1}^{\lambda},
 E_{13}v_{41,2}^{\lambda}=v_{49,1}^{\lambda}, E_{12}v_{46,1}^{\lambda}=-v_{49,1}^{\lambda},
 E_{25}v_{41,2}^{\lambda}=-v_{54,1}^{\lambda},
 E_{35}v_{46,1}^{\lambda}=-v_{54,1}^{\lambda}
.$
Consider the equation
  $$t_{35,45}=-v_{26,3}^{\lambda}-v_{26,4}^{\lambda}+E_{53}v_{41,2}^{\lambda}
    -
   E_{52}v_{46,1}^{\lambda}+
   E_{51}v_{49,1}^{\lambda},\eqno(5.65)
   $$
For $wt(\xi_{0,3})=\lambda+
  \overrightarrow{w}_{46}^{2\omega_1+\omega_4}, \overrightarrow{w}_{50}^{2\omega_1+\omega_4},$
 the equation $E_{25}t_{35,45}=0$ implies $
 2+ \lambda_2+\lambda_3+\lambda_4=0$.
For $wt(\xi_{0,3})=\lambda+
  \overrightarrow{w}_{49}^{2\omega_1+\omega_4},
  \overrightarrow{w}_{53}^{2\omega_1+\omega_4}, $
 the equation $E_{15}t_{35,45}=0$ implies $3+ |\lambda|=0$.
 For $wt(\xi_{0,3})=\lambda+
  \overrightarrow{w}_{54}^{2\omega_1+\omega_4}, $
 the equation $E_{35}E_{25}t_{35,45}=0$ implies $(1+ \lambda_2)v_{54,1}^{\lambda}=0$.
 For $wt(\xi_{0,3})=\lambda+
  \overrightarrow{w}_{55}^{2\omega_1+\omega_4}, $
 the equation $E_{35}E_{15}t_{35,45}=0$ implies $\lambda_1+ \lambda_2+2=0$.
\vspace{0.1cm}

Case 14  \quad  $wt(\xi_{0,3})=\lambda+ \overrightarrow{w}_{21}^{2\omega_1+\omega_4}$\vspace{0.3cm}

We have
$v_{21,3}^{\lambda}+2v_{21,4}^{\lambda}+E_{45}v_{15,1}^{\lambda}=0$.
And the equations $t_{25,45}=0$ and $t_{24,45}=0$ induce that $ v_{21,3}^{\lambda}+v_{21,4}^{\lambda} =0 $, $v_{15,1}^{\lambda}=0$.
Thus, $v_{21,3}^{\lambda}=v_{21,4}^{\lambda} =0 $.
Then
$$0=E_{25}.t_{24,25}=E_{25}.[-v_{4,1}^{\lambda}+E_{53}v_{17,1}^{\lambda}
   -E_{52}(2v_{21,1}^{\lambda}+v_{21,2}^{\lambda})]
      =0,\eqno(5.66)
   $$
$$0=E_{15}.t_{14,25}=E_{15}[-v_{2,1}^{\lambda}+E_{53}v_{12,1}^{\lambda}
   -E_{52}(v_{18,1}^{\lambda}+v_{18,2}^{\lambda}-v_{18,3}^{\lambda})+
   E_{51}v_{21,1}^{\lambda}]=0\eqno(5.67)$$
  yield
   $$2(1+\lambda_2+\lambda_3+\lambda_4)v_{21,1}^{\lambda}
   +(2+\lambda_2+\lambda_3+\lambda_4)v_{21,2}^{\lambda}=0, \
   (2+|\lambda|)v_{21,1}^{\lambda}=0 \eqno(5.68)$$i.e.
$v_{21,i}^{\lambda}=0$ for $i \in \overline{1,4}$.
A contradiction arises.\vspace{0.1cm}

Case 15  \quad  $wt(\xi_{0,3})=\lambda+ \overrightarrow{w}_{26}^{2\omega_1+\omega_4}$\vspace{0.1cm}

Note that $$t_{24,35}=-v_{8,1}^{\lambda}+E_{53}\sum\limits_{i=1}^{4}
          v_{21,i}^{\lambda}-E_{52}(v_{26,2}^{\lambda}+\sum\limits_{i=1}^{4}
          v_{26,i}^{\lambda}).\eqno(5.69)
                   $$
We have
$v_{26,2}^{\lambda}+v_{26,3}^{\lambda}+2v_{26,4}^{\lambda}
+E_{45}v_{19,1}^{\lambda}=0$,
$v_{26,2}^{\lambda}+v_{26,3}^{\lambda}
+E_{23}v_{21,3}^{\lambda}=0 $ and $v_{26,4}^{\lambda}
+E_{23}v_{21,4}^{\lambda}=0 $.
And the equations $t_{35,45}=0$ , $t_{34,45}=0$ and $t_{25,45}=0$ induce that $ v_{26,3}^{\lambda}+v_{26,4}^{\lambda} =0 $, $v_{19,1}^{\lambda}=0$ and
$v_{26,2}^{\lambda}+v_{26,3}^{\lambda}+v_{26,4}^{\lambda} =0 $. Thus, $v_{26,i}^{\lambda} =0 $ for $i \in \overline{2,4}$.
The equation $t_{24,35}=0$  induce that
$$E_{25}.[-v_{8,1}^{\lambda}
   +E_{53}(v_{21,1}^{\lambda}+v_{21,2}^{\lambda})-E_{52}v_{26,1}^{\lambda}]
      =0,\eqno(5.70)
   $$
  which means $\lambda_2+\lambda_3+\lambda_4=0$.
A contradiction arises.
\vspace{0.1cm}

Case 16  \quad  $wt(\xi_{0,3})=\lambda+ \overrightarrow{w}_{32}^{2\omega_1+\omega_4}$\vspace{0.1cm}

Note that
 $v_{32,1}^{\lambda}+v_{32,2}^{\lambda}
+E_{14}v_{18,3}^{\lambda}=0, v_{32,1}^{\lambda}+v_{32,4}^{\lambda}
+E_{14}v_{18,4}^{\lambda}=0, v_{32,1}^{\lambda}-2v_{32,3}^{\lambda}-v_{32,4}^{\lambda}
+E_{24}v_{21,3}^{\lambda}=0, v_{32,1}^{\lambda}+v_{32,2}^{\lambda}
+E_{24}v_{21,4}^{\lambda}=0 , 2v_{32,1}^{\lambda}
+E_{34}v_{26,3}^{\lambda}=0 , v_{32,1}^{\lambda}+v_{32,2}^{\lambda}
+E_{34}v_{26,4}^{\lambda}=0 .$
And the equations follow:
$$ t_{25,45} =v_{21,3}^{\lambda}+v_{21,4}^{\lambda} =0 , t_{35,45}=v_{26,3}^{\lambda}+v_{26,4}^{\lambda} =0 ,$$$$t_{34,45}=-v_{19,1}^{\lambda}+E_{53}
v_{32,3}^{\lambda}=0, t_{15,45}
=-v_{18,3}^{\lambda}-v_{18,4}^{\lambda}+E_{53}v_{30,1}^{\lambda} =0 .\eqno(5.71)$$
Thus,
$$-2v_{32,1}^{\lambda}+2v_{32,3}^{\lambda}-v_{32,2}^{\lambda}
+v_{32,4}^{\lambda}=0, \ 3v_{32,1}^{\lambda}+v_{32,2}^{\lambda}=0,
v_{32,2}^{\lambda}+(\lambda_3+\lambda_4-1)v_{32,3}^{\lambda}=0,
2v_{32,1}^{\lambda}+v_{32,2}^{\lambda}
+v_{32,4}^{\lambda}=0.\eqno(5.72)$$
 yield the contradiction:$2+\lambda_3+\lambda_4=0$.
\vspace{0.1cm}

Case 17  \quad  $wt(\xi_{0,3})=\lambda+ \overrightarrow{w}_{41}^{2\omega_1+\omega_4}$\vspace{0.1cm}

Observe that
$-v_{41,1}^{\lambda}-v_{41,2}^{\lambda}+v_{41,3}^{\lambda}+v_{41,4}^{\lambda}
+E_{24}v_{28,1}^{\lambda}=0 ,
v_{41,2}^{\lambda}+2v_{41,3}^{\lambda}+v_{41,4}^{\lambda}
+E_{34}v_{34,1}^{\lambda}=0 ,-v_{41,1}^{\lambda}-v_{41,2}^{\lambda}+v_{41,3}^{\lambda}
+E_{25}v_{21,3}^{\lambda}=0, v_{41,1}^{\lambda}+v_{41,2}^{\lambda}+v_{41,4}^{\lambda}
+E_{25}v_{21,4}^{\lambda}=0,
v_{41,1}^{\lambda}+2v_{41,2}^{\lambda}+v_{41,3}^{\lambda}
+E_{23}v_{35,1}^{\lambda}=0,
2v_{41,3}^{\lambda}
+E_{35}v_{26,3}^{\lambda}=0, -v_{41,2}^{\lambda}-v_{41,3}^{\lambda}+v_{41,4}^{\lambda}
+E_{35}v_{26,4}^{\lambda}=0.$
Then
$$0=E_{25}t_{25,45}=E_{25}[-v_{21,3}^{\lambda}-v_{21,4}^{\lambda}+E_{53}v_{35,1}^{\lambda}-
   E_{52}(v_{41,1}^{\lambda}+v_{41,2}^{\lambda}+v_{41,3}^{\lambda})]=0,\eqno(5.73)
     $$
$$0=E_{35}t_{35,45}=E_{35}[-v_{26,3}^{\lambda}-v_{26,4}^{\lambda}+E_{53}v_{41,2}^{\lambda}]=0.
    \eqno(5.74)  $$which mean
      $$v_{41,4}^{\lambda}-v_{41,1}^{\lambda}-2v_{41,2}^{\lambda}
      =(\sum\limits_{i=2}^{4}\lambda_i)
      (\sum\limits_{i=1}^{3}v_{41,i}^{\lambda}) , \ v_{41,3}^{\lambda}-v_{41,2}^{\lambda}+v_{41,4}^{\lambda}
      +(\lambda_3+\lambda_4)
      v_{41,2}^{\lambda}).\eqno(5.75)$$
      We claim that one of $v_{28,1}^{\lambda}$ and $ v_{34,1}^{\lambda}$ should be nonzero. Otherwise, the equation implies:
   $ -v_{41,1}^{\lambda}-v_{41,2}^{\lambda}+v_{41,3}^{\lambda}+v_{41,4}^{\lambda} =0$ and $ v_{41,2}^{\lambda}+2v_{41,3}^{\lambda}+v_{41,4}^{\lambda}=0$, which provide a contradiction.
Now,
 $0=E_{15}t_{15,35}=E_{15}( E_{53}v_{24,1}^{\lambda}
      -E_{51}v_{34,1}^{\lambda})=E_{13}v_{24,1}^{\lambda}-\sum\limits_{i=1}^
      {4}h_{i}v_{34,1}^{\lambda}=0$ and
$0=E_{15}t_{15,25}=E_{15}(E_{52}v_{24,1}^{\lambda}
      -E_{51}v_{28,1}^{\lambda})=E_{12}v_{24,1}^{\lambda}-\sum\limits_{i=1}^
      {4}h_{i}v_{34,1}^{\lambda}=0$.
      Thus $1+|\lambda|=0$. A contradiction arises.\vspace{0.1cm}

Case 18  \quad  $wt(\xi_{0,3})=\lambda+ \overrightarrow{w}_{18}^{2\omega_1+\omega_4}$\vspace{0.1cm}

Note that
 $E_{15}v_{1,1}^{\lambda}=v_{18,1}^{\lambda}-2v_{18,3}^{\lambda},
 E_{25}v_{2,1}^{\lambda}=-v_{18,1}^{\lambda}-v_{18,2}^{\lambda}
 -v_{18,3}^{\lambda},
 E_{13}v_{7,1}^{\lambda}=v_{18,1}^{\lambda}-
 2v_{18,2}^{\lambda}-2v_{18,4}^{\lambda}, E_{13}v_{3,1}^{\lambda}=-2v_{11,1}^{\lambda},
 E_{23}v_{12,1}^{\lambda}=-v_{18,1}^{\lambda}-2v_{18,2}^{\lambda}
 -v_{18,4}^{\lambda},
 E_{34}v_{11,1}^{\lambda}=-v_{18,2}^{\lambda}-
 v_{18,3}^{\lambda}-2v_{18,4}^{\lambda}
.$ We have $v_{11,1}^{\lambda}= 0$. Otherwise,
$0=E_{15}t_{13,15}=-(2+|\lambda|)v_{11,1}^{\lambda}=0$.
Hence,
$$0=-v_{18,2}^{\lambda}-
 v_{18,3}^{\lambda}-2v_{18,4}^{\lambda}=E_{34}v_{11,1}^{\lambda}
.\eqno(5.76)$$
Moreover,
 $$0=t_{15,45}=v_{18,3}^{\lambda}+v_{18,4}^{\lambda},
 0=E_{15}t_{14,15}=3v_{18,1}^{\lambda}+2v_{18,3}^{\lambda}-2
 v_{18,4}^{\lambda}+|\lambda|(v_{18,1}^{\lambda}-v_{18,4}^{\lambda}).\eqno(5.77)$$
Observe that$$t_{14,25}=-v_{2,1}^{\lambda}+E_{53}
          v_{12,1}^{\lambda}-E_{52}(v_{18,1}^{\lambda}+v_{18,2}^{\lambda}-
          v_{18,3}^{\lambda})+E_{51}(v_{21,1}^{\lambda}-v_{21,3}^{\lambda}-
          v_{21,4}^{\lambda})=0.\eqno(5.78)
                   $$
                   The equation $E_{25}t_{14,25}= 0$ induces
                   $$-v_{18,2}^{\lambda}+
 v_{18,3}^{\lambda}-v_{18,4}^{\lambda}-(\lambda_2+\lambda_3+\lambda_4)
 (v_{18,1}^{\lambda}+
 v_{18,2}^{\lambda}-v_{18,3}^{\lambda})=0.\eqno(5.79)$$
 All these equations yield the contradiction: $$\frac{3+|\lambda|}{4+|\lambda|}=-(\lambda_2+\lambda_3+\lambda_4).\eqno(5.80)$$
  \hspace{1cm} $\Box$\vspace{0.1cm}

\subsection{Singular vectors of degree four }

{\it Theorem 5.6}\quad {All the possible degree four  singular vectors  are listed in the following:
$$ d_{12}d_{13}d_{14}d_{15}v_{\lambda}, \ \lambda=(m,0,0,0), \ m \in \mathbb{N}.  $$}
{\it Proof} \quad The leading term of any singular vector of degree four
could be written as:
 $$\xi_{0,4}=\sum\limits_{j \in \overline{1,35}}v_{j}^{3\omega_1}
v_{j}^{\lambda},\eqno(5.81)$$which should satisfy $T_{5,123}.\xi_{0,4}=0$.
Since  $$[T_{5,123}-(E_{53}\partial_{y_{12}}+E_{52}(-1)^{1+|13|}\partial_{y_{13}}
 +E_{51}(-1)^{|23|}\partial_{y_{23}}),
 (x_{i}\partial_{x_{j}})_{0}']=0\eqno(5.82)$$ for $ 1 \leq j < i \leq 5$,
we have $$[T_{5,123}-(E_{53}\partial_{y_{12}}+E_{52}(-1)^{1+|13|}\partial_{y_{13}}
 +E_{51}(-1)^{|23|}\partial_{y_{23}})
 ].|_{ V(3\omega_1)\otimes V(\lambda)}=0
 .\eqno(5.83)$$Hence, $$T_{5,123}.\xi_{0,4}=[E_{53}\partial_{y_{12}}
 +E_{52}(-1)^{1+|13|}\partial_{y_{13}}
 +E_{51}(-1)^{|23|}\partial_{y_{23}}].\xi_{0,4}=0.\eqno(5.84)$$

Case 1 \quad $wt(\xi_{0,4})=\overrightarrow{w}_{1}^{3\omega_1}= \mbox{wt}(d_{12}d_{13}d_{14}d_{15}v_{\lambda})$\vspace{0.1cm}

The vector $d_{12}d_{13}d_{14}d_{15}v_{\lambda}$ is singular iff $\lambda=(m,0,0,0)$.\vspace{0.1cm}

Case 2  \quad $wt(\xi_{0,4})\in \{\lambda+ \overrightarrow{w}_{2}^{3\omega_1},
\lambda+ \overrightarrow{w}_{3}^{3\omega_1},\lambda+ \overrightarrow{w}_{5}^{3\omega_1},\lambda+ \overrightarrow{w}_{6}^{3\omega_1},
\lambda+ \overrightarrow{w}_{8}^{3\omega_1},\lambda+ \overrightarrow{w}_{10}^{3\omega_1},\lambda+ \overrightarrow{w}_{12}^{3\omega_1},
\lambda+ \overrightarrow{w}_{13}^{3\omega_1},\lambda+ \overrightarrow{w}_{15}^{3\omega_1},
\lambda+ \overrightarrow{w}_{16}^{3\omega_1},\lambda+ \overrightarrow{w}_{17}^{3\omega_1},
\lambda+ \overrightarrow{w}_{22}^{3\omega_1},
\lambda+ \overrightarrow{w}_{23}^{3\omega_1},\lambda+ \overrightarrow{w}_{27}^{3\omega_1},\lambda+ \overrightarrow{w}_{30}^{3\omega_1}
 \}$
\vspace{0.1cm}

In these cases, the equations are derived:
$E_{51}v_{i}^{\lambda}-E_{52}v_{j}^{\lambda}=0, \ E_{12}v_{j}^{\lambda}=-s_{ij}v_{i}^{\lambda}, \ \mbox{wt}( v_{i}^{\lambda})=\lambda, (i,j,s_{ij})\in \{(2,1,3),(3,2,2), (5,3,1), (8,6,1), (10,7,2),(12,9,1), (13,10,1),(15,11,2),(16,14,1),(17,15,1)\}.$
Then, the equation
$E_{15}.(E_{51}v_{i}^{\lambda}-E_{52}v_{j}^{\lambda})=(|\lambda|+s_{ij}
)v_{i}^{\lambda}=0$ yields
a contradiction.\vspace{0.1cm}

Case 3
 \quad $wt(\xi_{0,4})\in \{\lambda+ \overrightarrow{w}_{4}^{3\omega_1},
\lambda+ \overrightarrow{w}_{9}^{3\omega_1},\lambda+ \overrightarrow{w}_{14}^{3\omega_1},\lambda+ \overrightarrow{w}_{19}^{3\omega_1},
\lambda+ \overrightarrow{w}_{20}^{3\omega_1},\lambda+ \overrightarrow{w}_{21}^{3\omega_1},\lambda+ \overrightarrow{w}_{25}^{3\omega_1},
\lambda+ \overrightarrow{w}_{26}^{3\omega_1},\lambda+ \overrightarrow{w}_{29}^{3\omega_1},
\lambda+ \overrightarrow{w}_{33}^{3\omega_1}
 \}$
\vspace{0.1cm}

In these cases, we derive that:
$E_{51}v_{i}^{\lambda}-E_{53}v_{j}^{\lambda}=0, \ E_{13}v_{j}^{\lambda}=-t_{ij}v_{i}^{\lambda}, \ \mbox{wt}( v_{i}^{\lambda})=\lambda, (i,j,t_{ij})\in \{(4,1,3),(9,4,2), (14,7,2), (19,11,2), (20,9,1),(21,14,1), (25,18,1),(26,19,1),(29,24,1),(33,28,1)\}.
$Then,
$E_{15}.(E_{51}v_{i}^{\lambda}-E_{53}v_{j}^{\lambda})=(|\lambda|+t_{ij}
)v_{i}^{\lambda}=0$ yields
a contradiction.\vspace{0.1cm}

Case 4   \quad $wt(\xi_{0,4})\in \{\lambda+ \overrightarrow{w}_{7}^{3\omega_1},
\lambda+ \overrightarrow{w}_{11}^{3\omega_1},\lambda+ \overrightarrow{w}_{18}^{3\omega_1},\lambda+ \overrightarrow{w}_{24}^{3\omega_1},
\lambda+ \overrightarrow{w}_{28}^{3\omega_1}
 \}$\vspace{0.1cm}

In these cases, we have $\lambda_3+\lambda_4 > 0$. And
$E_{52}v_{\lambda}=E_{53}v_{\lambda}=0$, i.e. $\lambda=(m,0,0,0)$.   A contradiction arises.\vspace{0.1cm}

Case 5   \quad $wt(\xi_{0,4})\in \{\lambda+ \overrightarrow{w}_{31}^{3\omega_1},
\lambda+ \overrightarrow{w}_{32}^{3\omega_1},\lambda+ \overrightarrow{w}_{34}^{3\omega_1},\lambda+ \overrightarrow{w}_{35}^{3\omega_1}
 \}$\vspace{0.1cm}

Case 5.1  \quad $wt(\xi_{0,4})\in \{\lambda+ \overrightarrow{w}_{31}^{3\omega_1},
\lambda+ \overrightarrow{w}_{32}^{3\omega_1},\lambda+ \overrightarrow{w}_{34}^{3\omega_1}\}$\vspace{0.1cm}

In these cases, the following equations are derived:
$E_{51}v_{i}^{\lambda}-E_{53}v_{j}^{\lambda}=0, \ E_{13}v_{j}^{\lambda}=-q_{ij}v_{i}^{\lambda}, \ \mbox{wt}( v_{i}^{\lambda})=\lambda-\alpha_3,  \mbox{wt}( v_{j}^{\lambda})=\lambda-\alpha_{1}-\alpha_{2}-\alpha_{3}, (i,j,t_{ij})\in \{(25,18,1),(30,24,1), (34,28,1).$ Then,
$E_{15}.(E_{51}v_{i}^{\lambda}-E_{53}v_{j}^{\lambda})=(|\lambda|+q_{ij}
)v_{i}^{\lambda}=0$ yields
a contradiction.\vspace{0.1cm}

Case 5.2  \quad $wt(\xi_{0,4})=\lambda+ \overrightarrow{w}_{35}^{3\omega_1}$
\vspace{0.1cm}

Note that
$E_{51}v_{33}^{\lambda}-E_{53}v_{28}^{\lambda}=0,
 \ E_{51}v_{30}^{\lambda}-E_{52}v_{28}^{\lambda}=0,
\ \mbox{wt}( v_{33}^{\lambda})=\lambda-\alpha_3-\alpha_4, \
 \mbox{wt}(v_{30}^{\lambda})=\lambda-\alpha_{2}-\alpha_{3}-2\alpha_{4},
 \ \mbox{wt}( v_{28}^{\lambda})=\lambda-\alpha_1-\alpha_2-\alpha_3-\alpha_4, E_{13}v_{28}^{\lambda}=-v_{33}^{\lambda},
 E_{15}v_{28}^{\lambda}=-v_{35}^{\lambda},
 E_{35}v_{33}^{\lambda}=-v_{35}^{\lambda},
 E_{12}v_{28}^{\lambda}=-v_{30}^{\lambda},
 E_{25}v_{30}^{\lambda}=-v_{35}^{\lambda}.$
 Then $$E_{35}E_{15}(E_{51}v_{33}^{\lambda}-E_{53}v_{28}^{\lambda})=0=E_{35}(
 |\lambda|v_{33}^{\lambda}+E_{53}v_{35}^{\lambda})=-(\lambda_1+\lambda_2)
 v_{35}^{\lambda}\eqno(5.85)
 $$
 induces $\lambda_1=\lambda_2=0$.
 And
 $$E_{35}E_{15}(E_{51}v_{30}^{\lambda}-E_{52}v_{28}^{\lambda})=0=E_{35}[
 (|\lambda|-1)v_{30}^{\lambda}+E_{52}v_{35}^{\lambda}]=(1-\lambda_1)
 v_{35}^{\lambda}\eqno(5.86)
 $$
 induces $\lambda_1=1$.
A contradiction arises.
\hspace{1cm} $\Box$\vspace{0.5cm}

{\bf Acknowledgement:} Part of this work was completed while the author was visiting the Massachusetts Institute of
Technology in 2011. Se is very
grateful to all of the faculty and the staff members of the institution for their hospitality
and support during her visit.
The author would like to thank Professor Victor Kac
for his introducing  this interesting conjecture that motivates this work.
 The author would also like to thank Professor Dihua Jiang for his encouragements. Finally,  the author acknowledges the  financial support by CSC.

\begin{appendix}

\begin{table}
\caption{\bf Weights and weight vectors  for $\mbox{sl}_5$ module $V(\omega_2)$}
\begin{center}
\begin{tabular}{|c|p{4.5cm}|c|c|p{4.5cm}|c|}
\hline
$ i  $ &   $\overrightarrow{w}_{i}^{\omega_2}$  &
  $v_{i}^{\omega_2}$ & $ i  $ &   $\overrightarrow{w}_{i}^{\omega_2}$  &
  $v_{i}^{\omega_2}$ \\
\hline
1 & $(0,1,0,0)=\omega_2$ & $d_{12}$ & 6 & $(0,-1,0,1)=\omega_2-\alpha_1-2\alpha_2-\alpha_3$ & $d_{34}$
 \\
\hline
2 & $(1,-1,1,0)=\omega_2-\alpha_2$ & $d_{13}$ &
7 & $(1,0,0,-1)=\omega_2-\alpha_2-\alpha_3-\alpha_4$ & $d_{15}$
 \\
\hline
3 & $(-1,0,1,0)=\omega_2-\alpha_1-\alpha_2$ & $d_{23}$ & 8 & $(-1,1,0,-1)=\omega_2-\alpha_1-\alpha_2-\alpha_3-\alpha_4$ & $d_{25}$
 \\
\hline
4 & $(1,0,-1,1)=\omega_2-\alpha_2-\alpha_3$ & $d_{14}$  &
9 & $(0,-1,1,-1)=\omega_2-\alpha_1-2\alpha_2-\alpha_3-\alpha_4$ & $d_{35}$
\\
\hline
5 & $(-1,1,-1,1)=\omega_2-\alpha_1-\alpha_2-\alpha_3$ & $d_{24}$
& 10 & $(0,0,-1,0)=\omega_2-\alpha_1-2\alpha_2-2\alpha_3-\alpha_4$ & $d_{45}$ \\
\hline
\end{tabular}
\end{center}
\end{table}
\begin{table}
\caption{\bf Decomposition for wedge module $\Lambda^{k}W$ ($k \in \overline{1,10}$)}
\begin{center}
\begin{tabular}{|c|p{5.2cm}|p{7.2cm}|}
\hline
 $\Lambda^{k}W$ & irreducible components  for  $\Lambda^{k}W$ & maximal vector for the irreducible components
 \\
 \hline
$\wedge^{1}W$ & $V(\omega_2)$ & $d_{12}$
\\
\hline
$\wedge^{2}W$ & $V(\omega_1+\omega_3)$ &  ${ d_{12}}\wedge { d_{13}}$
\\
\hline
$ \wedge^{3}W$ & $V(2\omega_3)\oplus V(2\omega_1+\omega_4)$ & $
{ d_{12}}\wedge { d_{13}}\wedge { d_{23}},
{ d_{12}}\wedge { d_{13}}\wedge { d_{14}}$
\\
\hline
$ \wedge^{4}W$ & $V(3\omega_1)\oplus V(\omega_1+\omega_3+\omega_4)$ &
$
{ d_{12}}\wedge { d_{13}}\wedge { d_{14}}\wedge { d_{15}},
{ d_{12}}\wedge { d_{13}}\wedge { d_{23}}\wedge { d_{14}}$
 \\
\hline
$\wedge^{5}W$ & $V(2\omega_1+\omega_3)\oplus V(\omega_2+2\omega_4)$ &
${ d_{12}}\wedge { d_{13}}\wedge { d_{23}}\wedge { d_{14}}\wedge { d_{15}},
{ d_{12}}\wedge { d_{13}}\wedge { d_{23}}\wedge { d_{14}}\wedge { d_{24}}$
\\
\hline
$\wedge^{6}W$ & $V(3\omega_4)\oplus V(\omega_1+\omega_2+\omega_4)$ & $
{ d_{12}}\wedge { d_{13}}\wedge { d_{23}}\wedge { d_{14}}\wedge { d_{24}}\wedge { d_{34}}, { d_{12}}\wedge { d_{13}}\wedge {d_{23}}\wedge { d_{14}}\wedge { d_{24}}\wedge { d_{15}}$ \\
\hline
$\wedge^{7}W$ & $V(\omega_1+2\omega_4)\oplus V(2\omega_2)$ &
 ${ d_{13}}\wedge { d_{23}}\wedge { d_{14}}\wedge { d_{24}}\wedge {d_{34}}\wedge { d_{15}}, { d_{12}}\wedge { d_{13}}\wedge { d_{23}}\wedge { d_{14}}\wedge { d_{24}}\wedge { d_{15}}\wedge { d_{25}}$
 \\
 \hline
 $\wedge^{8}W$ & $V(\omega_2+\omega_4)$ & $ { d_{12}}\wedge { d_{13}}\wedge { d_{23}}\wedge { d_{14}}\wedge { d_{24}}\wedge { d_{34}}\wedge { d_{15}}\wedge { d_{25}}$
\\
\hline
$\wedge^{9}W$ & $V(\omega_3)$ & ${ d_{12}}\wedge { d_{13}}\wedge { d_{23}}\wedge { d_{14}}\wedge { d_{24}}\wedge { d_{34}}\wedge { d_{15}}\wedge { d_{25}}\wedge { d_{35}}$
\\
\hline
$\wedge^{10}W$ & $V(0)$ & ${ d_{12}}\wedge { d_{13}}\wedge { d_{23}}\wedge { d_{14}}\wedge { d_{24}}\wedge { d_{34}}\wedge { d_{15}}\wedge { d_{25}}\wedge { d_{35}}\wedge { d_{45}}$
\\
\hline
\end{tabular}
\end{center}
\end{table}
\begin{table}
\caption{\bf Tensor decomposition for $V(k\omega_4)\otimes V(\mu) $}
\begin{center}
\begin{tabular}{|c|p{13.2cm}|}
\hline
$\mu$ & highest  weight  in the decomposition $V(k\omega_4)\otimes V(\mu) $ \\
\hline
$ \omega_2$ &
$\omega_{(0,1,0,k)}, \omega_{(1,0,0,k-1)}$  \\
 \hline
 $ \omega_1+\omega_3$ &
$\omega_{(1,0,1,k)}, \omega_{(0,0,1,k-1)},\omega_{(1,1,0,k-1)},
\omega_{(0,1,0,k-2)}$ \\
\hline
$ 2\omega_3$ & $
\omega_{(0,0,2,k)},\omega_{ (0,1,1,k-1)},\omega_{ (0,2,0,k-2)}$
\\
\hline
$  2\omega_1+\omega_4$ & $
\omega_{(2,0,0,k+1)},\omega_{ (1,0,0,k)},\omega_{ (0,0,0,k-1)}, \omega_{(2,0,1,k-1)},\omega_{ (1,0,1,k-2)},\omega_{(0,0,1,k-3)}$
\\
\hline
$\omega_1+\omega_3+\omega_4$ &
$ \omega_{(1,0,1,k+1)},\omega_{ (0,0,1,k)},\omega_{ (1,1,0,k)}, \omega_{(0,1,0,k-1)},\omega_{ (1,0,2,k-1)},\omega_{ (0,0,2,k-2)}, \omega_{(1,1,1,k-2)},\omega_{ (0,1,1,k-3)}$ \\
\hline
$ 3\omega_1$ & $
\omega_{(3,0,0,k)},\omega_{ (2,0,0,k-1)},\omega_{ (1,0,0,k-2)},
\omega_{(0,0,0,k-3)}$ \\
\hline
$\omega_2+2\omega_4
$ & $ \omega_{(0,1,0,k+2)},\omega_{ (1,0,0,k+1)},\omega_{ (0,1,1,k)}, \omega_{(1,0,1,k-1)},\omega_{ (0,1,2,k-2)},\omega_{ (1,0,2,k-3)}$ \\
\hline
$2\omega_1+\omega_3
$ & $ \omega_{(2,0,1,k)},\omega_{(1,0,1,k-1)},\omega_{ (0,0,1,k-2)}, \omega_{(2,1,0,k-1)},\omega_{ (1,1,0,k-2)},\omega_{(0,1,0,k-3)}$
\\
\hline
$ 3\omega_4
$ & $ \omega_{(0,0,0,k+3)},\omega_{ (0,0,1,k+1)},\omega_{ (0,0,2,k-1)},
\omega_{ (0,0,3,k-3)}$ \\
\hline
$\omega_1+\omega_2+\omega_4$ &
$
\omega_{ (1,1,0,k+1)},\omega_{ (0,1,0,k)},\omega_{(2,0,0,k)},\omega_{ (1,0,0,k-1)},\omega_{ (1,1,1,k-1) },\omega_{ (0,1,1,k-2)},
 \omega_{(2,0,1,k-2)},\omega_{ (1,0,1,k-3)}$ \\
 \hline
 $\omega_1+2\omega_4
$ & $ \omega_{(1,0,0,k+2)},\omega_{ (0,0,0,k+1)},\omega_{ (1,0,1,k)}, \omega_{(0,0,1,k-1)},\omega_{(1,0,2,k-2)
},\omega_{ (0,0,2,k-3)}$ \\
\hline $
 2\omega_2$ & $ \omega_{(0,2,0,k)},\omega_{(1,1,0,k-1)},\omega_{(2,0,0,k-2)}$ \\
\hline
$\omega_2+\omega_4
$ & $ \omega_{(0,1,0,k+1)},\omega_{ (1,0,0,k)}, \omega_{(0,1,1,k-1)},\omega_{ (1,0,1,k-2)},$ \\
\hline $
 \omega_3
$ & $\omega_{(0,0,1,k)},\omega_{ (0,1,0,k-1)}$
\\
\hline
\end{tabular}
\end{center}
\end{table}
\begin{table}
\caption{\bf The leading term $l_{\mu}$}
\begin{center}
\begin{tabular}{|c|p{13.2cm}|}
\hline
$\mu$ & $l_{\mu}$ \\
\hline
$ \omega_2$ &
$d_{12}, d_{15}$  \\
 \hline
 $ \omega_1+\omega_3$ &
$d_{12}\wedge d_{13}, 2d_{12}\wedge d_{35}+
d_{23}\wedge d_{15}-d_{13}\wedge d_{25}, d_{12}\wedge d_{15}, d_{15}\wedge d_{25} $ \\
\hline
$ 2\omega_3$ & $
d_{12}\wedge d_{13}\wedge d_{23}, d_{12}\wedge d_{13}\wedge d_{25}-d_{12}\wedge d_{23}\wedge d_{15}, d_{12}\wedge d_{15}\wedge d_{25}$
\\
\hline
$  2\omega_1+\omega_4$ & $
d_{12}\wedge d_{13}\wedge d_{14}, d_{15}\wedge d_{25}\wedge d_{35}, 2 d_{12}\wedge d_{13}\wedge d_{45}-2d_{12}\wedge d_{14}\wedge d_{35}+
2d_{13}\wedge d_{14}\wedge d_{25}-
d_{13}\wedge d_{24}\wedge d_{15}-d_{23}\wedge d_{14}\wedge d_{15}
+d_{12}\wedge d_{34}\wedge d_{15}, d_{12}\wedge d_{35}\wedge d_{45}-d_{13}\wedge d_{25}\wedge d_{45}+
d_{23}\wedge d_{15}\wedge d_{45}+
d_{14}\wedge d_{25}\wedge d_{35}-d_{24}\wedge d_{15}\wedge d_{35}
+d_{34}\wedge d_{15}\wedge d_{25}, d_{12}\wedge d_{13}\wedge d_{15}, d_{12}\wedge d_{15}\wedge d_{35}-
d_{13}\wedge d_{15}\wedge d_{25}$
\\
\hline
$\omega_1+\omega_3+\omega_4$ &
$d_{12}\wedge d_{13}\wedge d_{23}\wedge d_{14}, d_{12}\wedge d_{13}\wedge d_{23}\wedge d_{15},2d_{12}\wedge d_{13}\wedge d_{23}\wedge d_{45}-
d_{12}\wedge d_{13}\wedge d_{24}\wedge d_{35}
+d_{12}\wedge d_{13}\wedge d_{34}\wedge d_{25}
+d_{12}\wedge d_{23}\wedge d_{14}\wedge d_{35}-d_{12}\wedge d_{23}\wedge d_{34}\wedge d_{15}
-d_{13}\wedge d_{23}\wedge d_{14}\wedge d_{25}
+d_{13}\wedge d_{23}\wedge d_{24}\wedge d_{15}, d_{12}\wedge d_{23}\wedge d_{14}\wedge d_{15}
-2d_{12}\wedge d_{13}\wedge d_{14}\wedge d_{25}
+d_{12}\wedge d_{13}\wedge d_{24}\wedge d_{15},d_{12}\wedge d_{13}\wedge d_{25}\wedge d_{45}
-d_{12}\wedge d_{23}\wedge d_{15}\wedge d_{45}
-d_{12}\wedge d_{14}\wedge d_{25}\wedge d_{35}
+d_{12}\wedge d_{24}\wedge d_{15}\wedge d_{35}
-d_{12}\wedge d_{34}\wedge d_{15}\wedge d_{25}, d_{12}\wedge d_{13}\wedge d_{25}\wedge d_{35}-
d_{12}\wedge d_{23}\wedge d_{15}\wedge d_{35}
+d_{13}\wedge d_{23}\wedge d_{15}\wedge d_{25},d_{12}\wedge d_{13}\wedge d_{15}\wedge d_{25}, d_{12}\wedge d_{15}\wedge d_{25}\wedge d_{35}$ \\
\hline
$ 3\omega_1$ & $
d_{12}\wedge d_{13}\wedge d_{14}\wedge d_{15}, d_{12}\wedge d_{13}\wedge d_{15}\wedge d_{45}
-d_{12}\wedge d_{14}\wedge d_{15}\wedge d_{35}
+d_{13}\wedge d_{14}\wedge d_{15}\wedge d_{25}, d_{12}\wedge d_{15}\wedge d_{35}\wedge d_{45}
 -
 d_{13}\wedge d_{15}\wedge d_{25}\wedge d_{45}
 +d_{14}\wedge d_{15}\wedge d_{25}\wedge d_{35}, d_{15}\wedge d_{25}\wedge d_{35}\wedge d_{45}$ \\
\hline
$\omega_2+2\omega_4
$ & $ d_{12}\wedge d_{13}\wedge d_{23}\wedge d_{14}\wedge d_{24}, 2d_{12}\wedge d_{13}\wedge d_{23}\wedge d_{14}\wedge d_{45}+
 2d_{12}\wedge d_{13}\wedge d_{14}\wedge d_{24}\wedge d_{35}
 -2d_{12}\wedge d_{13}\wedge d_{14}\wedge d_{34}\wedge d_{25}
 +d_{12}\wedge d_{13}\wedge d_{24}\wedge d_{34}\wedge d_{15}
 +d_{12}\wedge d_{23}\wedge d_{14}\wedge d_{34}\wedge d_{15}-
 d_{13}\wedge d_{23}\wedge d_{14}\wedge d_{24}\wedge d_{15},d_{12}\wedge d_{13}\wedge d_{23}\wedge d_{14}\wedge d_{25}-
 d_{12}\wedge d_{13}\wedge d_{23}\wedge d_{24}\wedge d_{15}, 3d_{12}\wedge d_{23}\wedge d_{14}\wedge d_{15}\wedge d_{35}-
 3d_{13}\wedge d_{23}\wedge d_{14}\wedge d_{15}\wedge d_{25}
 +4d_{12}\wedge d_{13}\wedge d_{23}\wedge d_{15}\wedge d_{45}
 -2d_{12}\wedge d_{13}\wedge d_{14}\wedge d_{25}\wedge d_{35}
 -d_{12}\wedge d_{13}\wedge d_{24}\wedge d_{15}\wedge d_{35}+
 d_{12}\wedge d_{13}\wedge d_{34}\wedge d_{15}\wedge d_{25}, d_{12}\wedge d_{13}\wedge d_{23}\wedge d_{15}\wedge d_{25}, d_{12}\wedge d_{13}\wedge d_{15}\wedge d_{25}\wedge d_{35}$ \\
\hline
$2\omega_1+\omega_3
$ & $ d_{12}\wedge d_{13}\wedge d_{23}\wedge d_{14}\wedge d_{15}, 4d_{12}\wedge d_{13}\wedge d_{23}\wedge d_{15}\wedge d_{45}-
 2d_{12}\wedge d_{13}\wedge d_{14}\wedge d_{25}\wedge d_{35}
 -d_{12}\wedge d_{13}\wedge d_{24}\wedge d_{15}\wedge d_{35}
 +3d_{12}\wedge d_{23}\wedge d_{14}\wedge d_{15}\wedge d_{35}
 -3d_{13}\wedge d_{23}\wedge d_{14}\wedge d_{15}\wedge d_{25}+
 d_{12}\wedge d_{13}\wedge d_{34}\wedge d_{15}\wedge d_{25}, 2d_{12}\wedge d_{13}\wedge d_{25}\wedge d_{35}\wedge d_{45}-
 2d_{12}\wedge d_{23}\wedge d_{15}\wedge d_{35}\wedge d_{45}
 +d_{12}\wedge d_{34}\wedge d_{15}\wedge d_{25}\wedge d_{35}
 +2d_{13}\wedge d_{23}\wedge d_{15}\wedge d_{25}\wedge d_{45}
 -d_{13}\wedge d_{24}\wedge d_{15}\wedge d_{25}\wedge d_{35}+
 d_{23}\wedge d_{14}\wedge d_{15}\wedge d_{25}\wedge d_{35}, d_{12}\wedge d_{13}\wedge d_{14}\wedge d_{15}\wedge d_{25}, d_{12}\wedge d_{15}\wedge d_{25}\wedge d_{35}\wedge d_{45}$
\\
\hline
$ 3\omega_4
$ & $ {\hat d_{15}}\wedge {\hat d_{25}}
\wedge
{\hat d_{35}}\wedge {\hat d_{45}}, {\hat d_{34}}\wedge
{\hat d_{15}}\wedge {\hat d_{25}}
\wedge
{\hat d_{45}}-{\hat d_{24}}\wedge {\hat d_{15}}
\wedge {\hat d_{35}}\wedge {\hat d_{45}}+
{\hat d_{14}}\wedge {\hat d_{25}}
\wedge {\hat d_{35}}\wedge {\hat d_{45}}, {\hat d_{14}}\wedge
{\hat d_{34}}\wedge {\hat d_{25}}
\wedge
{\hat d_{45}}-{\hat d_{24}}\wedge {\hat d_{34}}
\wedge {\hat d_{15}}\wedge {\hat d_{45}}-
{\hat d_{14}}\wedge {\hat d_{24}}
\wedge {\hat d_{35}}\wedge {\hat d_{45}}, {\hat d_{14}}\wedge {\hat d_{24}}
\wedge
{\hat d_{34}}\wedge {\hat d_{45}},$ \\
\hline
$\omega_1+\omega_2+\omega_4$ &
$
{\hat d_{34}}\wedge {\hat d_{25}}
\wedge
{\hat d_{35}}\wedge {\hat d_{45}}, {\hat d_{12}}\wedge {\hat d_{34}}
\wedge
{\hat d_{35}}\wedge {\hat d_{45}}-
{\hat d_{13}}\wedge {\hat d_{34}}
\wedge
{\hat d_{25}}\wedge {\hat d_{45}}
+{\hat d_{23}}\wedge {\hat d_{34}}
\wedge
{\hat d_{15}}\wedge {\hat d_{45}}
+{\hat d_{24}}\wedge {\hat d_{34}}
\wedge
{\hat d_{15}}\wedge {\hat d_{35}}
-{\hat d_{14}}\wedge {\hat d_{34}}
\wedge
{\hat d_{25}}\wedge {\hat d_{35}}, {\hat d_{23}}\wedge {\hat d_{24}}
\wedge
{\hat d_{35}}\wedge {\hat d_{45}}-
{\hat d_{23}}\wedge {\hat d_{34}}
\wedge
{\hat d_{25}}\wedge {\hat d_{45}}
-{\hat d_{24}}\wedge {\hat d_{34}}
\wedge
{\hat d_{25}}\wedge {\hat d_{35}}, {\hat d_{12}}\wedge {\hat d_{23}}
\wedge
{\hat d_{34}}\wedge {\hat d_{45}}+
{\hat d_{12}}\wedge {\hat d_{24}}
\wedge
{\hat d_{34}}\wedge {\hat d_{35}}
-{\hat d_{13}}\wedge {\hat d_{23}}
\wedge
{\hat d_{24}}\wedge {\hat d_{45}}-{\hat d_{13}}\wedge {\hat d_{24}}
\wedge
{\hat d_{34}}\wedge {\hat d_{25}}
+{\hat d_{23}}\wedge {\hat d_{14}}
\wedge
{\hat d_{24}}\wedge {\hat d_{35}}
-{\hat d_{23}}\wedge {\hat d_{14}}
\wedge
{\hat d_{34}}\wedge {\hat d_{25}}
+2{\hat d_{23}}\wedge {\hat d_{24}}
\wedge
{\hat d_{34}}\wedge {\hat d_{15}},{\hat d_{24}}\wedge {\hat d_{34}}
\wedge
{\hat d_{35}}\wedge {\hat d_{45}}, {\hat d_{13}}\wedge {\hat d_{24}}
\wedge
{\hat d_{34}}\wedge {\hat d_{45}}+2
{\hat d_{14}}\wedge {\hat d_{24}}
\wedge
{\hat d_{34}}\wedge {\hat d_{35}}-{\hat d_{23}}\wedge {\hat d_{14}}
\wedge
{\hat d_{34}}\wedge {\hat d_{45}}, {\hat d_{23}}\wedge {\hat d_{24}}
\wedge
{\hat d_{34}}\wedge {\hat d_{45}}, {\hat d_{23}}\wedge {\hat d_{14}}
\wedge
{\hat d_{24}}\wedge {\hat d_{34}}$ \\
 \hline
 $\omega_1+2\omega_4
$ & $ {\hat d_{25}}
\wedge
{\hat d_{35}}\wedge {\hat d_{45}}, {\hat d_{34}}
\wedge
{\hat d_{15}}\wedge {\hat d_{25}}-
{\hat d_{24}}
\wedge
{\hat d_{15}}\wedge {\hat d_{35}}+{\hat d_{14}}
\wedge
{\hat d_{25}}\wedge {\hat d_{35}}-
{\hat d_{23}}
\wedge
{\hat d_{15}}\wedge {\hat d_{45}}+{\hat d_{13}}
\wedge
{\hat d_{25}}\wedge {\hat d_{45}}-{\hat d_{12}}
\wedge
{\hat d_{35}}\wedge {\hat d_{45}}, {\hat d_{24}}
\wedge
{\hat d_{35}}\wedge {\hat d_{45}}, 2 {\hat d_{24}}
\wedge
{\hat d_{34}}\wedge {\hat d_{15}}-2
{\hat d_{14}}
\wedge
{\hat d_{34}}\wedge {\hat d_{25}}-{\hat d_{23}}
\wedge
{\hat d_{14}}\wedge {\hat d_{45}}+
{\hat d_{13}}
\wedge
{\hat d_{24}}\wedge {\hat d_{45}}-{\hat d_{12}}
\wedge
{\hat d_{34}}\wedge {\hat d_{45}}+2{\hat d_{14}}
\wedge
{\hat d_{24}}\wedge {\hat d_{35}}, {\hat d_{24}}
\wedge
{\hat d_{34}}\wedge {\hat d_{45}}
,
{\hat d_{14}}
\wedge
{\hat d_{24}}\wedge {\hat d_{34}},$ \\
\hline $
 2\omega_2$ & $ {\hat d_{34}}
\wedge
{\hat d_{35}}\wedge {\hat d_{45}}, {\hat d_{24}}
\wedge
{\hat d_{34}}\wedge {\hat d_{35}}+ {\hat d_{23}}
\wedge
{\hat d_{34}}\wedge {\hat d_{45}},{\hat d_{23}}
\wedge
{\hat d_{24}}\wedge {\hat d_{34}}, $ \\
\hline
$\omega_2+\omega_4
$ & ${\hat d_{35}}\wedge {\hat d_{45}}, {\hat d_{23}}\wedge {\hat d_{45}}
+{\hat d_{24}}\wedge {\hat d_{35}}-{\hat d_{34}}\wedge {\hat d_{25}}, {\hat d_{34}}\wedge {\hat d_{45}}, {\hat d_{24}}\wedge {\hat d_{34}}$ \\
\hline $
 \omega_3
$ & ${\hat d_{45}},{\hat d_{34}}$
\\
\hline
\end{tabular}
\end{center}
\end{table}
\begin{table}
\caption{\bf Weights and weight vectors  for $\mbox{sl}_5$ module $V(\omega_1+\omega_3)$}
\begin{center}
\begin{tabular}{|c|p{3.5cm}|p{1.5cm}|c|p{4.5cm}|p{1.5cm}|}
\hline
$ (i,j)  $ &   $\overrightarrow{w}_{i}^{\omega_1+\omega_3}$  &
  $v_{i,j}^{\omega_1+\omega_3}$ & $ (i,j)  $ &   $\overrightarrow{w}_{i}^{\omega_1+\omega_3}$  &
  $v_{i,j}^{\omega_1+\omega_3}$ \\
\hline
(1,1) & $(1,0,1,0)=\omega_1+\omega_3  $ & $d_{12}\wedge d_{13}$ &
(18,1) & $(-1,-1,1,1)=\omega_1+\omega_3-2\alpha_1-2\alpha_2-\alpha_3  $ & $d_{23}\wedge d_{34}
$   \\
\hline
(2,1) & $(-1,1,1,0)=\omega_1+\omega_3-\alpha_1  $ & $d_{12}\wedge d_{23}
$ &(19,1) & $(-2,1,1,-1)=\omega_1+\omega_3-2\alpha_1-\alpha_2-\alpha_3-\alpha_4  $ & $d_{23}\wedge d_{25}
$
 \\
\hline
(3,1) & $(1,1,-1,1)=\omega_1+\omega_3-\alpha_3  $ & $d_{12}\wedge d_{14}
$ & (20,1)& $(1,-1,-1,2) =\omega_1+\omega_3-\alpha_1-2\alpha_2-2\alpha_3 $ & $d_{14}\wedge d_{34}$
 \\
\hline
(4,1) & $(-1,2,-1,1)=\omega_1+\omega_3-\alpha_1-\alpha_3  $ & $
d_{12}\wedge d_{24}
$ & (21,1) & $(1,-1,0,0)==\omega_1+\omega_3-\alpha_1-2\alpha_2-2\alpha_3-\alpha_4  $ & $d_{14}\wedge d_{35}-d_{13}\wedge d_{45}
$
\\
\hline
(5,1) & $(1,1,0,-1)=\omega_1+\omega_3-\alpha_3-\alpha_4  $ & $d_{12}\wedge d_{15}$ &
(21,2) & $(1,-1,0,0)=\omega_1+\omega_3 -\alpha_1-2\alpha_2-2\alpha_3-\alpha_4  $ & $d_{14}\wedge d_{35}-d_{34}\wedge d_{15}
$  \\
\hline
(6,1) & $(0,-1,2,0)=\omega_1+\omega_3-\alpha_1-\alpha_2  $ &  $d_{13}\wedge d_{23}
$ &(21,3) & $(1,-1,0,0) =\omega_1+\omega_3 -\alpha_1-2\alpha_2-2\alpha_3-\alpha_4 $ & $d_{34}\wedge d_{15}+d_{13}\wedge d_{45}
$
\\
\hline
(7,1) & $(2,-1,0,1)=\omega_1+\omega_3-\alpha_2-\alpha_3  $ & $d_{13}\wedge d_{14}
$ & (22,1)& $(-1,-1,2,-1)=\omega_1+\omega_3-2\alpha_1-2\alpha_2-\alpha_3-\alpha_4  $ & $d_{23}\wedge d_{35}
$
 \\
\hline
(8,1) & $(0,0,0,1)=\omega_1+\omega_3-\alpha_1-\alpha_2-\alpha_3  $ & $d_{12}\wedge d_{34}-d_{13}\wedge d_{24}
$ & $(23,1) $ & $(-1,0,-1,2)=\omega_1+\omega_3-2\alpha_1-2\alpha_2-2\alpha_3  $ & $d_{24}\wedge d_{34}
$
 \\
\hline
(8,2) & $(0,0,0,1)=\omega_1+\omega_3-\alpha_1-\alpha_2-\alpha_3  $ & $d_{12}\wedge d_{34}+d_{23}\wedge d_{14}
$ & (24,1) & $(-2,2,-1,0)=\omega_1+\omega_3-2\alpha_1-\alpha_2-2\alpha_3-\alpha_4  $ & $d_{24}\wedge d_{25}
$
 \\
\hline
(8,3) & $(0,0,0,1)=\omega_1+\omega_3-\alpha_1-\alpha_2-\alpha_3  $ & $d_{13}\wedge d_{24}-d_{23}\wedge d_{14}
$  & (25,1) & $(0,1,0,-2)=\omega_1+\omega_3-\alpha_1-\alpha_2-2\alpha_3-2\alpha_4  $ & $d_{15}\wedge d_{25}
$
\\
\hline
(9,1) & $(-1,2,0,-1)=\omega_1+\omega_3-\alpha_1-\alpha_3-\alpha_4  $ & $d_{12}\wedge d_{25}
$  & (26,1) & $(-1,0,0,0)=\omega_1+\omega_3-2\alpha_1-2\alpha_2-2\alpha_3-\alpha_4  $ & $d_{34}\wedge d_{25}-d_{24}\wedge d_{35}$
\\
\hline
(10,1) & $(2,-1,1,-1)=\omega_1+\omega_3-\alpha_2-\alpha_3-\alpha_4  $ & $
d_{13}\wedge d_{15}
$ & (26,2) & $(-1,0,0,0)=\omega_1+\omega_3-2\alpha_1-2\alpha_2-2\alpha_3-\alpha_4  $ & $d_{23}\wedge d_{45}-d_{24}\wedge d_{35}$
 \\
\hline
(11,1) & $(0,0,1,-1)=\omega_1+\omega_3-\alpha_1-\alpha_2-\alpha_3-\alpha_4  $ & $d_{12}\wedge d_{35}-d_{13}\wedge d_{25}
$ & (26,3) & $(-1,0,0,0)=\omega_1+\omega_3-2\alpha_1-2\alpha_2-2\alpha_3-\alpha_4  $ & $d_{34}\wedge d_{25}+d_{23}\wedge d_{45}$
 \\
\hline
(11,2) & $(0,0,1,-1) =\omega_1+\omega_3-\alpha_1-\alpha_2-\alpha_3-\alpha_4 $ & $d_{12}\wedge d_{35}+d_{23}\wedge d_{15}
$ &
(27,1) & $(1,0,-2,1) =\omega_1+\omega_3-\alpha_1-2\alpha_2-3\alpha_3-\alpha_4 $ & $d_{14}\wedge d_{45}
$
 \\
\hline
(11,3) & $(0,0,1,-1) =\omega_1+\omega_3-\alpha_1-\alpha_2-\alpha_3-\alpha_4 $ & $d_{23}\wedge d_{15}-d_{13}\wedge d_{25}
$  & (28,1) & $(1,-1,1,-2)=\omega_1+\omega_3-\alpha_1-2\alpha_2-2\alpha_3-2\alpha_4  $ & $d_{15}\wedge d_{35}
$
\\
\hline
(12,1) & $(1,-2,1,1)=\omega_1+\omega_3-\alpha_1-2\alpha_2-\alpha_3  $ & $d_{13}\wedge d_{34}
$   & (29,1) & $(-1,1,-2,1)=\omega_1+\omega_3-2\alpha_1-2\alpha_2-3\alpha_3-\alpha_4  $ & $d_{24}\wedge d_{45}
$
\\
\hline
(13,1) & $(-2,1,0,1)=\omega_1+\omega_3-2\alpha_1-\alpha_2-\alpha_3  $ & $d_{23}\wedge d_{24}
$ & (30,1)& $(0,-2,1,0)=\omega_1+\omega_3-2\alpha_1-3\alpha_2-2\alpha_3-\alpha_4  $ & $d_{34}\wedge d_{35}
$
 \\
\hline
(14,1) & $(0,1,-2,2)=\omega_1+\omega_3-\alpha_1-\alpha_2-2\alpha_3  $ & $d_{14}\wedge d_{24}$
& (31,1) & $(1,0,-1,-1)=\omega_1+\omega_3-\alpha_1-2\alpha_2-3\alpha_3-2\alpha_4  $ & $d_{15}\wedge d_{45}$
\\
\hline
(15,1) & $(2,0,-1,0)=\omega_1+\omega_3-\alpha_2-2\alpha_3-\alpha_4  $ & $d_{14}\wedge d_{15}
$ & (32,1)& $(-1,0,1,-2)=\omega_1+\omega_3-2\alpha_1-2\alpha_2-2\alpha_3-2\alpha_4  $ & $d_{25}\wedge d_{35}
$  \\
\hline
(16,1) & $(0,1,-1,0)=\omega_1+\omega_3-\alpha_1-\alpha_2-2\alpha_3-\alpha_4  $ & $d_{12}\wedge d_{45}-
d_{14}\wedge d_{25}
$
& (33,1)& $(0,-1,-1,1)=\omega_1+\omega_3-2\alpha_1-3\alpha_2-3\alpha_3-\alpha_4  $ & $d_{34}\wedge d_{45}
$
 \\
\hline
(16,2) & $(0,1,-1,0)=\omega_1+\omega_3-\alpha_1-\alpha_2-2\alpha_3-\alpha_4  $ & $d_{12}\wedge d_{45}+
d_{24}\wedge d_{15}
$  & (34,1) & $(-1,1,-1,-1)=\omega_1+\omega_3-2\alpha_1-2\alpha_2-3\alpha_3-2\alpha_4  $ & $d_{25}\wedge d_{45}
$
\\
\hline
(16,3) & $(0,1,-1,0)=\omega_1+\omega_3-\alpha_1-\alpha_2-2\alpha_3-\alpha_4  $  & $d_{14}\wedge d_{25}-
d_{24}\wedge d_{15}
$ & (35,1) & $(0,-1,0,-1) =\omega_1+\omega_3-2\alpha_1-3\alpha_2-3\alpha_3-2\alpha_4 $ & $d_{35}\wedge d_{45}$
 \\
\hline
(17,1) & $(1,-2,2,-1) =\omega_1+\omega_3-\alpha_1-2\alpha_2-\alpha_3-\alpha_4 $ & $d_{13}\wedge d_{35}
$  & & &
 \\
\hline
\end{tabular}
\end{center}
\end{table}
\begin{table}
\caption{\bf Weights and weight vectors  for $\mbox{sl}_5$ module $V(2\omega_{1}+\omega_4)$}
\begin{center}
\begin{tabular}{|c|c|p{5.8cm}|}
\hline
$ (i,j)  $ &   $\overrightarrow{w}_{i}^{2\omega_{1}+\omega_4}$  &
  $v_{i,j}^{2\omega_{1}+\omega_4}$  \\
\hline
(1,1) & $(2,0,0,1)= 2\omega_{1}+\omega_4$  & $d_{12}\wedge d_{13}\wedge d_{14}$ \\
\hline
(2,1) & $(0,1,0,1)= 2\omega_{1}+\omega_4 -\alpha_1$ & $d_{12}\wedge d_{23}\wedge d_{14}+d_{12}\wedge d_{13}\wedge d_{24}$ \\
\hline
(3,1) & $(2,0,1,-1)=2\omega_{1}+\omega_4 -\alpha_4$ & $d_{12}\wedge d_{13}\wedge d_{15}$ \\
\hline
(4,1) & $(-2,2,0,1)=2\omega_{1}+\omega_4 -2\alpha_1 $ & $d_{12}\wedge d_{23}\wedge d_{24}$ \\
\hline
(5,1) & $(1,-1,1,1)=2\omega_{1}+\omega_4 -\alpha_1-\alpha_2$ & $d_{13}\wedge d_{23}\wedge d_{14}+d_{12}\wedge d_{13}\wedge d_{34}$ \\
\hline
(6,1) & $(0,1,1,-1)=2\omega_{1}+\omega_4 -\alpha_1-\alpha_4$ & $d_{12}\wedge d_{23}\wedge d_{15}+d_{12}\wedge d_{13}\wedge d_{25}$ \\
\hline
(7,1) & $(2,1,-1,0)=2\omega_{1}+\omega_4 -\alpha_3-\alpha_4 $ & $d_{12}\wedge d_{14}\wedge d_{15}$ \\
\hline
(8,1) & $ (-1,0,1,1)= 2\omega_{1}+\omega_4 -2\alpha_1-\alpha_2 $ & $d_{13}\wedge d_{23}\wedge d_{24}+d_{12}\wedge d_{23}\wedge d_{34}$ \\
\hline
(9,1) & $ (-2,2,1,-1)=2\omega_{1}+\omega_4 -2\alpha_1-\alpha_4  $ & $d_{12}\wedge d_{23}\wedge d_{25}$ \\
\hline
(10,1) & $ (1,0,-1,2)=2\omega_{1}+\omega_4 -\alpha_1-\alpha_2-\alpha_3 $  & $-d_{13}\wedge d_{14}\wedge d_{24}+d_{12}\wedge d_{14}\wedge d_{34}$ \\
\hline
(11,1) & $ (1,-1,2,-1)=2\omega_{1}+\omega_4 -\alpha_1-\alpha_2-\alpha_4$  & $d_{13}\wedge d_{23}\wedge d_{15}+d_{12}\wedge d_{13}\wedge d_{35}$ \\
\hline
(12,1) & $(0,2,-1,0)=2\omega_{1}+\omega_4 -\alpha_1-\alpha_3-\alpha_4 $ & $d_{12}\wedge d_{24}\wedge d_{15}+d_{12}\wedge d_{14}\wedge d_{25}$ \\
\hline
(13,1) & $(3,-1,0,0)=2\omega_{1}+\omega_4 -\alpha_2-\alpha_3-\alpha_4$ & $d_{13}\wedge d_{14}\wedge d_{15}$ \\
\hline
(14,1) & $(0,-2,2,1)=2\omega_{1}+\omega_4 -2\alpha_1-2\alpha_2$ & $d_{13}\wedge d_{23}\wedge d_{34}$ \\
\hline
(15,1) & $(-1,1,-1,2)=2\omega_{1}+\omega_4 -2\alpha_1-\alpha_2 -\alpha_3$ & $d_{12}\wedge d_{24}\wedge d_{34}-d_{23}\wedge d_{14}\wedge d_{24}$ \\
\hline
(16,1) & $(-1,0,2,-1)=2\omega_{1}+\omega_4 -2\alpha_1-\alpha_2-\alpha_4 $ & $d_{12}\wedge d_{23}\wedge d_{35}+d_{13}\wedge d_{23}\wedge d_{25}$ \\
\hline
(17,1) & $(-2,3,-1,0)= 2\omega_{1}+\omega_4 -2\alpha_1-\alpha_3-\alpha_4$ & $d_{12}\wedge d_{24}\wedge d_{25}$ \\
\hline
(18,1) & $ (1,0,0,0)=2\omega_{1}+\omega_4 -\alpha_1-\alpha_2-\alpha_3-\alpha_4$   &
$d_{13}\wedge d_{24}\wedge d_{15}+d_{13}\wedge d_{14}\wedge d_{25}+
d_{23}\wedge d_{14}\wedge d_{15}, $ \\
\hline
(18,2) & $ (1,0,0,0)=2\omega_{1}+\omega_4 -\alpha_1-\alpha_2-\alpha_3-\alpha_4$  &
$
 d_{13}\wedge d_{24}\wedge d_{15}+d_{13}\wedge d_{14}\wedge d_{25}+
d_{12}\wedge d_{14}\wedge d_{35}+d_{12}\wedge d_{34}\wedge d_{15},$ \\
\hline
(18,3) & $ (1,0,0,0)=2\omega_{1}+\omega_4 -\alpha_1-\alpha_2-\alpha_3-\alpha_4$  & $ d_{13}\wedge d_{24}\wedge d_{15}-
d_{13}\wedge d_{14}\wedge d_{25}-
d_{12}\wedge d_{34}\wedge d_{15}
+d_{12}\wedge d_{14}\wedge d_{35} $ \\
\hline
(18,4) & $ (1,0,0,0)=2\omega_{1}+\omega_4 -\alpha_1-\alpha_2-\alpha_3-\alpha_4$  & $ -d_{23}\wedge d_{14}\wedge d_{15}+
d_{13}\wedge d_{24}\wedge d_{15}+
d_{12}\wedge d_{14}\wedge d_{35}
+d_{12}\wedge d_{13}\wedge d_{45} $ \\
\hline
(19,1) & $ (0,-1,0,2)=2\omega_{1}+\omega_4 -2\alpha_1-2\alpha_2-\alpha_3$  & $d_{13}\wedge d_{24}\wedge d_{34}-d_{23}\wedge d_{14}\wedge d_{34}$ \\
\hline
(20,1) & $ (0,-2,3,-1)=2\omega_{1}+\omega_4 -2\alpha_1-2\alpha_2-\alpha_4$  & $d_{13}\wedge d_{23}\wedge d_{35}$ \\
\hline
(21,1) & $(-1,1,0,0)=2\omega_{1}+\omega_4 -2\alpha_1-\alpha_2-\alpha_3-\alpha_4$ & $d_{12}\wedge d_{34}\wedge d_{25}+
d_{12}\wedge d_{24}\wedge d_{35}+2d_{13}\wedge d_{24}\wedge d_{25}
+
d_{23}\wedge d_{14}\wedge d_{25}+
d_{23}\wedge d_{24}\wedge d_{15}
$ \\
\hline
(21,2) & $(-1,1,0,0)=2\omega_{1}+\omega_4 -2\alpha_1-\alpha_2-\alpha_3-\alpha_4$ & $d_{12}\wedge d_{24}\wedge d_{35}+
d_{12}\wedge d_{34}\wedge d_{25}+
d_{13}\wedge d_{24}\wedge d_{25}
$ \\
\hline
(21,3) & $(-1,1,0,0)=2\omega_{1}+\omega_4 -2\alpha_1-\alpha_2-\alpha_3-\alpha_4$ & $-d_{23}\wedge d_{14}\wedge d_{25}+
d_{13}\wedge d_{24}\wedge d_{25}
+d_{12}\wedge d_{24}\wedge d_{35}+
d_{12}\wedge d_{23}\wedge d_{45}
$ \\
\hline
(21,4) & $(-1,1,0,0)=2\omega_{1}+\omega_4 -2\alpha_1-\alpha_2-\alpha_3-\alpha_4$ & $d_{12}\wedge d_{24}\wedge d_{35}
-d_{12}\wedge d_{34}\wedge d_{25}
+
d_{23}\wedge d_{24}\wedge d_{15}-
d_{23}\wedge d_{14}\wedge d_{25}
$ \\
\hline
(22,1) & $(2,-2,1,0)=2\omega_{1}+\omega_4 -\alpha_1-2\alpha_2-\alpha_3-\alpha_4$ & $d_{13}\wedge d_{34}\wedge d_{15}+d_{13}\wedge d_{14}\wedge d_{35}
$ \\
\hline
(23,1) & $(1,1,-2,1)=2\omega_{1}+\omega_4 -\alpha_1-\alpha_2-2\alpha_3 -\alpha_4$ & $d_{12}\wedge d_{14}\wedge d_{45}+d_{14}\wedge d_{24}\wedge d_{15}
$ \\
\hline
(24,1) & $(1,0,1,-2)=2\omega_{1}+\omega_4 -\alpha_1-\alpha_2-\alpha_3-2\alpha_4 $  & $d_{12}\wedge d_{15}\wedge d_{35}-d_{13}\wedge d_{15}\wedge d_{25}
$ \\
\hline
(25,1) & $(0,0,-2,3)=2\omega_{1}+\omega_4 -2\alpha_1-2\alpha_2-2\alpha_3$ & $d_{14}\wedge d_{24}\wedge d_{34}$ \\
\hline
(26,1) & $(0,-1,1,0)=2\omega_{1}+\omega_4 -2\alpha_1-2\alpha_2-\alpha_3-\alpha_4$ & $d_{13}\wedge d_{24}\wedge d_{35}+d_{13}\wedge d_{34}\wedge d_{25}+
d_{23}\wedge d_{14}\wedge d_{35}+d_{23}\wedge d_{34}\wedge d_{15}$ \\
\hline
(26,2) & $(0,-1,1,0)=2\omega_{1}+\omega_4 -2\alpha_1-2\alpha_2-\alpha_3-\alpha_4$ & $-d_{23}\wedge d_{14}\wedge d_{35}+d_{13}\wedge d_{34}\wedge d_{25}+
2d_{13}\wedge d_{24}\wedge d_{35}+d_{12}\wedge d_{34}\wedge d_{35}+
d_{13}\wedge d_{23}\wedge d_{45}$ \\
\hline
(26,3) & $(0,-1,1,0)=2\omega_{1}+\omega_4 -2\alpha_1-2\alpha_2-\alpha_3-\alpha_4$ & $d_{13}\wedge d_{24}\wedge d_{35}+d_{13}\wedge d_{23}\wedge d_{45}-d_{23}\wedge d_{14}\wedge d_{35}
 $ \\
\hline
(26,4) & $(0,-1,1,0)=2\omega_{1}+\omega_4 -2\alpha_1-2\alpha_2-\alpha_3-\alpha_4$ & $d_{23}\wedge d_{34}\wedge d_{15}-d_{23}\wedge d_{14}\wedge d_{35}-
d_{13}\wedge d_{34}\wedge d_{25}+d_{13}\wedge d_{24}\wedge d_{35}$ \\
\hline
(27,1) & $(-1,2,-2,1)=2\omega_{1}+\omega_4 -2\alpha_1-\alpha_2-2\alpha_3-\alpha_4$ & $d_{14}\wedge d_{24}\wedge d_{25}+d_{12}\wedge d_{24}\wedge d_{45}$ \\
\hline
(28,1) & $(-1,1,1,-2)=2\omega_{1}+\omega_4 -2\alpha_1-\alpha_2-\alpha_3-2\alpha_4$ & $d_{12}\wedge d_{25}\wedge d_{35}-d_{23}\wedge d_{15}\wedge d_{25}$ \\
\hline
(29,1) & $(2,-1,-1,1)=2\omega_{1}+\omega_4 -\alpha_1-2\alpha_2-2\alpha_3-\alpha_4$ & $d_{14}\wedge d_{34}\wedge d_{15}+d_{13}\wedge d_{14}\wedge d_{45}$ \\
\hline
(30,1) & $(1,1,-1,-1)=2\omega_{1}+\omega_4 -\alpha_1-\alpha_2-2\alpha_3-2\alpha_4$ & $d_{12}\wedge d_{15}\wedge d_{45}-d_{14}\wedge d_{15}\wedge d_{25}$ \\
\hline
\end{tabular}
\end{center}
\end{table}

\begin{table}
\caption{\bf Weights and weight vectors  for $\mbox{sl}_5$ module $V(2\omega_{1}+\omega_4)$}
\begin{center}
\begin{tabular}{|c|c|p{4.5cm}|}
\hline
$ (i,j)  $ &   $\overrightarrow{w}_{i}^{2\omega_{1}+\omega_4}$  &
  $v_{i,j}^{2\omega_{1}+\omega_4}$  \\
\hline
(31,1) & $(-3,2,0,0) =2\omega_{1}+\omega_4 -3\alpha_1-\alpha_2-\alpha_3-\alpha_4$& $d_{23}\wedge d_{24}\wedge d_{25}$ \\
\hline
(32,1) & $(0,0,-1,1)=2\omega_{1}+\omega_4 -2\alpha_1-2\alpha_2-2\alpha_3-\alpha_4$ & $
d_{13}\wedge d_{24}\wedge d_{45}+2d_{14}\wedge d_{24}\wedge d_{35}-
d_{14}\wedge d_{34}\wedge d_{25}-d_{23}\wedge d_{14}\wedge d_{45}
+d_{24}\wedge d_{34}\wedge d_{15}
$ \\
\hline
(32,2) & $(0,0,-1,1)=2\omega_{1}+\omega_4 -2\alpha_1-2\alpha_2-2\alpha_3-\alpha_4$ & $d_{14}\wedge d_{24}\wedge d_{35}-d_{14}\wedge d_{34}\wedge d_{25}+
d_{24}\wedge d_{34}\wedge d_{15}
$ \\
\hline
(32,3) & $(0,0,-1,1)=2\omega_{1}+\omega_4 -2\alpha_1-2\alpha_2-2\alpha_3-\alpha_4$ & $
d_{12}\wedge d_{34}\wedge d_{45}+d_{13}\wedge d_{24}\wedge d_{45}+
d_{14}\wedge d_{34}\wedge d_{25}+d_{14}\wedge d_{24}\wedge d_{35}
$ \\
\hline
(32,4) & $(0,0,-1,1)=2\omega_{1}+\omega_4 -2\alpha_1-2\alpha_2-2\alpha_3-\alpha_4$ &
$d_{13}\wedge d_{24}\wedge d_{45}+
d_{23}\wedge d_{14}\wedge d_{45}+d_{14}\wedge d_{34}\wedge d_{25}+d_{24}\wedge d_{34}\wedge d_{15}
$ \\
\hline
(33,1) & $(-2,0,1,0)=2\omega_{1}+\omega_4 -3\alpha_1-2\alpha_2-\alpha_3-\alpha_4$ & $d_{23}\wedge d_{24}\wedge d_{35}+
d_{23}\wedge d_{34}\wedge d_{25}
$ \\
\hline
(34,1) & $(0,-1,2,-2)=2\omega_{1}+\omega_4 -2\alpha_1-2\alpha_2-\alpha_3-2\alpha_4$ &  $d_{13}\wedge d_{25}\wedge d_{35}-
d_{23}\wedge d_{15}\wedge d_{35}
$ \\
\hline
(35,1) & $(-1,2,-1,-1)=2\omega_{1}+\omega_4 -2\alpha_1-\alpha_2-2\alpha_3-2\alpha_4$ &  $d_{12}\wedge d_{25}\wedge d_{45}-
d_{24}\wedge d_{15}\wedge d_{25}
$ \\
\hline
(36,1) & $(2,-1,0,-1 )=2\omega_{1}+\omega_4 -\alpha_1-2\alpha_2-2\alpha_3-2\alpha_4$ &  $d_{13}\wedge d_{15}\wedge d_{45}-
d_{14}\wedge d_{15}\wedge d_{35}
$ \\
\hline
(37,1) & $(1,-3,2,0)=2\omega_{1}+\omega_4 -2\alpha_1-3\alpha_2-\alpha_3-\alpha_4$ &  $d_{13}\wedge d_{34}\wedge d_{35}
$ \\
\hline
(38,1) & $(-1,-2,2,0)=2\omega_{1}+\omega_4 -3\alpha_1-3\alpha_2-\alpha_3-\alpha_4$ & $d_{23}\wedge d_{34}\wedge d_{35}
$ \\
\hline
(39,1) & $ (-2,1,-1,1)=2\omega_{1}+\omega_4 -3\alpha_1-2\alpha_2-2\alpha_3-\alpha_4$ & $d_{23}\wedge d_{24}\wedge d_{45}+d_{24}\wedge d_{34}\wedge d_{25}
$ \\
\hline
(40,1) & $(1,-2,0,1)=2\omega_{1}+\omega_4 -2\alpha_1-3\alpha_2-2\alpha_3-\alpha_4$ & $d_{13}\wedge d_{34}\wedge d_{45}+d_{14}\wedge d_{34}\wedge d_{35}
$ \\
\hline
(41,1) & $(0,0,0,-1)=2\omega_{1}+\omega_4 -2\alpha_1-2\alpha_2-2\alpha_3-2\alpha_4$ & $d_{13}\wedge d_{25}\wedge d_{45}+d_{23}\wedge d_{15}\wedge d_{45}-
d_{14}\wedge d_{25}\wedge d_{35}-d_{24}\wedge d_{15}\wedge d_{35}
$ \\
\hline
(41,2) & $(0,0,0,-1)=2\omega_{1}+\omega_4 -2\alpha_1-2\alpha_2-2\alpha_3-2\alpha_4$ &
$d_{12}\wedge d_{35}\wedge d_{45}+
d_{13}\wedge d_{25}\wedge d_{45}-d_{24}\wedge d_{15}\wedge d_{35}-d_{34}\wedge d_{15}\wedge d_{25}
$ \\
\hline
(41,3) & $(0,0,0,-1)=2\omega_{1}+\omega_4 -2\alpha_1-2\alpha_2-2\alpha_3-2\alpha_4$ & $d_{13}\wedge d_{25}\wedge d_{45}-d_{23}\wedge d_{15}\wedge d_{45}
+d_{14}\wedge d_{25}\wedge d_{35}
-d_{24}\wedge d_{15}\wedge d_{35}
$ \\
\hline
(41,4) & $(0,0,0,-1)=2\omega_{1}+\omega_4 -2\alpha_1-2\alpha_2-2\alpha_3-2\alpha_4$  & $d_{14}\wedge d_{25}\wedge d_{35}-d_{24}\wedge d_{15}\wedge d_{35}
+d_{34}\wedge d_{15}\wedge d_{25}
$ \\
\hline
(42,1) & $ (0,1,-3,2)=2\omega_{1}+\omega_4 -2\alpha_1-2\alpha_2-3\alpha_3-\alpha_4$ & $d_{14}\wedge d_{24}\wedge d_{45}
$ \\
\hline
(43,1) & $ (0,0,1,-3)=2\omega_{1}+\omega_4 -2\alpha_1-2\alpha_2-2\alpha_3-3\alpha_4$ & $ d_{15}\wedge d_{25}\wedge d_{35}$ \\
\hline
(44,1) & $(-1,-1,0,1)=2\omega_{1}+\omega_4 -3\alpha_1-3\alpha_2-2\alpha_3-\alpha_4$ & $d_{23}\wedge d_{34}\wedge d_{45}+d_{24}\wedge d_{34}\wedge d_{35}
$ \\
\hline
(45,1) & $(-2,1,0,-1)=2\omega_{1}+\omega_4 -3\alpha_1-2\alpha_2-2\alpha_3-2\alpha_4$  & $d_{23}\wedge d_{25}\wedge d_{45}-d_{24}\wedge d_{25}\wedge d_{35}
$ \\
\hline
(46,1) & $(1,-2,1,-1)=2\omega_{1}+\omega_4 -2\alpha_1-3\alpha_2-2\alpha_3-2\alpha_4$ & $d_{13}\wedge d_{35}\wedge d_{45}-d_{34}\wedge d_{15}\wedge d_{35}
$ \\
\hline
(47,1) & $(0,1,-2,0)=2\omega_{1}+\omega_4 -2\alpha_1-2\alpha_2-3\alpha_3-2\alpha_4$ & $d_{14}\wedge d_{25}\wedge d_{45}-d_{24}\wedge d_{15}\wedge d_{45}
$ \\
\hline
(48,1) & $(1,-1,-2,2)=2\omega_{1}+\omega_4 -2\alpha_1-3\alpha_2-3\alpha_3-\alpha_4$ & $d_{14}\wedge d_{34}\wedge d_{45}
$ \\
\hline
(49,1) & $ (-1,-1,1,-1)=2\omega_{1}+\omega_4 -3\alpha_1-3\alpha_2-2\alpha_3-2\alpha_4$ & $d_{23}\wedge d_{35}\wedge d_{45}-d_{34}\wedge d_{25}\wedge d_{35}
$ \\
\hline
(50,1) & $ (1,-1,-1,0)=2\omega_{1}+\omega_4 -2\alpha_1-3\alpha_2-3\alpha_3-2\alpha_4$ & $d_{14}\wedge d_{35}\wedge d_{45}-d_{34}\wedge d_{15}\wedge d_{45}
$ \\
\hline
(51,1) & $ (-1,0,-2,2)=2\omega_{1}+\omega_4 -3\alpha_1-3\alpha_2-3\alpha_3-\alpha_4$ &  $d_{24}\wedge d_{34}\wedge d_{45}
$ \\
\hline
(52,1) & $ (0,1,-1,-2)=2\omega_{1}+\omega_4 -2\alpha_1-2\alpha_2-3\alpha_3-3\alpha_4$ & $d_{15}\wedge d_{25}\wedge d_{45}$  \\
\hline
(53,1) & $ (-1,0,-1,0)=2\omega_{1}+\omega_4 -3\alpha_1-3\alpha_2-3\alpha_3-2\alpha_4$ & $d_{24}\wedge d_{35}\wedge d_{45}-d_{34}\wedge d_{25}\wedge d_{45}
$ \\
\hline
(54,1) & $ (1,-1,0,-2)=2\omega_{1}+\omega_4 -2\alpha_1-3\alpha_2-3\alpha_3-3\alpha_4$ & $d_{15}\wedge d_{35}\wedge d_{45}
$ \\
\hline
(55,1) & $(-1,0,0,-2)=2\omega_{1}+\omega_4 -3\alpha_1-3\alpha_2-3\alpha_3-3\alpha_4$ & $d_{25}\wedge d_{35}\wedge d_{45}
$ \\
\hline
\end{tabular}
\end{center}
\end{table}
\begin{table}
\caption{\bf Weights and weight vectors  for $\mbox{sl}_5$ module $V(3\omega_1)$}
\begin{center}
\begin{tabular}{|c|p{6cm}|p{9cm}|}
\hline
$ i  $ &   $\overrightarrow{w}_{i}^{3\omega_1}$  &
  $v_{i}^{3\omega_1}$  \\
\hline
1 & $(3,0,0,0)=3\omega_1    $ & $d_{12}\wedge d_{13}\wedge d_{14}\wedge d_{15}$  \\
\hline
2 & $(1,1,0,0)=3\omega_1-\alpha_1    $ & $d_{12}\wedge d_{23}\wedge d_{14}\wedge d_{15}+
d_{12}\wedge d_{13}\wedge d_{24}\wedge d_{15}+d_{12}\wedge d_{13}\wedge d_{14}\wedge d_{25}
$  \\
\hline
3 & $(-1,2,0,0)=3\omega_1-2\alpha_1    $ & $d_{12}\wedge d_{23}\wedge d_{24}\wedge d_{15}
+d_{12}\wedge d_{23}\wedge d_{14}\wedge d_{25}+d_{12}\wedge d_{13}\wedge d_{24}\wedge d_{25}
$  \\
\hline
4 & $(2,-1,1,0) =3\omega_1-\alpha_1-\alpha_2   $ & $d_{13}\wedge d_{23}\wedge d_{14}\wedge d_{15}+
d_{12}\wedge d_{13}\wedge d_{34}\wedge d_{15}+
d_{12}\wedge d_{13}\wedge d_{14}\wedge d_{35}
$  \\
\hline
5 & $(-3,3,0,0)=3\omega_1-3\alpha_1   $ & $d_{12}\wedge d_{23}\wedge d_{24}\wedge d_{25}$  \\
\hline
6 & $(0,0,1,0) =3\omega_1-2\alpha_1-\alpha_2   $ &  $d_{13}\wedge d_{23}\wedge d_{24}\wedge d_{15}
+
d_{13}\wedge d_{23}\wedge d_{14}\wedge d_{25}
+d_{12}\wedge d_{23}\wedge d_{34}\wedge d_{15}
+d_{12}\wedge d_{13}\wedge d_{34}\wedge d_{25}
+d_{12}\wedge d_{23}\wedge d_{14}\wedge d_{35}
+d_{12}\wedge d_{13}\wedge d_{24}\wedge d_{35}
$ \\
\hline
7 & $(2,0,-1,1)=3\omega_1-\alpha_1-\alpha_2-\alpha_3  $ & $-d_{13}\wedge d_{14}\wedge d_{24}\wedge d_{15}+
d_{12}\wedge d_{13}\wedge d_{14}\wedge d_{45}+d_{12}\wedge d_{14}\wedge d_{34}\wedge d_{15}
$  \\
\hline
8 & $(-2,1,1,0)=3\omega_1-3\alpha_1-\alpha_2   $ & $d_{13}\wedge d_{23}\wedge d_{24}\wedge d_{25}+
d_{12}\wedge d_{23}\wedge d_{34}\wedge d_{25}+d_{12}\wedge d_{23}\wedge d_{24}\wedge d_{35}
$  \\
\hline
9 & $(1,-2,2,0)=3\omega_1-2\alpha_1-2\alpha_2    $ & $d_{13}\wedge d_{23}\wedge d_{34}\wedge d_{15}+
d_{13}\wedge d_{23}\wedge d_{14}\wedge d_{35}+d_{12}\wedge d_{13}\wedge d_{34}\wedge d_{35}
$  \\
\hline
10 & $(0,1,-1,1)=3\omega_1-2\alpha_1-\alpha_2-\alpha_3   $ & $-d_{23}\wedge d_{14}\wedge d_{24}\wedge d_{15}-
d_{13}\wedge d_{14}\wedge d_{24}\wedge d_{25}+d_{12}\wedge d_{23}\wedge d_{14}\wedge d_{45}
+d_{12}\wedge d_{13}\wedge d_{24}\wedge d_{45}+d_{12}\wedge d_{24}\wedge d_{34}\wedge d_{15}
+d_{12}\wedge d_{14}\wedge d_{34}\wedge d_{25}
$  \\
\hline
11 & $(2,0,0,-1)=3\omega_1-\alpha_1-\alpha_2-\alpha_3-\alpha_4    $ & $d_{13}\wedge d_{14}\wedge d_{15}\wedge d_{25}
+d_{12}\wedge d_{13}\wedge d_{15}\wedge d_{45}-d_{12}\wedge d_{14}\wedge d_{15}\wedge d_{35}
$  \\
\hline
12 & $(-1,-1,2,0)=3\omega_1-3\alpha_1-2\alpha_2    $ & $d_{13}\wedge d_{23}\wedge d_{34}\wedge d_{25}
+d_{13}\wedge d_{23}\wedge d_{24}\wedge d_{35}+d_{12}\wedge d_{23}\wedge d_{34}\wedge d_{35}
$  \\
\hline
13 & $(-2,2,-1,1) =3\omega_1-3\alpha_1-\alpha_2-\alpha_3   $ & $-d_{23}\wedge d_{14}\wedge d_{24}\wedge d_{25}
+d_{12}\wedge d_{24}\wedge d_{34}\wedge d_{25}+d_{12}\wedge d_{23}\wedge d_{24}\wedge d_{45}
$  \\
\hline
14 & $(1,-1,0,1)=3\omega_1-2\alpha_1-2\alpha_2-\alpha_3    $ & $-d_{23}\wedge d_{14}\wedge d_{34}\wedge d_{15}
-d_{13}\wedge d_{14}\wedge d_{24}\wedge d_{35}+d_{13}\wedge d_{23}\wedge d_{14}\wedge d_{45}
+d_{12}\wedge d_{13}\wedge d_{34}\wedge d_{45}+d_{13}\wedge d_{24}\wedge d_{34}\wedge d_{15}
+d_{12}\wedge d_{14}\wedge d_{34}\wedge d_{35}$  \\
\hline
15 & $(0,1,0,-1)=3\omega_1-2\alpha_1-\alpha_2-\alpha_3-\alpha_4    $ & $d_{23}\wedge d_{14}\wedge d_{15}\wedge d_{25}
+d_{13}\wedge d_{24}\wedge d_{15}\wedge d_{25}+d_{12}\wedge d_{23}\wedge d_{15}\wedge d_{45}
+d_{12}\wedge d_{13}\wedge d_{25}\wedge d_{45}-d_{12}\wedge d_{24}\wedge d_{15}\wedge d_{35}
-d_{12}\wedge d_{14}\wedge d_{25}\wedge d_{35}
$  \\
\hline
16 & $(-1,0,0,1)=3\omega_1-3\alpha_1-2\alpha_2-\alpha_3 $ & $-d_{23}\wedge d_{14}\wedge d_{34}\wedge d_{25}-
d_{23}\wedge d_{14}\wedge d_{24}\wedge d_{35}+d_{13}\wedge d_{24}\wedge d_{34}\wedge d_{25}
+d_{13}\wedge d_{23}\wedge d_{24}\wedge d_{45}+d_{12}\wedge d_{24}\wedge d_{34}\wedge d_{35}
+d_{12}\wedge d_{23}\wedge d_{34}\wedge d_{45}
$  \\
\hline
17 & $(-2,2,0,-1)=3\omega_1-3\alpha_1-\alpha_2-\alpha_3-\alpha_4    $ & $d_{23}\wedge d_{24}\wedge d_{15}\wedge d_{25}-
d_{12}\wedge d_{24}\wedge d_{25}\wedge d_{35}+d_{12}\wedge d_{23}\wedge d_{25}\wedge d_{45}
$  \\
\hline
18 & $(1,0,-2,2) =3\omega_1-2\alpha_1-2\alpha_2-2\alpha_3  $ & $d_{14}\wedge d_{24}\wedge d_{34}\wedge d_{15}
-d_{13}\wedge d_{14}\wedge d_{24}\wedge d_{45}+d_{12}\wedge d_{14}\wedge d_{34}\wedge d_{45}
$  \\
\hline
19 & $(1,-1,1,-1) =3\omega_1-2\alpha_1-2\alpha_2-\alpha_3-\alpha_4   $ & $d_{23}\wedge d_{14}\wedge d_{15}\wedge d_{35}
+d_{13}\wedge d_{34}\wedge d_{15}\wedge d_{25}
+d_{13}\wedge d_{23}\wedge d_{15}\wedge d_{45}+d_{12}\wedge d_{13}\wedge d_{35}\wedge d_{45}
-d_{12}\wedge d_{34}\wedge d_{15}\wedge d_{35}
-d_{13}\wedge d_{14}\wedge d_{25}\wedge d_{35}
$  \\
\hline
20 & $(0,-3,3,0)=3\omega_1-3\alpha_1-3\alpha_2   $ & $d_{13}\wedge d_{23}\wedge d_{34}\wedge d_{35}$  \\
\hline
21 & $(0,-2,1,1) =3\omega_1-3\alpha_1-3\alpha_2-\alpha_3   $ & $-d_{23}\wedge d_{14}\wedge d_{34}\wedge d_{35}
+d_{13}\wedge d_{23}\wedge d_{34}\wedge d_{45}+d_{13}\wedge d_{24}\wedge d_{34}\wedge d_{35}
$  \\
\hline
22 & $(-1,1,-2,2) =3\omega_1-3\alpha_1-2\alpha_2-2\alpha_3   $ & $d_{14}\wedge d_{24}\wedge d_{34}\wedge d_{25}-
d_{23}\wedge d_{14}\wedge d_{24}\wedge d_{45}+d_{12}\wedge d_{24}\wedge d_{34}\wedge d_{45}
$  \\
\hline
23 & $(1,0,-1,0) =3\omega_1-3\alpha_1-2\alpha_2-\alpha_3-\alpha_4   $ & $d_{23}\wedge d_{24}\wedge d_{15}\wedge d_{35}
+d_{23}\wedge d_{34}\wedge d_{15}\wedge d_{25}+d_{13}\wedge d_{23}\wedge d_{25}\wedge d_{45}
+d_{12}\wedge d_{23}\wedge d_{35}\wedge d_{45}-d_{12}\wedge d_{34}\wedge d_{25}\wedge d_{35}
-d_{13}\wedge d_{24}\wedge d_{25}\wedge d_{35}
$  \\
\hline
24 & $(0,-1,-1,2)=3\omega_1-2\alpha_1-2\alpha_2-2\alpha_3-\alpha_4    $ & $-d_{14}\wedge d_{24}\wedge d_{15}\wedge d_{35}
+d_{14}\wedge d_{34}\wedge d_{15}\wedge d_{25}+d_{13}\wedge d_{24}\wedge d_{15}\wedge d_{45}
+d_{12}\wedge d_{14}\wedge d_{35}\wedge d_{45}-d_{12}\wedge d_{34}\wedge d_{15}\wedge d_{45}
-d_{13}\wedge d_{14}\wedge d_{25}\wedge d_{45}
$  \\
\hline
25 & $(0,-2,2,-1)=3\omega_1-3\alpha_1-3\alpha_2-2\alpha_3    $ & $d_{14}\wedge d_{24}\wedge d_{34}\wedge d_{35}-
d_{23}\wedge d_{14}\wedge d_{34}\wedge d_{45}+d_{13}\wedge d_{24}\wedge d_{34}\wedge d_{45}
$  \\
\hline
26 & $(0,-2,2,-1) =3\omega_1-3\alpha_1-3\alpha_2-\alpha_3-\alpha_4   $ & $d_{23}\wedge d_{34}\wedge d_{15}\wedge d_{35}
+d_{13}\wedge d_{23}\wedge d_{35}\wedge d_{45}-d_{13}\wedge d_{34}\wedge d_{25}\wedge d_{35}
$  \\
\hline
27 & $(-1,1,-1,0)=3\omega_1-3\alpha_1-2\alpha_2-2\alpha_3-\alpha_4    $ & $d_{23}\wedge d_{24}\wedge d_{15}\wedge d_{45}
+d_{24}\wedge d_{34}\wedge d_{15}\wedge d_{25}-d_{23}\wedge d_{14}\wedge d_{25}\wedge d_{45}
+d_{12}\wedge d_{24}\wedge d_{35}\wedge d_{45}-d_{12}\wedge d_{34}\wedge d_{25}\wedge d_{45}
-d_{14}\wedge d_{24}\wedge d_{25}\wedge d_{35}
$  \\
\hline
28 & $(1,0,0,-2)=3\omega_1-2\alpha_1-2\alpha_2-2\alpha_3-2\alpha_4    $ & $d_{14}\wedge d_{15}\wedge d_{25}\wedge d_{35}
-d_{13}\wedge d_{15}\wedge d_{25}\wedge d_{45}+d_{12}\wedge d_{15}\wedge d_{35}\wedge d_{45}
$  \\
\hline
29 & $(0,-1,0,0)=3\omega_1-3\alpha_1-3\alpha_2-2\alpha_3-\alpha_4    $ & $d_{24}\wedge d_{34}\wedge d_{15}\wedge d_{35}
-d_{14}\wedge d_{34}\wedge d_{25}\wedge d_{35}+d_{23}\wedge d_{34}\wedge d_{15}\wedge d_{45}
-d_{23}\wedge d_{14}\wedge d_{35}\wedge d_{45}-d_{13}\wedge d_{34}\wedge d_{25}\wedge d_{45}
+d_{13}\wedge d_{24}\wedge d_{35}\wedge d_{45}
$  \\
\hline
30& $(-1,1,-1,0) =3\omega_1-3\alpha_1-2\alpha_2-2\alpha_3-\alpha_4   $ & $-d_{23}\wedge d_{15}\wedge d_{25}\wedge d_{45}
+d_{24}\wedge d_{15}\wedge d_{25}\wedge d_{35}+d_{12}\wedge d_{25}\wedge d_{35}\wedge d_{45}
$  \\
\hline
31 & $(0,0,-3,3)=3\omega_1-3\alpha_1-3\alpha_2-3\alpha_3   $ & $d_{14}\wedge d_{24}\wedge d_{34}\wedge d_{45}$  \\
\hline
32& $(0,0,-2,1) =3\omega_1-3\alpha_1-3\alpha_2-3\alpha_3-\alpha_4   $ & $d_{24}\wedge d_{34}\wedge d_{15}\wedge d_{45}-
d_{14}\wedge d_{34}\wedge d_{25}\wedge d_{45}+d_{14}\wedge d_{24}\wedge d_{35}\wedge d_{45}
$  \\
\hline
33& $(0,-1,1,-2)=3\omega_1-3\alpha_1-3\alpha_2-2\alpha_3-2\alpha_4    $ & $-d_{23}\wedge d_{15}\wedge d_{35}\wedge d_{45}
+d_{34}\wedge d_{15}\wedge d_{25}\wedge d_{35}+d_{13}\wedge d_{25}\wedge d_{35}\wedge d_{45}
$  \\
\hline
34 & $(0,0,-1,-1)=3\omega_1-3\alpha_1-3\alpha_2-3\alpha_3-2\alpha_4    $ & $-d_{24}\wedge d_{15}\wedge d_{35}\wedge d_{45}+
d_{14}\wedge d_{25}\wedge d_{35}\wedge d_{45}+d_{34}\wedge d_{15}\wedge d_{25}\wedge d_{45}
$  \\
\hline
35 & $(0,0,0,-3)=3\omega_1-3\alpha_1-3\alpha_2-3\alpha_3-3\alpha_4  $ & $d_{15}\wedge d_{25} \wedge d_{35}\wedge  d_{45}$  \\
\hline
\end{tabular}
\end{center}
\end{table}

\begin{table}
\caption{\bf Weights  for $\mbox{sl}_5$ module $V(\omega_1+\omega_2)$}
\begin{center}
\begin{tabular}{|c|p{6cm}|c|p{6cm}|}
\hline
$ i  $ &  $\overrightarrow{w}_{i}^{\omega_2}$   & $ i  $ &   $\overrightarrow{w}_{i}^{\omega_2}$ \\
\hline
1 & $(1,1,0,0)=\omega_1+\omega_2$ &  16 & $(1,0,-2,2)=\omega_1+\omega_2-\alpha_1-2\alpha_2-2\alpha_3$
 \\
\hline
2 & $(-1,2,0,0)=\omega_1+\omega_2-\alpha_1$ &
17 & $(1,-1,1,-1)=\omega_1+\omega_2-\alpha_1-2\alpha_2-\alpha_3-\alpha_4$
 \\
\hline
3 & $(2,-1,1,0)=\omega_1+\omega_2-\alpha_2$ &  18 & $(0,-2,1,1)=\omega_1+\omega_2-2\alpha_1-3\alpha_2-\alpha_3$
 \\
\hline
4 & $(0,0,1,0)=\omega_1+\omega_2-\alpha_1-\alpha_2$ &
19 & $(-1,1,-2,2)=\omega_1+\omega_2-2\alpha_1-2\alpha_2-2\alpha_3$
\\
\hline
5 & $(2,0,-1,1)=\omega_1+\omega_2-\alpha_2-\alpha_3$ &  20 & $(-1,0,1,-1)=\omega_1+\omega_2-2\alpha_1-2\alpha_2-\alpha_3-\alpha_4$  \\
\hline
6 & $(-2,1,1,0)=\omega_1+\omega_2-2\alpha_1-\alpha_2$ &  21 & $(1,0,-1,0)=\omega_1+\omega_2-\alpha_1-2\alpha_2-2\alpha_3-\alpha_4$  \\
\hline
7 & $(1,-2,2,0)=\omega_1+\omega_2-\alpha_1-2\alpha_2$ &  22 & $(0,-1,-1,2)=\omega_1+\omega_2-2\alpha_1-3\alpha_2-2\alpha_3$  \\
\hline
8 & $(0,1,-1,1)=\omega_1+\omega_2-\alpha_1-\alpha_2-\alpha_3$ &  23 & $(0,-2,2,-1)=\omega_1+\omega_2-2\alpha_1-3\alpha_2-\alpha_3-\alpha_4$  \\
\hline
9 & $(2,0,0,-1)=\omega_1+\omega_2-\alpha_2-\alpha_3-\alpha_4$ &  24 & $(-1,1,-1,0)=\omega_1+\omega_2-2\alpha_1-2\alpha_2-2\alpha_3-\alpha_4$  \\
\hline
10 & $(-1,-1,2,0)=\omega_1+\omega_2-2\alpha_1-2\alpha_2$ &  25 & $(1,0,0,-2)=\omega_1+\omega_2-\alpha_1-2\alpha_2-2\alpha_3-2\alpha_4$  \\
\hline
11 & $(-2,2,-1,1)=\omega_1+\omega_2-2\alpha_1-\alpha_2-\alpha_3$ &  26 & $(0,-1,0,0)=\omega_1+\omega_2-2\alpha_1-3\alpha_2-2\alpha_3-\alpha_4$  \\
\hline
12 & $(1,-1,0,1)=\omega_1+\omega_2-\alpha_1-2\alpha_2-\alpha_3$ &  27 & $(-1,1,0,-2)=\omega_1+\omega_2-2\alpha_1-2\alpha_2-2\alpha_3-2\alpha_4$  \\
\hline
13 & $(0,1,0,-1)=\omega_1+\omega_2-\alpha_1-\alpha_2-\alpha_3-\alpha_4$ &  28 & $(0,0,-2,1)=\omega_1+\omega_2-2\alpha_1-3\alpha_2-3\alpha_3-\alpha_4$  \\
\hline
14 & $(-1,0,0,1)=\omega_1+\omega_2-2\alpha_1-2\alpha_2-\alpha_3$ &  29 & $(0,-1,1,-2)=\omega_1+\omega_2-2\alpha_1-3\alpha_2-2\alpha_3-2\alpha_4$  \\
\hline
15 & $(-2,2,0,-1)=\omega_1+\omega_2-2\alpha_1-\alpha_2-\alpha_3-\alpha_4$ &  30 & $(0,0,-1,-1)=\omega_1+\omega_2-2\alpha_1-3\alpha_2-3\alpha_3-2\alpha_4$  \\
\hline
\end{tabular}
\end{center}
\end{table}

\end{appendix}

\end{document}